%% file: Maria.tex
\documentclass[11pt,a4paper,oneside,titlepage]{book}
\usepackage[english]{babel}
\usepackage[latin1]{inputenc}
\usepackage[T1]{fontenc}
\usepackage{indentfirst}
\usepackage{fancyhdr}
\usepackage{amsfonts}
\usepackage{pst-eps}
\usepackage{auto-pst-pdf}
\usepackage{amsmath}
\usepackage{amsthm}
\usepackage{amssymb}
\usepackage{stmaryrd} 
\usepackage{mathrsfs}
\usepackage{graphics,graphicx}
\usepackage{subfig}
\newtheorem{theorem}{Theorem}[section]
\newtheorem{definition}[theorem]{Definition}
\newtheorem{lemma}[theorem]{Lemma}
\newtheorem{rem}[theorem]{Remark}
\newtheorem{corollary}[theorem]{Corollary}
\newtheorem{notation}[theorem]{Notation}
\newtheorem{prop}[theorem]{Proposition}
\newtheorem{ex}[theorem]{Examples}
\begin{document}
\pagestyle{fancy} 
\renewcommand{\chaptermark}[1]{\markboth{#1}{}}
\renewcommand{\sectionmark}[1]{\markright{\thesection \ #1}}
\fancyhf{} 
\fancyhead[R]{\bfseries\thepage}
\fancyhead[L]{\sc \bfseries\rightmark}
\headheight= 14pt
\thispagestyle{empty}
\include{copertina}

\thispagestyle{empty}
\include{Dedica}
\frontmatter
\include{Abstract}

\include{Sommario}
\tableofcontents
\listoffigures
\mainmatter
\addcontentsline{toc}{chapter}{Introduction}
\input{Intro}

\input{Cap1}

\input{Cap2}

\input{Cap3}

\input{Concl}

\backmatter
\addcontentsline{toc}{chapter}{Bibliography}

\end{document}

%% file: copertina.tex
\begin{center}
{\Huge\sc Universit\`a della Calabria} \\
\vspace{0.3cm}
{ Dipartimento di Matematica} \\
\vspace{0.3cm}
{\large\bf Dottorato di Ricerca in Matematica ed Informatica} \\
\vspace{0.3cm}
{\sc xxiii ciclo}

\rule[0.1cm]{\textwidth}{0.1mm}
{ Settore Disciplinare MAT/05 -- ANALISI MATEMATICA}

\vspace{2 cm}
{\sc Tesi di Dottorato}

\vspace{1.0 cm} 
{\LARGE\bf\sc Uniform distribution of}\\
\vspace{0.2 cm}
{\LARGE\bf\sc sequences of points and partitions}\\
\vspace{2.5cm}
{Maria Infusino}

\vspace{4 cm}
\begin{tabular}{lcl}
{\bf Supervisore} & \hspace{2cm} & {\bf Coordinatore} \vspace{.25cm} \\
Prof. Aljo\v sa Vol\v ci\v c & & Prof. Nicola Leone \vspace{.25cm}
\end{tabular}

\vspace{2 cm}
\rule[0.1cm]{\textwidth}{0.1mm}

A.A. 2009 -- 2010
\end{center}

%% file: Dedica.tex
\null\vspace{\stretch{1}}
\begin{flushright}
\textbf{\emph{ To my family }}
\end{flushright}
\vspace{\stretch{2}}\null

%% file: Abstract.tex
\chapter*{Abstract}
The interest for uniformly distributed (u.d.) sequences of points, in particular for low discrepancy sequences, arises from various applications, especially in the field of numerical integration. The basic idea in numerical integration is trying to approximate the integral of a function $f$ by a weighted average of the function evaluated at a set of points $\{x_1,\ldots, x_N\}$
$$\int_{I^d}f(x)dx\approx\frac{1}{N}\sum_{i=1}^Nw_if(x_i),$$
where $I^d$ is the $d-$dimensional unit hypercube, the $x_i$'s are $N$ points in $I^d$ and $w_i>0$ are weights such that $\sum\limits_{i=1}^Nw_i=N$. In some cases it is assumed $w_i=1$ for every $1\leq i\leq N$, as for instance in the classical Monte Carlo method where the points $x_1,\ldots, x_N$ are picked from a sequence of random or pseudorandom elements in $I^d$. Another possibility is to use deterministic sequences with given distribution properties. This procedure is known as Quasi-Monte Carlo method and it is more advantageous than many other approximation techniques. In fact, as the Koksma-Hlawka inequality states, the quality of the approximation provided by the Quasi-Monte Carlo method is linked directly to the discrepancy of the $x_i$'s. The better the nodes are distributed in $I^d$, the faster the approximation is expected. Hence, a good choice for the integration points is the initial segment of a sequence with small discrepancy.

In this context the construction of u.d.\! sequences with low discrepancy in various spaces is of crucial importance. The objectives of this thesis are related to this main topic of uniform distribution theory and can be summarized as follows:
\begin{description}
\item[(A)] The research of explicit techniques for introducing new classes of u.d.\! sequences of points and of partitions on $[0,1]$ and also on fractal sets,
\item[(B)] A quantitative analysis of the distribution behaviour of a class of generalized Kakutani's sequences on $[0,1]$ through the study of their discrepancy.
\end{description}

To achieve these purposes, a fundamental role is played by the concept of u.d.\! sequences of partitions. In fact when we deal with fractals, and in particular with fractals generated by an Iterated Function System (IFS), partitions turn out to be a convenient tool for introducing a uniform distribution theory. In this thesis we extend to certain fractals the notion of u.d.\! sequences of partitions, introduced by Kakutani in 1976 for the unit interval and we employ it to construct van der Corput type sequences on a whole class of IFS fractals. More precisely in Chapter 2, where we develop the objective $\textbf{(A)}$, we present a general algorithm to produce u.d.\! sequences of partitions and of points on the class of fractals generated by a system of similarities on $\mathbb{R}^d$ having the same ratio and verifying the open set condition. We also provide an estimate for the elementary discrepancy of these sequences.

Generalized Kakutani's sequences of partitions of $[0,1]$ are extremely useful in the extension of these results to a wider class of fractals obtained by eliminating the restriction that all the similarities defining the fractal have the same ratio. According to a remark by Mandelbrot, which allows to see $[0,1]$ as the attractor of an IFS, the simplest setting for this problem is the unit interval.
Perfectly fitting our problem is a recent generalization of Kakutani's splitting procedure on $[0,1]$, namely the technique of $\rho-$refinements. Consequently, in Chapter 3 we deal with objective \textbf{(B)} and focus on deriving bounds for the discrepancy of the sequences generated by this technique.

Our approach is based on a tree representation of any sequence of partitions constructed by successive $\rho-$refinements, which is exactly the parsing tree generated by Khodak's coding algorithm. This correspondence allows to give bounds of the discrepancy for all the sequences generated by successive $\rho-$refinements, when $\rho$ is a partition of $[0,1]$ consisting of $m$ subintervals of lenghts $p_1,\ldots,p_m$ such that $\log\left(\frac{1}{p_1}\right),\ldots,\log\left(\frac{1}{p_m}\right)$ are rationally related. This result applies also to a countable family of classical Kakutani's sequences and provides estimates of their discrepancy, not known in the existing literature. 
Moreover, we are also able to cover several situations in the irrational case, which means that at least one of the fractions $\frac{\log p_i}{\log p_j}$ is irrational. More precisely, we discuss some instances of the irrational case when the initial probabilities are $p$ and $q=1-p$. In this case we obtain weaker upper bounds for the discrepancy, since they depend heavily on Diophantine approximation properties of the ratio $\frac{\log p}{\log q}$. Finally, we prove bounds for the elementary discrepancy of the sequences of partitions constructed through an adaptation of the $\rho-$refinements method to the new class of fractals.

%% file: Sommario.tex
\chapter*{Sommario}
L'interesse per le successioni di punti uniformemente distribuite (u.d.) emerge da svariate applicazioni specialmente nell'ambito dell'integrazione numerica. Un approccio tipico di questa disciplina è l'approssimazione dell'integrale di una funzione $f$ con la media pesata dei valori assunti dalla funzione in un insieme di punti $\{x_1,\ldots, x_N\}$
$$\int_{I^d}f(x)dx\approx\frac{1}{N}\sum_{i=1}^Nw_if(x_i),$$
dove $I^d$ è l'ipercubo unitario $d-$dimensionale, gli $x_i$ sono $N$ elementi di $I^d$ e i pesi $w_i>0$ sono tali che $\sum\limits_{i=1}^Nw_i=N$. In alcuni casi si assume che $w_i=1$ per ogni $1\leq i\leq N$, come ad esempio nel metodo classico di Monte Carlo in cui i punti $x_1,\ldots, x_N$ sono selezionati da una successione casuale o pseudo-casuale di elementi in $I^d$. Un'altra possibilità è effettuare la scelta degli $x_i$ all'interno di successioni deterministiche con proprietà di distribuzione fissate. Questa procedura è nota come metodo di Quasi-Monte Carlo ed è più vantaggiosa di molte altre tecniche d'approssimazione numerica. Infatti, la disuguaglianza di Koksma-Hlawka stabilisce che la qualità dell'approssimazione fornita dal metodo di Quasi-Monte Carlo è strettamente legata alla discrepanza degli $x_i$. Pertanto, risulta conveniente scegliere come insieme dei punti di integrazione il segmento iniziale di una successione a bassa discrepanza. 

La ricerca di successioni di punti u.d.\!\! con bassa discrepanza è dunque di importanza cruciale in ambito applicativo. Gli obiettivi di questo lavoro si collocano all'interno di questo filone di ricerca e interessano due tematiche fondamentali:
\begin{description}
\item[(A)] la ricerca di tecniche esplicite che consentano di costruire successioni u.d.\! di punti e di partizioni su $[0,1]$ e su insiemi frattali,
\item[(B)] l'analisi del comportamento asintotico della discrepanza di una classe di successioni di partizioni di Kakutani generalizzate.
\end{description}

Nei risultati proposti uno strumento essenziale è il concetto di successione di partizioni u.d.\!. Infatti quando si lavora con i frattali, ed in particolare con frattali generati da un Sistema di Funzioni Iterate (IFS), le partizioni risultano essere più convenienti delle successioni di punti in relazione alla teoria della distribuzione uniforme. Pertanto abbiamo esteso ai frattali la definizione di successione di partizioni u.d., introdotta da Kakutani nel 1976 per partizioni di $[0,1]$, ed abbiamo sfruttato questo concetto per costruire successioni di tipo van der Corput su un'intera classe di frattali IFS. Più precisamente nel Capitolo 2, in cui viene affrontata la tematica \textbf{(A)}, presentiamo un algoritmo per generare successioni u.d.\! di punti e di partizioni sui frattali individuati da un numero finito di similitudini su $\mathbb{R}^d$, aventi tutte lo stesso rapporto di similitudine e che soddifano la condizione dell'insieme aperto. Inoltre abbiamo ricavato una stima della discrepanza elementare delle successioni prodotte.

La seconda problematica studiata è l'estensione dei risultati ottenuti a una classe più ampia di frattali, eliminando la restrizione che le similitudini dell'IFS abbiano tutte lo stesso rapporto. Secondo un'osservazione dovuta a Mandelbrot, che consente di vedere $[0,1]$ come attrattore di infiniti IFS, l'ambientazione più semplice per tale problema è proprio l'intervallo unitario. Una tecnica che si adatta perfettamente alle caratteristiche della nuova classe di attrattori è una recente generalizzazione della procedura di Kakutani: la tecnica dei $\rho$-raffinamenti. Pertanto, nel Capitolo~3 affrontiamo la tematica~\textbf{(B)} con l'obiettivo di determinare stime della discrepanza delle successioni di partizioni di $[0,1]$ prodotte tramite tale tecnica.

L'approccio che usiamo è basato su una rappresentazione ad albero di questa classe di successioni che produce lo stesso albero costruito secondo l'algoritmo di Khodak. Questa corrispondenza consente di ricavare stime della discrepanza delle successioni generate dai successivi $\rho-$raffinamenti dell'intervallo unitario, quando $\rho$ è una partizione costituita da $m$ intervalli di lunghezza $p_1,\ldots,p_m$ tali che $\log\left(\frac{1}{p_1}\right),\ldots$\\$\ldots,\log\left(\frac{1}{p_m}\right)$ siano razionalmente correlati. Questo caso include una classe numerabile di successioni di Kakutani classiche, per le quali otteniamo stime della discrepanza ancora non presenti in letteratura. Per quanto concerne il caso irrazionale, cioè quando almeno uno dei rapporti $\frac{\log p_i}{\log p_j}$ non è razionale, sono state osservate diverse complicazioni. In questo lavoro analizziamo la situazione in cui $\rho$ è costituita da due intervalli di lunghezza $p$ e $q=1-p$. Tuttavia, le stime della discrepanza ottenute in questo sottocaso sono più deboli, in quanto dipendono fortemente dalle proprietà di approssimazione diofantea del rapporto $\frac{\log p}{\log q}$. Infine, introduciamo alcuni risultati sulla discrepanza elementare delle successioni di partizioni costruite tramite un adattamento del metodo dei $\rho-$raffinamenti alla nuova classe di frattali.

%% file: Intro.tex
\chapter*{Introduction}
\markboth{INTRODUCTION}{INTRODUCTION}
The theory of uniform distribution was developed extensively within and among several mathematical disciplines and numerous applications. In fact, the main root of this theory is number theory and diophantine approximation, but there are strong connections to various fields of mathematics such as measure theory, probability theory, harmonic analysis, summability theory, discrete mathematics and numerical analysis.
 
The central goals of this theory are the assessment of uniform distribution and the construction of uniformly distributed (u.d.)\! sequences in various mathematical spaces. The objectives of this thesis are related to these main topics.
In particular, the aim of this work is to introduce new classes of u.d.\! sequences of points and of partitions on $[0,1]$ and also on fractal sets. Moreover, we intend to present a quantitative analysis of the distribution behaviour of the new sequences produced studying their discrepancy.
\\

The problem of finding explicit methods for constructing u.d.\! sequences was originally investigated in the setting of sequences of points. In fact, the starting point of the development of the theory was just the study of u.d.\! sequences of points on the unit interval. The result which marked the beginning of the theory was the discovery that the fractional parts of the multiples of an irrational number are u.d.\! in the unit interval or, equivalently, on the unit circle. This was a refinement of an approximation theorem due to Kronecker who had already proved the density of this special sequence in the unit interval. So, at the beginning of the last century, many authors independently proposed the theorem about uniform distribution of Kronecker's sequence such as Bohl \cite{B}, Sierpi\' nski \cite{S} and Weyl \cite{W1}. The latter was the first to estabilish a systematic treatment of uniform distribution theory in his famous paper \cite{W3}, where the formal definition of u.d.\! sequences of points in $[0,1]$ was given for the first time. Moreover, in that paper the theory of u.d.\! sequences of points was generalized to the higher-dimensional unit cube. 

The uniform distribution of a sequence of points means that the empirical distribution of the sequence is asymptotically equal to the uniform distribution. Therefore in the twenties and thirties several authors began to study u.d.\! sequences of points from a quantitative point of view introducing the discrepancy \cite{Berg, VdCPis, W3}. This quantity is the classical measure of the deviation of a sequence from the ideal uniform distribution. Consequently, having a precise estimate of the discrepancy is very useful for applications but it is not a trivial problem. Proving general lower bounds for the discrepancy is a subject still having open questions nowadays.
\\

The interest for u.d.\! sequences of points, in particular for low discrepancy sequences, arises from various applications in areas like numerical integration, random number generation, stochastic simulation and approximation theory. Indeed, numerical integration was one of the first applications of uniform distribution theory \cite{Hlaw3}. The basic problem considered by numerical integration is to compute an approximate solution to a definite integral. The classical quadrature formulae are less and less efficient the higher the dimension is. To overcome this problem, a typical approach is trying to approximate the integral of a function $f$ by a weighted average of the function evaluated at a set of points $\{x_1,\ldots, x_N\}$
$$\int_{I^d}f(x)dx\approx\frac{1}{N}\sum_{i=1}^Nw_if(x_i),$$
where $I^d$ is the $d-$dimensional unit hypercube, the $x_i$'s are $N$ points in $I^d$ and $w_i>0$ are weights such that $\sum\limits_{i=1}^Nw_i=N$. In some cases it is assumed $w_i=1$ for every $1\leq i\leq N$, as for instance in the classical Monte Carlo method where the points $x_1,\ldots, x_N$ are picked from a sequence of random or pseudorandom elements in $I^d$. The advantage of the Monte Carlo method is that it is less sensitive to the increase of the dimension.

Another possibility is to use deterministic sequences with given distribution properties for the choice of the $x_i$'s. This procedure is known as Quasi-Monte Carlo method and it is more advantageous than many other approximation techniques. In fact, the Koksma-Hlawka inequality (\ref{KHIneq}) shows that the error of such a method can be bounded by the product of a term only depending on the discrepancy of $\{x_1,\ldots, x_N\}$ and one only depending on the function. Therefore it is convenient to choose the initial segment of a low discrepancy sequence as the set of integration points in the Quasi-Monte Carlo method. These are sequences with a discrepancy of order $\frac{(\log N)^d}{N}$, where $d$ is the dimension of the space in which we take the sequence. Hence, by using low discrepancy sequences, the Quasi-Monte Carlo method has a faster rate of convergence than a corresponding Monte Carlo method, since in the latter case the point sets do not have necessarily minimal discrepancy. Infact, it behaves, in average, as $\frac{1}{\sqrt{N}}$. Indeed, the Monte Carlo method yields only a probabilistic bound on the integration error. Neverthless, both Monte Carlo and Quasi-Monte Carlo methods offer the advantage to add further points without recalculating the values of the function in the previous points and this is a big step forward compared to classical methods. Quasi-Monte Carlo methods have an important role in financial and actuary mathematics, where high-dimensional integrals occur. During the last twenty years all these applications have been a rapidly growing area of research \cite{NiedBook, Hel-Larch}. 
\\

One of the best known techiniques for generating low discrepancy sequences of points in the unit interval was introduced by van der Corput in 1935 (see \cite{VdC}). Successively, van der Corput's procedure was extended to the higher-dimensional case by Halton \cite{Hal}. Moreover, a generalization of van der Corput sequences is due to Faure who introduced the permuted or generalized van der Corput sequences. They are also very interesting because there exist formulae for the discrepancy of these sequences which show their good asymptotic behaviour \cite{F1,F2,F3}. 

The study of van der Corput type sequences has not been limited to the classical setting of the unit interval in one dimension or the unit hypercube in higher dimensions, but interesting extensions have been made to more abstract spaces such as fractals. In fact, the theory of uniform distribution with respect to a given measure has been generalized in several ways: sequences of points in compact and locally compact spaces \cite{KN, Nied, Helm}, sequences of probability measures on a separable compact space \cite{Sch2}, in particular sequences of discrete measures associated to partitions of a compact interval \cite{Ka} and to partitions of a separable metric space \cite{CV}. In the following we use the basic definitions of uniform distribution theory in compact Hausdorff spaces and in a particular class of fractal compact sets. 
\\

Fractals are involved in several applications because they are a powerful tool to describe effectively a variety of phenomena in a large number of fields. To exploit Quasi-Monte Carlo methods on these sets it is essential to study discrepancy bounds for sequences of points on fractals. One of the earlier papers devoted to uniform distribution on fractals is \cite{Grab-Tichy}, where this theory is developed on the Sierpi\' nski gasket. In this paper the notion of discrepancy on fractals has been introduced for the first time. The authors define several concepts of discrepancy for sequences of points on the Sierpi\' nski gasket by choosing different kinds of partitions on this fractal. Successively, these notions were generalized also to other fractals, such as the $d-$dimensional Sierpi\' nski carpet in \cite{CT, CT2}. In particular, in \cite{CT2} a van der Corput type construction is considered to generate u.d.\! sequences of points on the $d-$dimensional Sierpi\' nski carpet and the exact order of convergence of various notions of discrepancy is determined for these sequences.

In this work we get a more general result by constructing van der Corput type sequences on a whole class of fractals generated by an Iterated Function System (IFS). More precisely, we are going to study fractals defined by a system of similarities on $\mathbb{R}^d$ having the same ratio and verifying a natural separation condition of their components, namely the Open Set Condition (OSC). This class includes the most popular fractals, but also the unit interval $[0,1]$ which can be seen as the attractor of infinitely many different IFS. Starting from this remark, which goes back to Mandelbrot \cite{M2}, we present an alternative construction of the classical van der Corput sequences of points on $[0,1]$. By imitating this approach, we introduce an explicit procedure to define u.d.\! sequences of points on our special class of fractals (see Subsection \ref{subsec: Algo}). So we call these sequences of \emph{van der Corput type}, just to emphasize the particular order given to the points by our algorithm. It is important to underline that as probability on a fractal $F$ of our class we take the normalized $s$-dimensional Hausdorff measure, where $s$ is the Hausdorff dimension of $F$. This is the most natural choice for a probability measure on this kind of fractals, also because the OSC guarantees the existence of an easy formula for evaluating the Hausdorff dimension of these fractals (see Theorem \ref{MoranThm}). A crucial role in the proof of the uniform distribution of the sequences constructed is played by the elementary sets, i.e.\! the family of all sets generated by applying our algorithm to the whole fractal $F$. In this way our technique produces also u.d.\! sequences of partitions of the fractals belonging to the considered class. 
\\

The concept of u.d.\! sequence of partitions on fractals is just one of the most important aspects of this thesis. When we deal with fractals, and in particular with IFS fractals, partitions turn out to be a more convenient tool in relation to the uniform distribution theory. Consequently we extend the notion of u.d.\! sequences of partitions, introduced by Kakutani in 1976 for the unit interval in \cite{Ka}, to our class of fractals.

The construction ideated by Kakutani, called Kakutani's splitting procedure, allows to construct a whole class of u.d.\! sequences of partitions of $[0,1]$ and it is based on the concept of $\alpha-$refinement of a partition. For a fixed $\alpha\in]0,1[$, the $\alpha-$refinement of a partition $\pi$ is obtained by splitting all the intervals of $\pi$ having maximal lenght in two parts, proportional to $\alpha$ and $1-\alpha$ respectively. Kakutani proved that the sequence of partitions generated through successive $\alpha-$refinements of the trivial partition $\omega=\{[0,1]\}$ is u.d.. This result received a considerable attention in the late seventies, when other authors provided different proofs of Kakutani's theorem \cite{AF} and of its stochastic versions, in which the intervals of maximal lenght are splitted according to certain probability distributions \cite{VZ, LO1, LO2, BD, PvZ}.
Recently different generalizations of Kakutani's technique have been introduced. A result in this direction is the extension of Kakutani's splitting procedure to the multidimensional case with a construction which is intrinsically higher-dimensional~\cite{Vol-Car}. Moreover, in a recent paper of Vol\v{c}i\v{c}, Kakutani's technique is extended also in the one dimensional case introducing the concept of $\rho-$refinement of a partition, which generalizes Kakutani's $\alpha-$refinement. Actually, the $\rho-$refinement of a partition $\pi$ is obtained by splitting the longest intervals of $\pi$ into a finite number of parts homothetically to a given finite partition $\rho$ of $[0,1]$. The author has proved that the technique of successive $\rho-$refinements allows to construct new families of u.d.\! sequences of partitions of $[0,1]$ in \cite{Vol}. The last paper also investigates the connections of the theory of u.d.\! sequences of partitions to the well-estabilished theory of u.d.\! sequences of points, showing how it is possible to associate u.d.\! sequences of points to any u.d.\! sequence of partitions.
\\

Generalized Kakutani's sequences on $[0,1]$ are a fundamental tool in the extension of the results obtained on our class of fractals. The first attempt of enlarging the class of fractals considered in our previous analysis consists in eliminating the restriction that all the similarities defining the fractal have the same ratio.  

The procedure of successive $\rho$-refinements fits perfectly to the problem of generating u.d.\! sequences of partitions on this new class of fractals. Let $\psi=\{\psi_1,\ldots,\psi_m\}$ be a system of $m$ similarities on $\mathbb{R}^d$ having ratio $c_1,\ldots,c_m\in\ ]0,1[$ respectively and such that they verify the OSC. Let $F$ be the attractor of $\psi$ and let $s$ be its Hausdorff dimension. Applying successively the $m$ similarities to the fractal $F$, we get a first partition consisting of $m$ subsets of $F$ each of probability $p_i=c_i^s$ (where for probability we again mean the normalized $s-$dimensional Hausdorff measure). At the second step we choose the susbsets with the highest probability and we apply to each of them the $m$ similarities in the same order, and so on. Iterating this procedure, which exploits the same basic idea of $\rho-$refinements, we obtain a sequence of partitions of $F$. Now the problem is the assessment of the uniform distribution of these sequences and the estimation of their discrepancy.

According to the Mandelbrot's remark the simplest setting for this problem is the unit interval. In fact, if we consider $[0,1]$ as the attractor of $m$ similarities $\varphi_1,\ldots,\varphi_m$ having different ratios and satisfying the OSC and we apply the procedure described above, then we get exactly the sequence of $\rho-$refinements $(\rho^n\omega)$, where $\rho=\{\varphi_1([0,1]),\ldots,\varphi_m([0,1])\}$ and $\omega=\{[0,1]\}$.
\\

In the second part of this work we focus on deriving bounds for the discrepancy of the generalized Kakutani's sequences of partitions of $[0,1]$ generated through the techinique of successive $\rho-$refinements. The problem of estimating the asymptotic behaviour of the discrepancy of these sequences has been posed for the first time in \cite{Vol}. At the moment the only known discrepancy bounds for a class of such sequences have been given by Carbone in \cite{Car}. In this paper the author considered the so-called $LS$-sequences which are generated by successive $\rho-$refinements where $\rho$ is a partition with $L$ subintervals of $[0,1]$ of length $\alpha$ and $S$ subintervals of length $\alpha^2$ (where $\alpha$ is given by the equation $L\alpha + S\alpha^2 = 1$).

To study this problem in more generality we use a correspondence between the procedure of successive $\rho-$refinements and Khodak's algorithm \cite{Kho}. This new approach is based on a parsing tree related to Khodak's coding algorithm, which represents the successive $\rho$-refinements. We introduce improvements of the results obtained in \cite{Drmota} to provide significative bounds of the discrepancy for all the sequences generated by successive $\rho-$refinements, when $\rho$ is a partition of $[0,1]$ consisting of $m$ subintervals of lenghts $p_1,\ldots,p_m$ such that $\log\left(\frac{1}{p_1}\right),\ldots,\log\left(\frac{1}{p_m}\right)$ are rationally related. This result applies also to a countable family of classical Kakutani's sequences and provides, for the first time after thirty years, quantitative estimates of their discrepancy. Moreover, the class of generalized Kakutani's sequences belonging to this rational case also includes the $LS-$sequences.

In the following we are also able to cover several situations in the irrational case, which means that at least one of the fractions $\frac{\log p_i}{\log p_j}$ is irrational. This case is much more involved than the rational one. In this work we discuss some instances of the irrational case when the initial probabilities are two, namely $p$ and $q=1-p$. The upper bounds for the discrepancy that we obtain in this subcase are weaker, since they depend heavily on Diophantine approximation properties of the ratio $\frac{\log p}{\log q}$. Furthermore, if the initial partition is composed of more than two intervals, then the analysis of the behaviour of the discrepancy is even more complicated, as evident by comparing with \cite{Flajolet-Vallee}.

The approach applied for achieving these bounds of the discrepancy of generalized Kakutani's sequences on $[0,1]$ can be also used for the sequences of partitions constructed on fractals defined by similarities which do not have the same ratio and satisfing the OSC. In fact, we have described above an analogue of the method of successive $\rho-$refinements which allows to produce sequences of partitions on this new class of fractals. We actually introduce a new correspondence between nodes of the tree associated to Khodak's algorithm and the subsets belonging to the partitions generated on the fractal. Consequently, with a technique similar to the one used on $[0,1]$ we prove bounds for the elementary discrepancy of these sequences of partitions, too.
\\

Let us give a brief outline of the thesis.\\
\indent\textbf{Chapter 1} provides the basic background knowledge on the areas of uniform distribution theory that are investigated in this thesis. The first part of the chapter deals with the classical part of the theory. Basic definitions and properties of u.d.\! sequences of points on the unit interval are introduced and specific examples of u.d.\! sequences of points are described throughout. Then a whole section is devoted to the more recent theory of u.d.\! sequences of partitions, which plays an essential role in this work. Some extensions of uniform distribution theory are also touched on in this chapter, such as the theory in the unit hypercube and the theory in Hausdorff compact spaces. 

\textbf{Chapter 2} regards the uniform distribution on a special class of fractals. More precisely, we are concerned with fractals generated by an iterated function system of similarities having the same ratio and satisfying the open set condition. We propose an algorithm for generating u.d.\! sequences of partitions and of points on this class of fractals. Furthermore, in the last part of this chapter we study the order of convergence of the elementary discrepancy of the van der Corput type sequences constructed on these fractals. The results presented in this chapter have been first published in \cite{Inf-Vol}.

In \textbf{Chapter 3} we extend the results given in the second chapter to a wider class of fractals by using a new approach, which allows to derive bounds for the discrepancy of a class of generalized Kakutani's sequences of partitions of $[0,1]$, constructed through successive $\rho-$refinements. We present the recent technique of $\rho-$refinements and the generalization of Kakutani's theorem to the class of sequences of partitions generated by this procedure. Then, we analyze the behaviour of the discrepancy of these sequences from a new point of view. The crucial idea is a tree representation of any sequence of partitions constructed by successive $\rho-$refinements, which is precisely the parsing tree generated by Khodak's coding algorithm. The correspondence between the two techniques allows not only to give optimal upper bounds in the so-called rational case on $[0,1]$ but also to extend the results obtained in the second chapter to a wider class of fractals. Moreover, we study the irrational case which is more involved than the rational one. Finally, we give some examples and applications of the results achieved so far. The new contributions presented in this chapter are collected in \cite{Drm-Inf}.

The thesis concludes by reviewing, in \textbf{Chapter 4}, the main results
we have obtained and indicating open problems and directions of future research.

%% file: Cap1.tex
\chapter{Preliminary topics}
This chapter is meant to give a short overview of known results about uniform distribution theory not only in the classical setting of $[0,1]$ but also in more general spaces. First we intend to mention some necessary definitions and basic results concerning u.d.\! sequences of points in $[0,1]$. Then we will introduce the more recent theory of u.d.\! sequences of partitions which is fundamental in the development of this work. Finally, we will point out the main aspects of uniform distribution theory on the unit hypercube and on compact spaces.

\section{Uniformly distributed sequences of points in~$[0,1]$}
In this section we develop the classical part of uniform distribution theory. The standard references for this topic are \cite{KN} and \cite{DT}. We start introducing the basic concepts related to u.d.\! sequences of points and then we proceed to consider the quantitative aspect of the theory. Moreover, a whole subsection is devoted to a special class of sequences with certain advantageous distribution properties, namely the van der Corput sequences.
\subsection{Definitions and basic properties}
First of all, let us state the main definition of the theory.

\begin{definition}\label{def-ud}\ \\
A sequence $(x_n)$ of points in $[0,1]$ is said to be \emph{uniformly distributed (u.d.)} if for any real number $a$ such that $0<a\leq 1$ we have
\begin{equation}\label{succ-ud}
\lim_{N\rightarrow{\infty}}\frac{1}{N}\sum_{n=1}^{N}{\chi_{[0,a[}{(x_n)}}=a
\end{equation}
where $\chi_{[0,a[}$ is the characteristic function of the interval $[0,a[$.
\end{definition}

Let us introduce some concepts which are very useful to characterize u.d.\! sequences of points.
\begin{definition} \ \\
A class $\mathcal{F}$ of Riemann-integrable functions on $[0,1]$ is said to be \emph{determining} for the uniform distribution of sequences of points, if for any sequence $(x_n)$ in $[0,1]$ the validity of the relation
\begin{equation}\label{famDet}
\lim_{N\rightarrow{\infty}}\frac{1}{N}\sum_{n=1}^{N}{f(x_n)}=\int_{0}^{1}f(x)\ dx
\end{equation}
for all $f\in\mathcal{F}$ already implies that $(x_n)$ is u.d.. In particular, a system of subsets of $[0,1]$ such that the family of their characteristic functions is determining is called \emph{discrepancy system}.
\end{definition}

Hence, we can restate the Definition \ref{def-ud} saying that the family of all characteristic functions $\chi_{[0,a[}$ for $0<a\leq 1$ is determining or that the system of all sets $[0,a[$ for $0<a\leq 1$ is a discrepancy system.

An important determining class is the family of all continuous (real or complex-valued) functions on $[0,1]$. This result is due to Weyl and it is very useful to extend the theory to more general spaces \cite{W2, W3}.
\begin{theorem}[Weyl's Theorem]\label{continueDet}\ \\
A sequence $(x_n)$ of points in $[0,1]$ is u.d.\! if and only if for any real-valued continuous function $f$ defined on $[0,1]$ the equation (\ref{famDet}) holds.
\end{theorem}
\proof\ \\ Let $(x_n)$ be u.d.\! and let $f$ be a step function
\begin{equation}\label{fSemplice}
f(x)=\sum_{i=0}^{k-1}{c_i\chi_{[a_i,a_{i+1}[}{(x)}}
\end{equation}  
where $0=a_0<a_1<\ldots<a_k=1$ and $c_i\in\mathbb{R}$ for $i=0,\ldots,k-1$. Then it follows from~(\ref{succ-ud}) and~(\ref{fSemplice}) that
\begin{eqnarray*}\label{dis}
\lim_{N\rightarrow{\infty}}\frac{1}{N}\sum_{n=1}^{N}{f(x_n)}&=&\lim_{N\rightarrow{\infty}}
\frac{1}{N}\sum_{n=1}^{N}\sum_{i=0}^{k-1}{c_i\chi_{[a_i,a_{i+1}[}{(x_n)}} 
\\&=&\sum_{i=0}^{k-1}c_i\lim_{N\rightarrow{\infty}}\frac{1}{N} \left( \sum_{n=1}^{N}\chi_{[0,a_{i+1}[}{(x_n)}
-\sum_{n=1}^{N}\chi_{[0,a_{i}[}{(x_n)}\right)
\\&=&\sum_{i=0}^{k-1}c_i \left( \int_{0}^{1}\chi_{[0,a_{i+1}[}(x)\ dx
-\int_{0}^{1}\chi_{[0,a_{i}[}(x)\ dx\right)
\\&=&\sum_{i=0}^{k-1}\int_{0}^{1}{c_i\chi_{[a_i,a_{i+1}[}{(x)}}\ dx
\\&=&\int_{0}^{1}{f(x)}\ dx.
\end{eqnarray*}
Now, assume that $f$ is a real-valued function defined on $[0,1]$. Fixed $\varepsilon>0$, by the definition of the Riemann integral, there exist two step functions $f_1$ and $f_2$ such that $$f_1(x)\leq f(x)\leq f_2(x) \ ,\ \forall x\in[0,1]$$and$$ \int_{0}^{1}{\left(f_2(x)-f_1(x)\right)}\ dx\leq\varepsilon.$$
Then we have the following chain of inequalities
\begin{eqnarray*}
\int_{0}^{1}{f(x)}\ dx-\varepsilon&\leq&\int_{0}^{1}{f_2(x)}\ dx-\varepsilon\leq\int_{0}^{1}{f_1(x)}\ dx=\lim_{N\rightarrow{\infty}}\frac{1}{N}\sum_{n=1}^{N}{f_1(x_n)}\nonumber \\&\leq&\liminf_{N\rightarrow{\infty}}\frac{1}{N}\sum_{n=1}^{N}{f(x_n)} \leq\limsup_{N\rightarrow{\infty}}\frac{1}{N}\sum_{n=1}^{N}{f(x_n)}\nonumber \\&\leq&\lim_{N\rightarrow{\infty}}\frac{1}{N}\sum_{n=1}^{N}{f_2(x_n)} =\int_{0}^{1}{f_2(x)}\ dx\nonumber \\&\leq&\int_{0}^{1}{f_1(x)}\ dx+\varepsilon\leq\int_{0}^{1}{f(x)}\ dx+\varepsilon.
\end{eqnarray*}
So the relation (\ref{famDet}) holds for all continuous functions on $[0,1]$.\\
\indent Conversely, let $(x_n)$ be a sequence of points in $[0,1]$ such that the (\ref{famDet}) holds for every real-valued continuous function $f$ defined on $[0,1]$. Let $a\in ]0,1[$, then  for any $\varepsilon>0$ there exist two continuous functions $g_1$ and $g_2$ such that $$g_1(x)\leq\chi_{[0,a[}{(x)}\leq g_2(x)\ ,\ \forall x\in [0,1]$$ and $$ \int_{0}^{1}{\left(g_2(x)-g_1(x)\right)}\ dx\leq\varepsilon.$$ Then we have
\begin{eqnarray*}
a-\varepsilon&\leq&\int_{0}^{1}{g_2(x)}\ dx-\varepsilon\leq\int_{0}^{1}{g_1(x)}\ dx=\lim_{N\rightarrow{\infty}}\frac{1}{N}\sum_{n=1}^{N}{g_1(x_n)}\nonumber \\&\leq&
\liminf_{N\rightarrow{\infty}}\frac{1}{N}\sum_{n=1}^{N}{\chi_{[0,a[}{(x_n)}}
\leq\limsup_{N\rightarrow{\infty}}\frac{1}{N}\sum_{n=1}^{N}{\chi_{[0,a[}{(x_n)}}
\nonumber \\&\leq&\lim_{N\rightarrow{\infty}}\frac{1}{N}\sum_{n=1}^{N}{g_2(x_n)}
=\int_{0}^{1}{g_2(x)}\ dx\nonumber \\&\leq&\int_{0}^{1}{g_1(x)}\ dx+\varepsilon\leq\ a+\varepsilon.
\end{eqnarray*}
Since $\varepsilon$ is arbitrarily small, we have (\ref{succ-ud}).\\
\endproof
Moreover, we can state a more general result. 
\begin{theorem}\ \\
A sequence $(x_n)$ of points in $[0,1]$ is u.d.\! if and only if for any Riemann-integrable function $f$ defined on $[0,1]$ the equation (\ref{famDet}) holds.
\end{theorem}
\proof\ \\
The sufficiency follows directly from the previous theorem, because every continuous function is Riemann-integrable. The other implication was shown by De Bruijn and Post \cite{DBP}, who proved that if $f$ is defined on $[0,1]$ and if the averages $\frac{1}{N}\sum\limits_{n=1}^{N}{f(x_n)}$ admit limit for any $(x_n)$ u.d., then $f$ is Riemann-integrable.\\
\endproof
The problem of finding the largest reasonable determining classes has been addressed also in \cite{CV} and \cite{SV}.\\

Other examples of determing classes are the following ones.
\paragraph{Examples}
\begin{itemize}
\item The class of all characteristic functions of open (closed or half-open) subintervals of $[0,1]$ is determining.
\item The class of the characteristic functions of all intervals of the type $[0,q]$ with $q\in\mathbb{Q}$ is determining.
\item The class of all  step functions, i.e. functions given by finite linear combinations of characteristic functions of half-open subintervals of $[0,1]$ is determining.
\item The class of all  continuous (real or complex-valued) functions $g$ on $[0,1]$ such that $g(0)=g(1)$ is determining.
\item The class of all  polynomials with rational coefficients is determining.
\end{itemize}

Now, consider all functions of the type $f(x)=e^{2\pi ihx}$ where $h$ is a non-zero integer. One of the most important facts of uniform distribution theory is that these functions give a criterion to determine if a sequence of points is u.d..
\begin{theorem}[Weyl's Criterion]\ \\
The sequence $(x_n)$ is u.d.\! if and only if
\begin{equation*}
\lim_{N\rightarrow{\infty}}\frac{1}{N}\sum_{n=1}^{N}{e^{2\pi ihx_n}}=0
\end{equation*} 
for all integers $h\not=0$.
\end{theorem}
This important result was proved for the first time by Weyl in \cite{W3}, but a lot of proofs can be find in literature. Moreover, this criterion has a variety of applications in uniform distribution theory and also in the estimation of exponential sums. In particular, Weyl applied this theorem to the special sequence $(\{n\theta\})$, with $\theta$ irrational, to give a new proof of the following theorem.

Let us recall that for any $x\in\mathbb{R}$, we denote by $\{x\}$ the fractional part of $x$, which satisfies $\{x\}=x-[x]$, where $[x]$ is the integral part of $x$ (i.e the greatest integer less or equal to $x$).
\begin{theorem}\label{KronSeq}\ \\
Let $\theta$ be an irrational number. Then the sequence $(\{n\theta\})$ is u.d..
\end{theorem}
This result was independently estabilished by Weyl \cite{W1}, Bohl \cite{B} and Sierpi\' nski~\cite{S} in 1909-1910. The problem of the distribution of this special sequence has its origin in the theory of secular perturbations in astronomy and signs the beginning of the theory of u.d.\! sequences of points. Theorem~\ref{KronSeq} improves a previous theorem due to Kronecker, who proved that the points $e^{in\theta}$ are dense in the unit circle, whenever $\theta$ is an irrational multiple of $\pi$ (Kronecker's approximation theorem). For this reason the sequence $(\{n\theta\})$ with $\theta$ irrational is called \emph{Kronecker's sequence}.

Finally, it is important to underline that uniform distribution has also a measure-theoretic aspect. In fact, if we look at Definition \ref{def-ud}, we realize that a sequence $(x_n)$ of points in $[0,1]$ is u.d.\! if and only if the sequence of discrete measures  $\left(\frac{1}{n}\sum\limits_{i=1}^{n}\delta_{x_i}\right)$ converges weakly to the Lebesgue measure $\lambda$ on $[0,1]$, where $\delta_{t}$ is the Dirac measure concentrated in $t$.

The notion of weak convergence of measures represents the link between u.d.\! sequence of points and u.d.\! sequence of partitions. 

\subsection{Discrepancy of sequences}
As a quantitative measure of the distribution behaviour of a u.d.\! sequence we consider the so-called discrepancy, that is the maximal deviation between the empirical distribution of the sequence and the uniform distribution. This notion was studied for the first time in a paper of Bergstr\"om, who used the term ``Intensit\"atdispersion''(see~\cite{Berg}). The term discrepancy was probably coined by van der Corput. Moreover, the first intensive study of discrepancy is due to van der Corput and Pisot in \cite{VdCPis}.
\begin{definition}[Discrepancy]\label{def-dis}\ \\
Let $\omega_N=\{x_1,\ldots,x_N\}$ be a finite set of real numbers in $[0,1]$. The number
$$D_N(\omega_N)=\sup_{0\leq a<b\leq 1}\Bigg|\frac{1}{N}\sum_{i=1}^N\chi_{[a,b[}(x_i)-(b-a)\Bigg|
$$
is called the \emph{discrepancy} of the given set $\omega_N$.
\end{definition}
If $(x_n)$ is an infinite sequence of points, we associate to it the sequence of positive real numbers $D_N(\{x_1, x_2, \dots x_N\})$. So, the symbol $D_N(x_n)$ denotes the discrepancy of the initial segment $\{x_1, x_2, \dots x_N\}$ of the infinite sequence.

The importance of the concept of discrepancy in uniform distribution theory is revealed by the following fact (see \cite{W3} for more details).
\begin{theorem}\ \\
A sequence $(x_n)$ of points in $[0,1]$ is u.d.\! if and only if
$$\lim_{N\to\infty}D_N(x_n)=0.$$
\end{theorem} 
Sometimes it is useful to restrict the family of intervals considered in the definition of discrepancy. The most important type of restriction is to consider only intervals of the form $[0,a[$ with $0< a\leq 1$.
\begin{definition}[Star discrepancy]\label{def-stardis}\ \\
Let $\omega_N=\{x_1,\ldots,x_N\}$ be a finite set of real numbers in $[0,1]$, we define \emph{star discrepancy} of $\omega_N$ the quantity
$$D^*_N(\omega_N)=\sup_{0< a\leq 1}\Bigg|\frac{1}{N}\sum_{i=1}^N\chi_{[0,a[}(x_i)-a\Bigg|.
$$
\end{definition}
The definition $D^*_N$ is extended to the infinite sequence in the same way as we did for $D_N$. Moreover, the discrepancy and the star discrepancy are related by the following inequality.
\begin{theorem}\label{relDisStarDis}\ \\
For any sequence $(x_n)$ of points in $[0,1]$ we have
$$D^*_N(x_n)\leq D_N(x_n)\leq 2 D_N^*(x_n).$$
\end{theorem}
The most prominent open problem in theory of irregularities of distribution is to determine the optimal lower bound for the discrepancy. A first trivial lower bound is given by the following proposition.
\begin{prop}\label{trivialLowBound}\ \\
For any finite set $\omega=\{x_1,\ldots,x_N\}$ in $[0,1]$ we have that
$$\frac{1}{N}\leq D_N(\omega)\leq 1.$$
\end{prop}
The finite set $x_n=\frac{n}{N}$, $n=1,\ldots,N$ satisfies $D_N(\{x_1,\ldots,x_N\})=\frac{1}{N}$. But sequences of this kind can only exist in the one-dimensional case by a theorem due to Roth \cite{R} and this shows that the lower bound is optimal.
Moreover, in this example it is easy to see that for every $N$ a new set $\{x_1,\ldots, x_N\}$ is constructed. So the natural question is if there exists an infinite sequence $(x_n)$ in $[0,1]$ such that $D_N(x_n)=\mathcal{O}\left(\frac{1}{N}\right)$ as $N\to \infty$. Van der Corput made the conjecture that there are no sequences of this kind in the unit interval and this was proved by van Aardenne-Ehrenfest in \cite{vAE1, vAE2}. 
But the van der Corput conjecture was completely solved also from a quantitative point of view with the following important result due to Schmidt \cite{Sch}.
\begin{theorem}[Schmidt's Theorem]\label{SchmidtThm}\ \\
For any sequence $(x_n)$ in $[0,1]$ we have that
$$ND_N(x_n)>c\log N$$
for infinitely many positive integers $N$, where $c>0$ is an absolute constant.
\end{theorem}
This lower bound is the best possible in the one-dimensional case.\\

Usually, sequences having discrepancy of the order $\mathcal{O}\left(\frac{\log N}{N}\right)$ are called \emph{low discrepancy sequences} and they are very important for several applications. An interesting example of this kind of sequences are the van der Corput sequences.

\subsection{The van der Corput sequence}\label{sec:vdC}
In 1935 van der Corput introduced a procedure to generate low discrepancy sequences on $[0,1]$ (see \cite{VdC}). These sequences are considered the best distributed on $[0,1]$, because no infinite sequence has yet been found with discrepancy of smaller order of magnitude than the van der Corput sequences. The technique of van der Corput is based on a very simple idea.
First of all we have to define the radical inverse function which is at the basis of this construction.
\begin{definition}[Radical-inverse function]\label{RIFunct}\ \\
Let $b\geq 2$ an integer and let $n=\sum\limits_{k=0}^ra_kb^k$ be the digital expansion of the integer $n\geq 1$ in base $b$, $a_k\in\{0,\ldots,b-1\}$. The function
$$\gamma_b(n)=\sum_{k=0}^ra_kb^{-k-1}$$
is called radical inverse function in base $b$.
\end{definition}
The radical inverse function $\gamma_b(n)$ represents the fraction lying between $0$ and $1$ constructed by reversing the order of the digits in the $b-$adic expansion of $n$.
\begin{definition}[van der Corput sequences]\ \\
Let $b\geq 2$ a fixed prime integer. 
The sequence $(x_n)_{n\geq 1}$, where
$$x_n=\gamma_b(n-1),$$
is called \emph{van der Corput sequence in base $b$}.
\end{definition}

For example, the van der Corput sequence in base $b=2$ is given by
$$0, \frac{1}{2},\ \frac{1}{4},\ \frac{3}{4},\ \frac{1}{8},\ \frac{5}{8},\ \frac{3}{8},\ \frac{7}{8},\ \ldots$$
The construction of these points is explicitely showed in the following table.
\begin{center}
\begin{tabular}{|c|c|c|c|}
\hline
$\mathbf{n-1}$ & $(\mathbf{n-1})$ \textbf{in base 2}  & $\mathbf{x_n}$ \textbf{in base 2} &$\mathbf{x_n}$ \\
\hline
0 & $0.0$  & $0.0$ & 0 \\
\hline
 1 & $1.0$  & $ 0.1$ & $\frac{1}{2}$ \\
\hline
 2 & $10.0$ &  $0.01$ & $\frac{1}{4}$  \\
\hline
 3 & $11.0$ & $ 0.11$ & $\frac{3}{4}$ \\
\hline
 4 & $100.0$ &  $ 0.001$ &  $\frac{1}{8}$  \\
\hline
 5 & $101.0$ &  $0.101$ & $\frac{5}{8}$ \\
\hline
 6 & $110.0$ &  $0.011$ & $\frac{3}{8}$ \\
\hline
 7 & $111.0$ &  $0.111$ & $\frac{7}{8}$ \\
\hline
\end{tabular}
\end{center}
and so on.  
Let us introduce the estimate of the convergence order of the discrepancy of the van der Corput sequence in base $2$. But before proving this result, we need some further notions (see~\cite{KN}).
\begin{definition}\label{AlmArit}\ \\
For $0\leq\delta<1$ and $\varepsilon>0$, a finite set $x_1<x_2<\cdots<x_N$ in $[0,1]$ is said to be an \emph{almost-arithmetic progression} if there exists an $\eta$ with $0<\eta\leq\varepsilon$ such that the following conditions are satisfied:
\begin{itemize}
\item $0\leq x_1\leq \eta+\delta\eta$
\item $\eta-\eta\delta\leq x_{n+1}-x_n\leq \eta+\eta\delta$  for $1\leq n\leq N-1$
\item $1-\eta-\delta\eta\leq x_N<1$.
\end{itemize}
\end{definition}
It is clear that if $\delta=0$, then we have a genuine arithmetic progression with difference $\eta$.
\begin{theorem}\label{ThmAlmArit}\ \\
Let $x_1<x_2<\cdots<x_N$ be an almost-arithmetic progression in $[0,1]$ and let $\eta$ be the parameter in Definition \ref{AlmArit}. Then
$$D^*_N(x_1,\ldots,x_N)\leq \frac{1}{N}+\frac{\delta}{1+\sqrt{1-\delta^2}}\quad \text{for}\ \delta>0$$
$$D^*_N(x_1,\ldots,x_N)\leq \min\left\{\eta,\frac{1}{N}\right\}\quad \text{for}\ \delta=0.$$
\end{theorem}
\begin{theorem}\label{ThmDiscrDecomp}\ \\
Let $\omega=\{x_1,\ldots, x_N\}$ be a finite set of $N$ points in $[0,1]$. For $1\leq j\leq r$ let $\omega_j$ be a subset of $\omega$ consisting of $N_j$ elements such that its discrepancy is $D_{N_j}(\omega_j)$, its star discrepancy is $D^*_{N_j}(\omega_j)$, $\omega_j\cap \omega_i=\emptyset$ for all $j\neq i $ and $N= N_1+\ldots+N_r$. Then
$$ 
D_N(\omega)\leq\sum_{j=1}^{r}\frac{N_j}{N}D_{N_j}(\omega_j)
$$ 
and also
$$ 
D^*_N(\omega)\leq\sum_{j=1}^{r}\frac{N_j}{N}D^*_{N_j}(\omega_j).
$$ 
\end{theorem}
Now, we are ready to prove the following result.
\begin{theorem}\ \\
The discrepancy $D_N(x_n)$ of the van der Corput sequence in base $2$ satisfies
$$D_N(x_n)\leq c\left( \frac{\log(N+1)}{N}\right)$$
where $c>0$ is an absolute constant.
\end{theorem}
\proof\ \\
Let $N\geq 1$. We represent $N$ by its dyadic expansion 
$$ N=2^{h_1}+\ldots+2^{h_r}\qquad\text{with}\quad h_1>h_2>\ldots>h_r\geq 0.$$
Partition the interval $[1,N]\cap\mathbb{N}$ of integers in $r$ subsets $M_1,\ldots,M_r$ defined as follows
$$M_j=[2^{h_1}+\ldots+2^{h_{j-1}}+1, 2^{h_1}+\ldots+2^{h_{j-1}}+2^{h_j}]\cap\mathbb{N}\qquad\text{for}\quad 1< j\leq r$$
and put $M_1=[0, 2^{h_1}]\cap\mathbb{N}$.

An integer $n\in M_j$ can be written in the form 
$$n=1+2^{h_1}+\ldots+2^{h_{j-1}}+\sum_{i=0}^{h_j-1}a_i2^i, \quad\text{with}\quad a_i\in\{0,1\}.$$
In fact, we get all $2^{h_j}$ integers in $M_j$ if we let the $a_j$ run through all the possible combinations of $0$ and $1$. It follows that the point $x_n$ of the van der Corput sequence is given by
$$x_{n}=2^{-h_1-1}+\ldots+2^{-h_{j-1}-1}+\sum_{i=0}^{h_j-1}a_i2^{-i-1}=y_j+\sum_{i=0}^{h_j-1}a_i2^{-i-1}$$
where $y_j$ only depends on $j$ and not on $n$.

If $n$ runs through $M_j$, then the sum $\sum\limits_{i=0}^{h_j-1}a_i2^{-i-1}$ runs through all fractions $0, 2^{-h_j},\ldots, (2^{h_j}-1)\cdot 2^{-h_j}$. Moreover, we can note that $0\leq y_j< 2^{-h_j}$.

We conclude that if the elements $x_n$ with $n\in M_j$ are ordered according to their magnitude, then we obtain a sequence $\omega_j$ consisting of $N_j=2^{h_j}$ elements that is an arithmetic progression with parameters $\delta=0$ and $\eta=2^{-h_j}$, (see Definition~\ref{AlmArit}). By Theorem~\ref{ThmAlmArit}, we have that
$$D^*_{N_j}(\omega_j)\leq\min\left\{\eta,\frac{1}{N_j}\right\}=2^{-h_j}.$$

The set of the first $N$ terms of the van der Corput sequence, i.e. $\omega=\{x_1,\ldots,x_N\}$, can be decomposed in the $r$ subset $\omega_j$ defined above, since $ N=N_1+\cdots+N_r=2^{h_1}+\ldots+2^{h_r}$. Hence, by Theorem \ref{ThmDiscrDecomp} we have
\begin{equation}\label{e1}
D^*_N(\omega)\leq \sum_{j=1}^r\frac{N_j}{N}D^*_{N_j}(\omega_j)\leq\sum_{j=1}^r\frac{1}{N}=\frac{r}{N}.
\end{equation}

It remains to estimate $r$ in terms of $N$. Since $h_1>h_2>\ldots>h_r\geq 0$ then we have that
$h_r\geq 0\ ,\ h_{r-1}\geq 1 \ ,\ h_{r-2}\geq 2,\ldots, h_1\geq r-1.$ So we have that
$$N= 2^{h_1}+\ldots+2^{h_r}\geq 2^{r-1}+\ldots+2^{0}= 2^r-1,$$
and so 
\begin{equation}\label{e2}
r\leq\frac{\log(N+1)}{\log 2}.
\end{equation}

Finally, by combining (\ref{e1}) and (\ref{e2}) we have
$$D^*_N(\omega)\leq\frac{\log(N+1)}{N\log 2}$$ 
and since Theorem \ref{relDisStarDis} holds, we have
$$D_N(\omega)\leq\left(\frac{2}{\log 2}\right)\cdot\left(\frac{\log(N+1)}{N}\right).$$ 
\endproof

\section{Uniformly distributed sequences of partitions on $[0,1]$}\label{sec:partitions}
In this section, we will consider u.d.\! sequences of partitions of~$[0,1]$, a concept which has been introduced in 1976 by Kakutani in \cite{Ka}. In particular, we will sketch the theory of u.d.\! sequences of partitions introducing the significant example constructed by Kakutani. In the second part of this section, we will investigate the relation between u.d.\! sequences of partitions and u.d.\! sequences of points. This topic is analyzed more thoroughly in \cite{Vol}.\\

Firstly, let us give the basic definitions.
\begin{definition}\ \\ 
Let $(\pi_n)$ be a sequence of partitions of $[0,1]$, where $\pi_n=\{[t_{i-1}^n, t_i^n] : 1\leq i\leq k(n)\}$. The sequence $(\pi_n)$ is said to be uniformly distributed (u.d.)\! if for any continuous function $f$ on $[0,1]$ we have
\begin{equation}\label{def-udpart}
\lim_{n\rightarrow\infty}\frac{1}{k(n)}\sum_{i=1}^{k(n)}f(t_i^n)=\int_0^1 {f(t)\ dt}.  
\end{equation}
\end{definition}
Equivalently, $(\pi_n)$ is u.d.\! if the sequence of 
discrepancies
\begin{equation}\label{discrPart}
D_n=\sup_{0\leq a<b\leq 1}\bigg|\frac{1}{k(n)}\sum_{i=1}^{k(n)}
\chi_{[a,b[}(t_i^{(n)})- (b-a)\bigg|
\end{equation}
tends to $0$ as $n\to\infty$.

Similarly to the sequences of points, we can note that the uniform distribution of the sequence of partitions $(\pi_n)$ is equivalent to the weak convergence to $\lambda$ of the associated sequences of measures $(\nu_n)$, with  
\begin{equation}\label{AssMeas}
\nu_n=\frac{1}{k(n)}\sum_{i=1}^{k(n)}\delta_{t_i^n}.
\end{equation}

Moreover, it is easy to see that the uniform distribution of the sequence of partitions $(\pi_n)$ is equivalent to each of the following two conditions:
\begin{enumerate}
\item 
For any choice of the points $\tau_i^n\in[t_{i-1}^n, t_i^n]$ we have
$$
\lim_{n\rightarrow\infty}\frac{1}{k(n)}\sum_{i=1}^{k(n)}f(\tau_i^n)=\int_0^1 {f(t)\ dt}  
$$
for any continuous function $f$ on $[0,1]$.
\item For any choice of the points $\tau_i^n\in[t_{i-1}^n, t_i^n]$ we have that the sequence of measures
$$\frac{1}{k(n)}\sum_{i=1}^{k(n)}\delta_{\tau_i^n}$$
converges weakly to the Lebesgue measure $\lambda$ on $[0,1]$.
\end{enumerate}

\subsection{Kakutani's splitting procedure}\label{KakProc}
Let us describe a particular technique which allows to construct a whole class of u.d.\! sequences of partitions of $[0,1]$. This procedure was introduced by Kakutani in 1976 and works through successive $\alpha-$refinements of the unit interval \cite{Ka}.
\begin{definition}\ \\
If $\alpha\in]0,1[$ and $\pi=\{[t_{i-1}, t_i] : 1\leq i\leq k\}$ is any partition of
$[0, 1]$, then Kakutani's $\alpha$-refinement of $\pi$ (which will be denoted by $\alpha\pi$) is obtained by splitting only the intervals of $\pi$ having maximal lenght in two parts, proportional to $\alpha$ and $\beta=1-\alpha$ respectively.
\end{definition}
We will denote by $\alpha^2\pi$ the $\alpha$-refinement of $\alpha\pi$ and, in general, by $\alpha^n\pi$ the $\alpha-$refinement of $\alpha^{n-1}\pi$. Starting with the trivial partition $\omega$ of $[0,1]$, i.e. $\omega=\{[0,1]\}$, we get Kakutani's sequence of partitions $\kappa_n=\alpha^n\omega$. 

For example, if $\alpha<\beta$ we have that\\
$\kappa_1=\{[0,\alpha],[\alpha, 1]\}$\\
$\kappa_2=\{[0,\alpha],[\alpha, \alpha+\alpha\beta],[\alpha+\alpha\beta,1]\}$\\
and so on.

About this splitting procedure Kakutani proved the following result.
\begin{theorem}\label{KakThm}\ \\
For every $\alpha\in]0,1[$ the sequence of partitions $(\kappa_n)$ of $[0,1]$ is u.d..
\end{theorem} 
The most transparent proof of this theorem is due to Adler and Flatto and follows from a combination of classical results from ergodic theory \cite{AF}. Indeed, Kakutani's procedure caught the attention of several authors in the late seventies also from a stochastic point of view. In fact, Kakutani's theorem was a partial answer to the following question posed by the physicist H. Araki, which regarded random splitting of the interval $[0,1]$. Let $X_1$ be choosen randomly with respect to the uniform distribution on $[0,1]$. Once $X_1,\ldots, X_n$ have been choosen, let $X_{n+1}$ be a point picked at random and accordingly to the uniform distribution in the largest of the $n+1$ intervals determined by the previous $n$ points. Kakutani had been originally asked whether the associated sequence of empirical distribution functions converges uniformly, with probability 1, to the distribution function of the uniform random variable on $[0,1]$. 

This question has been studied in \cite{VZ, LO1, LO2, BD} and later in \cite{PvZ}. It is important to note that in the probabilistic setting the possibility that the partition obtained at the $n-$th step has more than one interval of maximal lenght can be neglected, since it is an event which has probability equal to zero. On the other hand, in Kakutani's splitting procedure for every $\alpha$ the partition $\alpha^n\omega$ has more than one interval of maximal lenght for infinitely many values of $n$.

Recently, some new results and ideas revived the interest for this subject. In fact, Kakutani's technique has been generalized in several directions. In \cite{Vol-Car} the splitting procedure has been extended to higher dimensions, providing a sequence of nodes in the hypercube $[0,1]^d$ which is proved to be u.d.. In \cite{Vol-Car2} a von Neumann type theorem is presented for sequences of partitions of $[0,1]$. More precisely, u.d.\! sequences of partitions of the unit interval are constructed starting from sequences of partitions $\pi_n$ whose diameter tends to zero for $n\to\infty$. In \cite{Vol} the concept of $\alpha-$refinement is generalized and it is introduced a new splitting procedure for constructing a larger class of u.d.\! sequences of partitions on $[0,1]$. Moreover, in this paper it is analyzed the deep relation between the theory of u.d.\! sequences of partitions and the theory of u.d.\! sequences of points. This strong connection between the two theories makes more interesting the study of u.d.\! sequences of partitions in view of possible applications to Quasi-Monte Carlo methods. 

\subsection{Associated uniformly distributed sequences of points}\label{subsec: AssSeqPoints}
In the following, we intend to study the problem of associating to a u.d.\! sequence of partitions a u.d.\! sequence of points. Before investigating this problem, let us note that the converse problem results to be easier in many cases. 
\begin{theorem}\ \\
If $(x_n)$ is a u.d.\! sequence of points in $[0,1]$ such that $x_n\neq x_m$ when $n\neq m$ and $x_n\notin\{0,1\}$ for any $n\in\mathbb{N}$, then the sequence of partitions $(\pi_n)$, where each $\pi_n$ is determined by the points $\{0, 1, x_k \  \text{with}\   k\leq n\}$ ordered by magnitude, is u.d..
\end{theorem}
\proof\ \\
By using the assumption that $(x_n)$ is u.d.\! and Theorem \ref{continueDet}, it follows that the relation (\ref{def-udpart}) holds for any continuous function $f$ defined on $[0,1]$.\\
\endproof
The requirement that  $x_n\neq x_m$ when $n\neq m$ is important and it is not possible to avoid this assumption in the theorem as it is shown in the following example.

\paragraph{Example}\ \\
Consider the sequence $(x_n)$ defined by consecutive blocks of $4m$ points for $m\in\mathbb{N}$. Each block is defined as follows
$$\left\{\frac{1}{2m+1},\frac{1}{2m+1},\ldots,
\frac{m}{2m+1},\frac{m}{2m+1},\ldots,\frac{1}{2},\frac{2m+1}{4m},
\frac{2m+2}{4m},\ldots,\frac{4m-1}{4m}\right\}.$$
In each block the first $m$ points are repeated twice, while the others are all distinct. In this way, the points of the sequence have double density in the right half of $[0,1]$, but they have however a good distribution because of the repetition in the left half of $[0,1]$. So the sequence $(x_n)$ is u.d.. But when we take in consideration the sequence of partitions $(\pi_n)$ associated to $(x_n)$, according to the procedure described in the previous theorem, the repetitions are cancelled. Hence, we get a sequence $(\pi_n)$ having twice as many subintervals in $\big[\frac{1}{2},1\big]$ than in $\big[0,\frac{1}{2}\big[$ and so $(\pi_n)$ is not u.d.. \\
\vspace{0.5mm}\\
Now, consider our starting problem of associating a u.d.\! sequence of points to a fixed u.d.\! sequence of partitions. Let us introduce an important result proved by Vol\v{c}i\v{c} in \cite{Vol}, where a probabilistic answer to this problem is given. 

Suppose $(\pi_n)$ is a u.d.\! sequence of partitions in $[0,1]$ with $\pi_n=\{[t_{i-1}^n, t_i^n]: 1\leq i\leq k(n)\}$. The natural question is if it is possible to rearrange the points $t_i^n$ determining the partitions $\pi_n$, for $1\leq i\leq k(n)$, in order to get a u.d.\! sequence of points. Clearly, there exist many ways of reordering the points $t_i^n$. A natural restriction is that we first reorder all the points determining $\pi_1$ then those defining $\pi_2$, and so on. This kind of reorderings are called \emph{sequential reorderings}.\\

Before presenting the result of Vol\v{c}i\v{c}, we need some preliminaries. In particular, we introduce a version of the strong law of large numbers for negatively correlated random variables, which is attributed to Aleksander Rajchman and can be proved following the lines of Theorem 5.1.2 in \cite{Chu}.

\begin{lemma}\label{lemmaAlek}\ \\
Let $(\varphi_n)$ be a sequence of real, negatively correlated random variables with variances uniformly bounded by $V$ on the probability space $(W,P)$. Moreover, suppose that 
$$\lim_{i\to\infty}E(\varphi_i)=M.$$ Then
$$\lim_{n\to\infty}\frac{1}{n}\sum_{i=1}^n\varphi_i=M\qquad\text{almost surely}.$$
\end{lemma}
\proof\ \\
We may assume $E(\varphi_i)=0$ and remove afterwards this restriction by applying the conclusions to the sequence of random variables $\varphi_i-E(\varphi_i)$.

Put $S_n=\sum\limits_{i=1}^n\varphi_i$. For any $\varepsilon>0$, by using the \v{C}ebi\v{s}ev inequality we have
\begin{eqnarray*}
P\left(\frac{1}{n^2}S_{n^2}\geq\varepsilon\right)
&\leq&\frac{1}{\varepsilon^2} Var\left(\frac{1}{n^2}S_{n^2}\right)
=\frac{1}{n^4\varepsilon^2}E\left(S^2_{n^2}\right)\\
&=&\frac{1}{n^4\varepsilon^2}E\left(\left(\sum\limits_{i=1}^{n^2}\varphi_i\right)^2\right)\\
&=&\frac{1}{n^4\varepsilon^2}E\left(\sum\limits_{i=1}^{n^2}\varphi_i^2
+\sum\limits_{i=1}^{n^2}\sum\limits_{\stackrel{i\neq j}{j=1}}^{n^2}\varphi_i\varphi_j\right)\\
&=&\frac{1}{n^4\varepsilon^2}\left(\sum\limits_{i=1}^{n^2}E\left(\varphi_i^2\right)
+\sum\limits_{i=1}^{n^2}\sum\limits_{\stackrel{i\neq j}{j=1}}^{n^2}E\left(\varphi_i\varphi_j\right)\right).
\end{eqnarray*}
Now, because of the negative correlation of the $\varphi_i$'s we have that the terms $E\left(\varphi_i\varphi_j\right)$  for $i\neq j$ are not positive. So by using this fact and the bound for the variance, we get the estimate
$$P\left(\frac{1}{n^2}S_{n^2}\geq\varepsilon\right)
\leq\frac{1}{n^4\varepsilon^2}\left(\sum\limits_{i=1}^{n^2}E\left(\varphi_i^2\right)\right)
=\frac{1}{n^4\varepsilon^2}\left(\sum\limits_{i=1}^{n^2}Var\left(\varphi_i\right)\right)
\leq\frac{V}{n^2\varepsilon^2}. $$
Since the series of the upper bounds is convergent, the series $$\sum\limits_{n=1}^\infty P\left(\frac{1}{n^2}S_{n^2}\geq\varepsilon\right)$$ is convergent, too. Therefore by the Borel-Cantelli lemma, we have that
\begin{equation}\label{uno}
\lim_{n\to\infty}\frac{1}{n^2}S_{n^2}=0 \quad\text{a.s.}.
\end{equation}
Define now
$$L_n=\max_{n^2\leq j<(n+1)^2}\left|S_j-S_{n^2}\right|.$$
For the same $\varepsilon$, the \v{C}ebi\v{s}ev inequality implies that
\begin{eqnarray*}
P\left(\frac{L_{n}}{n^2}\geq\varepsilon\right)&\leq&\frac{1}{n^4\varepsilon^2}E\left(L^2_n\right)
\leq\frac{1}{n^4\varepsilon^2}E\left(\sum_{j=n^2+1}^{(n+1)^2-1}\left|S_j-S_{n^2}\right|^2\right)\\
&=&\frac{1}{n^4\varepsilon^2}E\left(\sum_{j=n^2+1}^{(n+1)^2-1}\left(\sum_{i=n^2+1}^j\varphi_i\right)^2\right)\\
&=&\frac{1}{n^4\varepsilon^2}E\left(\sum_{j=n^2+1}^{(n+1)^2-1}\left(\sum_{i=n^2+1}^j\varphi_i^2
+\sum_{i=n^2+1}^j\sum_{\stackrel{h=n^2+1}{i\neq h}}^j\varphi_i\varphi_h\right)\right)\\
&=&\frac{1}{n^4\varepsilon^2}\sum_{j=n^2+1}^{(n+1)^2-1}\left(\sum_{i=n^2+1}^jE(\varphi_i^2)
+\sum_{i=n^2+1}^j\sum_{\stackrel{h=n^2+1}{i\neq h}}^jE(\varphi_i\varphi_h)\right)\\
&\leq&\frac{1}{n^4\varepsilon^2}\sum_{j=n^2+1}^{(n+1)^2-1}\sum_{i=n^2+1}^jVar(\varphi_i)\\
&\leq&\frac{1}{n^4\varepsilon^2}\sum_{j=n^2+1}^{(n+1)^2-1}\sum_{i=n^2+1}^{(n+1)^2-1}Var(\varphi_i)\\
&\leq& \frac{V(2n-1)^2}{n^4\varepsilon^2}.
\end{eqnarray*}
Since the series of the upper bounds is convergent, the series $$\sum\limits_{n=1}^\infty P\left(\frac{1}{n^2}L_{n}\geq\varepsilon\right)$$ is convergent and therefore, again by the Borel-Cantelli lemma, we have
\begin{equation}\label{due}
\lim_{n\to\infty}\frac{1}{n^2}L_{n}=0 \quad\text{a.s.}.
\end{equation}
Since for any $m$ with $n^2\leq m<(n+1)^2$ we have
$$\frac{\left|S_m\right|}{m}\leq\frac{1}{n^2}\left(\left|S_{n^2}\right|+L_n\right)$$
the conclusion follows from (\ref{uno}) and (\ref{due}).\\
\endproof

Let $\varphi$ be the random variable taking with probability $\frac{1}{k}$ values in the sample space $W=\{w_i\in[0,1], 1\leq i\leq k\}$ with $k\geq 2$. We assume that $w_{i-1}<w_i$ for $1\leq i\leq k$. Denote by $\varphi_i$ the value assumed by $\varphi$ in the $i-$th draw from $W$ without replacement. Fix $c\in]0,1[$ and let $\psi_i=\chi_{[0,c[}(\varphi_i)$. Then the following property holds.
\begin{prop}\label{negcorr}\ \\
The variances of the random variables $\psi_i$, $1\leq i\leq k$, are bounded by $\frac{1}{4}$ and the $\psi_i$'s are negatively correlated.
\end{prop}
\proof\ \\
The expectation of $\psi_i$ is given by
$$E(\psi_i)=\frac{1}{k}\sum_{i=1}^{k}\chi_{[0,c[}(\omega_i),$$
so $E(\psi_i^2)=E(\psi_i)$. Then
$$Var(\psi_i)=E(\psi_i^2)-\left(E(\psi_i)\right)^2=E(\psi_i)\left(1-E(\psi_i)\right).$$
Now, it is easy to see that $\frac{1}{4}$ is an upper bound for the right-hand side and so we have that
$$Var(\psi_i)\leq\frac{1}{4}.$$ 
Since all pairs of distinct $\psi_i$'s have the same joint distribution, we
may evaluate just the covariance of $\psi_1$ and $\psi_2$.
Suppose that $w_i\in[0,c[$ if and only if $i\leq h$, with $0<h<k$. Then
\begin{eqnarray*}
Cov(\psi_1,\psi_2)&=&E(\psi_1\psi_2)-E(\psi_1)E(\psi_2)\\
&=&\frac{1}{k(k-1)}\sum_{i=1}^{h}\sum_{\stackrel{j=1}{i\neq j}}^h\chi_{[0,c[}(w_i)\chi_{[0,c[}(w_j)-
\left(\frac{1}{k}\sum_{i=1}^{h}\chi_{[0,c[}(w_i)\right)^2\\
&=&\frac{h(h-1)}{k(k-1)}-\frac{h^2}{k^2}=\frac{h(h-k)}{k^2(k-1)}< 0
\end{eqnarray*}
\endproof

Now, we are ready to introduce the result of Vol\v{c}i\v{c} (see \cite{Vol}). In the following, we consider the \emph{sequential random reordering} of the points $(t_i^n)$, defined as follows.
\begin{definition}\ \\
If $(\pi_n)$ is a u.d.\! sequence of partitions of $[0,1]$ with $\pi_n=\{[t_{i-1}^n, t_i^n]: 1\leq i\leq k(n)\}$, the \emph{sequential random reordering} of the points $t_i^n$ is a sequence $(\varphi_m)$ made up of consecutive blocks of random variables. The $n$-th block consists of $k(n)$ random variables which have the same law and represent the drawing, without replacement, from the sample space $W_n=\left\{t_1^n,\ldots,t_{k(n)}^n\right\}$ where each singleton has probability $\frac{1}{k(n)}$.
\end{definition}
Denote by $T_n$ the set of all permutations on $W_n$, endowed with the natural probability $P$ compatible with the uniform probability on $W_n$, i.e. $P(\tau_n)=\frac{1}{k(n)!}$ with $\tau_n\in T_n$.

Any sequential random reordering of $(\pi_n)$ corresponds to a random selection of $\tau_n\in T_n$ for each $n\in\mathbb{N}$. The permutation $\tau_n\in T_n$ identifies the reordered $k(n)$-tuple of random variables $\varphi_i$ with $K(n-1)\leq i\leq K(n)$, where $K(n)=\sum\limits_{i=1}^nk(i)$. Therefore, the set of all sequential random reorderings can be endowed with the natural product probability on the space $T=\prod\limits_{n=1}^\infty T_n$.

\begin{theorem}\label{RandReord}\ \\
If $(\pi_n)$ is a u.d.\! sequence of partitions of $[0,1]$, then the sequential random reordering of the points $t_i^n$ defining them is almost surely a u.d.\! sequence of points in $[0,1]$.
\end{theorem}

\proof \ \\
Let $(\varphi_m)$ be the sequential random reordering of $(\pi_n)$. First of all, note that if  $0<c<1$ and $\varphi_m$ belongs to the $n-$th block of $k(n)$ random variables, then
$$E(\chi_{[0,c[}(\varphi_m))=\frac{1}{k(n)}\sum_{i=1}^{k(n)}\chi_{[0,c[}(t_i^n)$$
and this quantity tends to $c$, when $m$ and hence $n$ tends to infinity, since $(\pi_n)$ is u.d.\! by assumption.

If we consider $\psi_m=\chi_{[0,c[}(\varphi_m)$ for $K(n-1)\leq m\leq K(n)$, then Proposition \ref{negcorr} holds and so the $\psi_m$'s are negatively correlated for $K(n-1)\leq m\leq K(n)$, i.e. when the $\varphi_m$ belong to the same block. On the other hand, the correlation is zero when the $\varphi_m$ belong to different blocks, since they are independent.

Let $\{c_h, h\in\mathbb{N}\}$ be a dense subset of $[0,1]$. Fix $h\in\mathbb{N}$ and consider the sequence $\left(\chi_{[0,c_h[}(\varphi_m)\right)$. Hence, we may apply the Lemma \ref{lemmaAlek} and get that
$$\lim_{n\to\infty}\frac{1}{n}\sum_{i=1}^n\chi_{[0,c_h[}(\varphi_i)=c_h\qquad\text{a.s.}$$
for any $c_h$. But this is a sufficient condition for the uniform distribution and so we have our conclusion.\\
\endproof

\section{Uniform distribution theory on $[0,1]^d$}\label{UDmultidim}
In this section we deal with the extension of uniform distribution theory to the unit hypercube. We will introduce the basic definitions and results of the theory with a particular attention to the study of discrepancy and to some special u.d.\! sequences of points in this space.
\subsection{Definitions and basic properties}
Let $d$ be an integer with $d\geq 2$. Let $J=[a_1,b_1[ \times \cdots \times [a_d,b_d[\subset \mathbb{R}^d$ be a rectangle with sides parallel to the axes in the $d-$dimensional space $\mathbb{R}^d$. If we denote by $\lambda_d$ the $d-$dimensional Lebesgue measure, then the volume of $J$ is given by $$\lambda_d(J)=\prod\limits_{i=1}^d(b_i-a_i).$$ Let us denote by $I^d$ the $d-$dimensional unit hypercube, i.e. $I^d=[0,1]^d$.
\begin{definition}\label{def-udmultidim}\ \\
A sequence $(x_n)$ of points in $I^d$ is said to be \emph{uniformly distributed (u.d.)} if for any rectangle $R$ of the form $R=[0,a_1[\times\cdots\times [0,a_d[\subset I^d$ we have
\begin{equation}\label{udmultidim}
\lim_{N\rightarrow{\infty}}\frac{1}{N}\sum_{n=1}^{N}{\chi_{R}{(x_n)}}=\lambda_d(R)
\end{equation}
where $\chi_{R}$ is the characteristic function of the rectangle $R$.
\end{definition}
As in the one-dimensional case we can introduce the concept of \emph{determining class} of functions.
\begin{definition}\ \\
A class $\mathcal{F}$ of Riemann-integrable functions on $I^d$ is said to be \emph{determining} for the uniform distribution of sequences of points, if for any sequence $(x_n)$ in $I^d$ the validity of the relation
\begin{equation}\label{famDetMultidim}
\lim_{N\rightarrow{\infty}}\frac{1}{N}\sum_{n=1}^{N}{f(x_n)}=\int_{I^d}f\ d\lambda_d \end{equation}
for all $f\in\mathcal{F}$ already implies that $(x_n)$ is u.d. .
\end{definition}

Weyl was the first to extend to the multidimensional case the uniform distribution theory. So, we can give also in this case his classical results \cite{W2,W3}.
\begin{theorem}[Weyl's Theorem]\ \\
A sequence $(x_n)$ of points in $I^d$ is u.d.\! if and only if for any (real or complex-valued) continuous function $f$ defined on $I^d$ the equation (\ref{famDetMultidim}) holds.
\end{theorem}

Moreover, let $x=(x_1,\ldots,x_d)$ and $y=(y_1,\ldots,y_d)$ be in $\mathbb{R}^d$ and let us denote by $x\cdotp y$ the usual inner product in $\mathbb{R}^d$, i.e. $x\cdotp y=\sum\limits_{i=1}^dx_iy_i$. Then we can give the generalization of the Weyl's Criterion.
\begin{theorem}[Weyl's Criterion]\ \\
The sequence $(x_n)$ in $I^d$ is u.d.\! if and only if
$$
\lim_{N\rightarrow{\infty}}\frac{1}{N}\sum_{n=1}^{N}{e^{2\pi ih\cdotp x_n}}=0
$$
for all non-zero integer lattice points $h\in\mathbb{Z}^d-\{(0,\ldots,0)\}$.
\end{theorem}
Weyl applied this theorem to Kronecker's sequence also in the multidimensional case for giving a new proof of Kronecker's approximation theorem in $\mathbb{R}^d$ (see \cite{W3}).
\begin{theorem}[Kronecker's Approximation Theorem]\label{KronMultidim} \ \\ 
Let $\theta=(\theta_1,\ldots, \theta_d)\in\mathbb{R}^d$ such that $1,\theta_1,\ldots, \theta_d$ are linearly independent over the rationals.\! Then the sequence of fractionals parts $(\{ n\theta\})$, where $\{ n\theta\}=(\{ n\theta_1\},\ldots,\{ n\theta_d\})$, is dense in $I^d$. 
\end{theorem}
Furthermore, Weyl's criterion implies that a sequence of the form $( n\theta)$ is u.d.\! if and only if $1,\theta_1,\ldots, \theta_d$ are linearly independent over $\mathbb{Q}$ . Hence it follows that $( n\theta)$ is u.d.\! if and only $(\{ n\theta\})$ is dense in $I^d$.

\subsection{Estimation of discrepancy}
Definitions \ref{def-dis} and \ref{def-stardis} may be extended to sequences of points in $I^d$ as follows.
\begin{definition}\ \\
Let $\omega_N=\{x_1,\ldots,x_N\}$ be a finite set of points in $I^d$. 
\begin{itemize}
\item The discrepancy of $\omega_N$ is defined by
$$D_N(\omega_N)=\sup_{J}\Bigg|\frac{1}{N}\sum_{i=1}^N\chi_{J}(x_i)-\lambda_d(J)\Bigg|,
$$
where $J$ runs through all rectangles in $I^d$ of the form $J=[a_1,b_1[ \times \cdots \times [a_d,b_d[$ with $0\leq a_i<b_i\leq 1$.
\item The star discrepancy of $\omega_N$ is defined by
$$D^*_N(\omega_N)=\sup_{R}\Bigg|\frac{1}{N}\sum_{i=1}^N\chi_{R}(x_i)-\lambda_d(R)\Bigg|,
$$
where $R$ runs through all rectangles in $I^d$ of the form $R=[0,a_1[\times\cdots\times [0,a_d[$ with $0< a_i\leq 1$.
\end{itemize}
\end{definition}
Moreover, the discrepancy and the star discrepancy are related by the following inequality.
\begin{theorem}\ \\
For any sequence $(x_n)$ of points in $I^d$ we have
$$D^*_N(x_n)\leq D_N(x_n)\leq 2^d D_N^*(x_n).$$
\end{theorem}
In the same way as in the one-dimensional case if $(x_n)$ is an infinite sequence of points, we associate to it the sequence of positive real numbers $D_N(\{x_1, x_2, \dots x_N\})$. So, the symbol $D_N(x_n)$ denotes the discrepancy of the initial segment $\{x_1, x_2, \dots x_N\}$ of the infinite sequence. It is easy to see that
\begin{theorem}\ \\
A sequence $(x_n)$ of points in $I^d$ is u.d.\! if and only if
$$\lim_{N\to\infty}D_N(x_n)=0.$$
Equivalently a sequence $(x_n)$ of points in $I^d$ is u.d.\! if and only if
$$\lim_{N\to\infty}D^*_N(x_n)=0.$$
\end{theorem} 

The immediate lower bound given in Proposition \ref{trivialLowBound} holds also in the higher- dimensional case. In fact, we get the following inequality.
\begin{prop}\ \\
For any finite set $\omega=\{x_1,\ldots,x_N\}$ of points in $I^d$ we have that
$$\frac{1}{N}\leq D_N(\omega)\leq 1.$$
\end{prop}
\proof\ \\
The right-hand side inequality is evident from the definition of discrepancy. Now, choose $\varepsilon>0$ and consider the first point of $\omega$, namely $x_1=\left(x_1^{(1)},\ldots, x_1^{(d)}\right)\in I^d$. Let $J=[x_1^{(1)}, x_1^{(1)}+\varepsilon[\times\cdots\times[x_1^{(d)}, x_1^{(d)}+\varepsilon[$. Since $x_1\in J$ then we have
$$D_N(\omega)\geq \frac{1}{N}-\lambda_d(J)=\frac{1}{N}-\varepsilon^d$$
and so the conclusion follows.\\
\endproof
As we have already said, only in the one-dimensional case we have examples of sequences such that $D_N(x_n)=\frac{1}{N}$. In fact, in the higher-dimensional case such examples cannot exist by Roths's theorem \cite{R}. So far this is the best known result for $d>3$.
\begin{theorem}[Roth's Theorem]\ \\
Let $d\geq 2$. Then the discrepancy $D_N(x_n)$ of the finite set $\omega=\{x_1,\ldots,x_N\}\subset I^d$ is bounded from below by
$$D_N(\omega)\geq c_d\left(\frac{(\log N)^{\frac{d-1}{2}}}{N}\right),$$
where $c_d>0$ is an absolute constant given by $c_d=\frac{1}{2^{4d}((d-1)\log 2)^\frac{d-1}{2}}$.
\end{theorem} 
For further information on bounds for the dimensions 2 and 3 and refinements of Roth's theorem we refer to \cite{DT}.

A well known conjecture states that for every dimension $d$ there exists a constant $c_d$ such that for any infinite sequence $(x_n)$ in $\mathbb{R}^d$ with $d\geq 1$ we have
$$D_N(x_n)\geq c_d \left(\frac{(\log N)^d}{N}\right)$$
for infinitely many $N$. This conjecture has been proved by Schmidt only for $d=1$ (see Theorem \ref{SchmidtThm}), while it is still open for $d\geq 2$.

Usually, sequences of points in $\mathbb{R}^d$ having discrepancy bounded from above by $\mathcal{O}\left(\frac{(\log N)^d}{N}\right)$ are called \emph{low discrepancy sequences}.
We have already described an important class of low discrepancy sequences in the one-dimensional case, that is the van der Corput sequences. In the following, we will introduce their higher-dimensional generalization. Before defining these special u.d.\! sequences, let us give a result that proves the important role played by low discrepancy sequences in numerical integration.

\subsection{The Koksma-Hlawka inequality}\label{subsec: KHIneq}
The concept of discrepancy gives a quantitative measure of the order of convergence in the relation (\ref{udmultidim}) defining the uniform distribution of a given sequence. Consequently, it is also very interesting to get information on the order of convergence in (\ref{famDetMultidim}). Referring to this problem, a very useful estimate is provided by the Koksma-Hlawka inequality. In fact, it states that the order of convergence of the difference between the actual value of the integral in (\ref{famDetMultidim}) and its approximation can be estimated in terms of the variation of the function and the star discrepancy. Before we can write down this result, we need to define the variation of a function $f: I^d\to\mathbb{R}$. 

By a partition $P$ of $I^d$ we mean a set of $d$ finite sequences $(\eta_i^{(0)},\ldots,\eta_i^{(m_i)})$ for $i=1,\ldots,d$ with $0=\eta_i^{(0)}\leq\eta_i^{(1)}\leq\cdots\leq\eta_i^{(m_i)}=1$. In connection with such a partition we define for each $i=1,\ldots,d$ an operator $\Delta_i$ by
\begin{eqnarray*}
\Delta_if(x_1,\ldots,x_{i-1},\eta_i^{(j)},x_{i+1},\ldots,x_d)&=&
f(x_1,\ldots,x_{i-1},\eta_i^{(j+1)},x_{i+1},\ldots,x_d)\\
&-&f(x_1,\ldots,x_{i-1},\eta_i^{(j)},x_{i+1},\ldots,x_d)
\end{eqnarray*} 
for $0\leq j<m_i$. Operators with different subscrites obviously commute and 
$\Delta_{i_1,\ldots,i_k}$ stands for $\Delta_{i_1}\cdots\Delta_{i_k}$. Such an operator commutes with summation over variables on which it does not act.

\begin{definition}[Function of bounded variation in the sense of Vitali]\ \\
For a function $f: I^d\to\mathbb{R}$ we set
$$V^{(d)}(f)=\sup_P\sum_{j_1=0}^{m_1-1}\cdots\sum_{j_d=0}^{m_d-1}
\left|\Delta_{1,\ldots,d}f(\eta_1^{(j_1)},\ldots,\eta_d^{(j_d)})\right|,$$
where the supremum is extended over all partitions $P$ of $I^d$. \\
If $V^{(d)}(f)$ is finite then $f$ is said to be of bounded variation on $I^d$ in the sense of Vitali.
\end{definition}

\begin{definition}[Function of bounded variation in the sense of Hardy and Krause]
Let $f: I^d\to\mathbb{R}$ and assume that $f$ is of bounded variation in the sense of Vitali.
If the restriction $f^{(F)}$ of $f$ to each face $F$ of $I^d$ of dimension $1,2,\ldots,d-1$ is of bounded variation on $F$ in the sense of Vitali, then $f$ is said to be of bounded variation on $I^d$ in the sense of Hardy and Krause.
\end{definition}

So we can state the following theorem.
\begin{theorem}[Koksma-Hlawka's Inequality]\ \\
Let $f$ be a function of bounded variation on $I^d$ in the sense of Hardy and Krause. Let $\omega=(x_1,\ldots, x_N)$ be a finite set of points in $I^d$. Let us denote by $\omega_l$ the projection of $\omega$ on the $(d-l)-$dimensional face $F_l$ of $I^d$ defined by $F_l=\{(u_1,\ldots,u_d)\in I^d : u_{i_1}=\cdots$ $\cdots=u_{i_l}=1\}$. Then we have
\begin{equation}\label{KHIn}
\left|\frac 1N \sum_{n=1}^N f(x_n)-\int_{I^d}f(x)dx\right|\leq
\sum_{l=0}^{d-1}\sum_{F_l}D_N^*(\omega_l)V^{(d-l)}(f^{(F_l)}),
\end{equation}
where the second sum is extended over all $(d-l)-$dimensional faces $F_l$ of the form $u_{i_1}=\cdots=u_{i_l}=1$. The discrepancy $D_N^*(\omega_l)$ is clearly computed in the face of $I^d$ in which $\omega_l$ is contained.
\end{theorem}
\begin{rem}\ \\
Trivially $D_N^*(\omega_l)$ can be bounded by $D_N^*(\omega)$. Hence we get from (\ref{KHIn}) that
\begin{equation}\label{KHIneq}
\left|\frac 1N \sum_{n=1}^N f(x_n)-\int_{I^d}f(x)dx\right|\leq V(f)D^*_N(\omega)
\end{equation}
where
$$V(f)=\sum_{l=0}^{d-1}\sum_{F_l}V^{(d-l)}(f^{(F_l)})$$
is called the \emph{variation of Hardy and Krause}.
\end{rem}
A proof can be found in \cite{KN}, but the original proof is given in \cite{Hlaw}. This relation provides a strong motivation for the choice of low discrepancy sequences in Quasi-Monte Carlo integration.
 
\subsection{The Halton and Hammersley sequences}
A very important application of u.d.\! sequences is numerical integration. In fact, given a function  $f$ on $I^d$, the basic idea of classical Monte Carlo integration is to approximate the integral
$$I(f)=\int_{I^d}fd\lambda_d$$
with the mean $$I_N(f)=\frac{1}{N}\sum_{i=1}^Nf(x_i)$$
where $x_1,\ldots,x_N$ are $N$ points choosen randomly or pseudorandomly in $I^d$.

For a large class of functions, Quasi-Monte Carlo methods have a faster rate of convergence than Monte Carlo methods. Indeed, the Quasi-Monte Carlo method works by choosing deterministically the $N$ integration points instead of actual random points. Therefore, it is essential that the nodes are well distributed on $I^d$. This means that it is convenient if their distribution is close to the uniform distribution. A good choice for the integration points is the initial segment of a sequence $(x_n)$ with small discrepancy, since the Koksma-Hlawka inequality holds, i.e
$$\left|I_N(f)-I(f)\right|\leq V(f)D_N^*(x_n), $$
where $V(f)$ is the variation of $f$ in the sense of Hardy-Krause (see Subsection \ref{subsec: KHIneq}).

Finally, the deterministic nature of Quasi-Monte Carlo methods provides many advantages with respect to Monte Carlo methods. First of all, the Quasi-Monte Carlo method allows to work with deterministic points rather than random samples and then it offers the availability of deterministic error bounds instead of the probabilistic Monte Carlo rate of convergence. Moreover, with the same computational effort, the Quasi-monte Carlo method achieves a significantly higher accuracy than the Monte Carlo method just thanks to the choice of the integration points with small discrepancy.

In this subsection, we want to introduce some important classes of sequences of points in $I^d$ with small discrepancy: the Halton sequences and the Hammersley sequences. Both constructions are based on the radical inverse function (see Definition~\ref{RIFunct}).
\begin{definition}[Halton sequence]\ \\
For a given dimension $d\geq 2$ the \emph{$d-$dimensional Halton sequence} $(x_n)$ in $I^d$ is defined by 
$$x_n=\left(\gamma_{b_1}(n),\ldots, \gamma_{b_d}(n)\right)$$
where $b_1,\ldots, b_d$ are given coprime integers.
\end{definition}

As it was shown in \cite{Hal}, the Halton sequence is a low discrepancy sequence. In fact, it  has a discrepancy of order $\mathcal{O}\left(\frac{(\log N)^d}{N}\right)$.

For $d=1$ we just get the van der Corput sequence (see Subsection~\ref{sec:vdC}). So, Halton's construction is a generalization of the van der Corput one to the higher-dimensional case.

For example, let us consider $b_1=2$ and $b_2=3$. By applying Halton's construction we first have to generate the van der Corput sequence in base~$2$ that is $(\gamma_2(n))$, i.e.
$$ \frac{1}{2}, \frac{1}{4}, \frac{3}{4}, \frac{1}{8}, \frac{5}{8}, \frac{3}{8}, \frac{7}{8}, \ldots$$
and then we have to generate the van der Corput sequence in base $3$ that is $(\gamma_3(n))$, i.e.
$$ \frac{1}{3}, \frac{2}{3}, \frac{1}{9}, \frac{4}{9}, \frac{7}{9}, \frac{2}{9}, \frac{5}{9}, \frac{8}{9},\ldots$$

Finally, the Halton sequence $(\gamma_2(n),\gamma_3(n))$ in the unit square $I^2$ is obtained by pairing up these two sequences
$$\left(\frac{1}{2}, \frac{1}{3}\right), \left(\frac{1}{4}, \frac{2}{3}\right), \left(\frac{3}{4}, \frac{1}{9}\right), \left(\frac{1}{8}, \frac{4}{9}\right), \left(\frac{5}{8}, \frac{7}{9}\right), \left(\frac{3}{8}, \frac{2}{9}\right), \left(\frac{7}{8}, \frac{5}{9}\right),\ldots$$ 

While the performance of standard Halton sequences is very good in low dimensions, problems with correlation have been observed among sequences generated from higher primes. This can cause serious problems in the estimation of models with high-dimensional integrals. In order to deal with this problem, various other methods have been proposed; one of the most prominent solutions is the technique of \emph{scrambled Halton sequence}, which uses permutations of the coefficients employed in the construction of the standard sequences \cite{Schl, BW, MorCal}.

\begin{definition}[Hammersley sequence]\ \\
For given integers $d\geq 2$ and $N$, the \emph{$d-$dimensional Hammersley sequence} $(x_n)$ of size $N$ in $I^d$ is defined by 
$$x_n=\left(\frac{n}{N}, \gamma_{b_1}(n),\ldots, \gamma_{b_{d-1}}(n)\right)$$
where $b_1,\ldots, b_{d-1}$ are given coprime integers.
\end{definition}
As it was shown in \cite{Ham}, the Hammersley sequence has a discrepancy of order $\mathcal{O}\left(\frac{(\log N)^{d-1}}{N}\right)$.

Note that the Hammersley sequence is a finite set of size $N$ which cannot be extended to an infinite sequence. So in the approximation of the integral $I(f)$, one should decide in advance the value of $N$ in order to perform the calculation, since the first coordinate depends on $N$. In the computational practice of Quasi-Monte Carlo integration it is often convenient to be able to increase the value of $N$ without losing the previously calculated function values. For this purpose, it is preferable to work with a whole low discrepancy sequence of nodes and then take its first $N$ terms whenever a value of $N$ has been selected. In this way, $N$ can be increased while all data from the earlier computations can be still used. Therefore, in several cases the Halton sequences are more convenient in Quasi-Monte Carlo integration than the Hammersley point sets.

\section{Uniform distribution theory in compact spaces}
A theory of uniform distribution can be developed in settings more abstract than the unit interval and the unit hypercube. In this section, we present its generalization to compact Hausdorff spaces. The study of this theory was intiated by Hlawka in \cite{Hlaw1,Hlaw2}. The notion of u.d.\! sequences in such spaces is related to a given non-negative regular normalized Borel measure, but for convenience we will consider a regular probability. 

Let $X$ be a compact Hausdorff space and let us denote by $\cal B$ the $\sigma$-algebra of Borel subsets of $X$. Suppose $\mu$ is a regular probability on $\cal B$. \begin{definition}[Regular Borel measure]\label{RegMis} \ \\
A positive measure $\mu$ defined on $\cal B$ is said to be \emph{regular} if
$$
\mu(E)=\sup\{\mu(C)\colon C\subseteq E,\ C \ \text{closed}\}=\inf\{\mu(D)\colon E\subseteq D,\ D \ \text{open}\}$$ for all $E\in{\cal B}$. 
\end{definition}

\begin{notation}
Let us denote by:
\begin{itemize}
\item $\mathcal{B}(X)$ the set of all bounded real-valued measurable functions defined on~$X$
\item $\mathcal{C}(X)$ the subset of $\mathcal{B}(X)$ consisting of all continuous real-valued functions defined on $X$.
\end{itemize}
\end{notation}

The space $\mathcal{B}(X)$, endowed with the norm $\|f\|_{\infty}=\sup\limits_{x\in X}|f(x)|$, is a Banach space.

Among the various characterizations of the concept of u.d.\! sequences of points in $[0,1]$ the most easily adaptable to this general situation is Weyl's Theorem (that is Theorem~\ref{continueDet}), which allows to give the following definition.

\begin{definition}[U.d.\! sequences of points]\ \\
A sequence $(x_i)$ of elements in $X$ is said to be \emph{uniformly distributed (u.d.)\! with respect to $\mu$}, if 
$$
    \lim_{N\rightarrow\infty}\frac{1}{N}\sum_{i=1}^N f(x_i)=\int_{X}f(t)\ d\mu(t)
$$
for all $f\in\mathcal{C}(X)$.
\end{definition}

In order to generalize to compact spaces the concepts of u.d.\! sequences of partitions, we need to introduce the notion of $\mu-$continuity set.

\begin{definition}[$\mu$-continuity set]\ \\
A Borel set $M\subset X$ is called a \emph{$\mu$-continuity set} if $\mu(\partial M)=0$, where $\partial M$ denotes the boundary of $M$ with respect to the relative topology on~$X$.
\end{definition}

\begin{definition}[U.d.\! sequences of partitions]\ \\
Let $(\pi_n)$ be a sequence of partitions of $X$, where $\pi_n=\Big\{A_1^n,A_2^n,\ldots,A_{k(n)}^n\Big\}$ and the $A^n_i$'s are $\mu$-continuity sets. The sequence $(\pi_n)$ is said to be \emph{uniformly distributed (u.d.)\! with respect to $\mu$} if for any $f\in {\cal C}(X)$, and any choice $t^n_i\in A^n_i$ we have
$$
\lim_{n\rightarrow\infty}\frac{1}{k(n)}\sum_{i=1}^{k(n)}f(t_i^n)=\int_X {f(t)\ d\mu(t)} .  
$$
\end{definition}

The existence of u.d.\! sequences of partitions in separable metric spaces has been addressed, but not completely solved, in \cite{CV}. On the other hand, the existence problem for u.d.\! sequences of points can be easily settled in compact Hausdorff spaces satisfying the second axiom of countability. In fact, if $X$ is a compact Hausdorff with countable basis, then almost all sequences in $X$ are u.d.\! with respect to $\mu$ \cite{Hlaw1,KN}. 
Neverthless, the existence problem is still open in the general setting of compact Hausdorff spaces. The strongest constructive result is due to Hedrl\' in, who showed that u.d.\! sequences of points exist in every compact metric space using an explicit construction in \cite{Hedr}. Interesting results on this topic are proposed in \cite{Nied}.

\begin{definition}[Determining functions for sequences of points]\ \\
A class $\cal F$ of Riemann-integrable functions is said to be {\it determining} for the uniform distribution of sequences of points with respect to $\mu$ if for any sequence $(x_i)$ in $X$ the validity of the relation
\begin{equation}\label{determinante}
\lim_{N\rightarrow\infty}\frac{1}{N}\sum_{i=1}^{N}f(x_i)=\int_X{f(t)\ d\mu(t)}\,,
\end{equation}
for all $f\in \cal F$ already implies that $(x_i)$ is u.d..
\end{definition}

Similarly, we can give the analogous definition for sequences of partitions.

\begin{definition}[Determining functions for sequences of partitions]\ \\
A class $\cal F$ of Riemann-integrable functions is said to be {\it determining} for the uniform distribution of sequences of partitions with respect to $\mu$ if for any sequence $(\pi_n)$, where $\pi_n=\Big\{A_1^n,A_2^n,\ldots,A_{k(n)}^n\Big\}$ and the $A^n_i$'s are $\mu$-continuity sets, the validity of the relation
$$
\lim_{n\rightarrow\infty}\frac{1}{k(n)}\sum_{i=1}^{k(n)}f(t_i^n)=\int_X {f(t)\ d\mu(t)}
$$
for all $f\in \cal F$ and for any choice $t^n_i\in A^n_i$ already implies that $(\pi_n)$ is u.d..
 \end{definition}

Observe that the determining classes for the sequences of points play the same role for the sequences of partitions and viceversa.

As in uniform distribution theory on $[0,1]$, a family of $\mu-$continuity set $\cal G$ such that the class $\mathcal{F}=\left\{\chi_M , M\in\cal G\right\}$ is determining is called \emph{discrepancy system}. Obviously, it is possible to define this notion independently of the concept of the determining class (see \cite{DT}).

\begin{definition}[Discrepancy system]\label{DisSystem}\ \\
A system $\cal G$ of $\mu-$continuity sets of $X$ is called \emph{discrepancy system} if 
$$\lim_{N\to\infty}\sup_{M\in\cal G}\left|\frac{1}{N}\sum_{n=1}^N\chi_M(x_n)-\mu(M)\right|=0$$
holds if and only if the sequence $(x_n)$ is u.d.\! with respect to $\mu$.
\end{definition}

For a family of real-valued functions $\mathcal{F}$, we will denote by \emph{$\mbox{span}(\mathcal{F})$} the linear space generated by $\mathcal{F}$ and by $\overline{\mbox{\emph{span}}(\mathcal{F})}$ its closure. The construction of many important determining classes is based on the following theorem.

\begin{theorem}\label{teoSpan}\ \\
Let $(x_n)$ be a sequence of points in $X$. If $\mathcal{F}$ is a class of functions from $\mathcal{B}(X)$ such that (\ref{determinante}) holds for all $f\in \cal F$ and $\overline{\mbox{\emph{span}}(\mathcal{F})}\supset\mathcal{C}(X)$, then $\mathcal{F}$ is a determining class for $(x_n)$.
\end{theorem}
\proof\ \\
Let us first show that (\ref{determinante}) holds for all $g\in\mbox{span}(\mathcal{F})$. In fact any $g\in\mbox{span}(\mathcal{F})$ is of the form $g=\alpha_1f_1+\ldots+\alpha_kf_k$ with $f_i\in\mathcal{F}$ and $\alpha_i\in\mathbb{R}$, $1\leq i\leq k$. Since (\ref{determinante}) holds for all $f\in\mathcal{F}$, in particular holds for all $f_i$. Therefore by linearity, the function~$g$ satisfies the relation~(\ref{determinante}).

Now, let us consider $f\in\mathcal{C}(X)$. Fixed $\varepsilon>0$, by the assumption of density there exists $h\in\mbox{span}(\mathcal{F})$ such that $\|f-h\|_{\infty}<\varepsilon$. Then we have
\begin{eqnarray}
\Bigg|\frac{1}{N}\sum_{n=1}^Nf(x_n)-\int_X f\ d\mu\Bigg|&\leq&
\Bigg|\frac{1}{N}\sum_{n=1}^N(f-h)(x_n)-\int_X (f-h)\ d\mu\Bigg|\nonumber\\&+&\Bigg|\frac{1}{N}\sum_{n=1}^Nh(x_n)-\int_X h\ d\mu\Bigg|\nonumber\\ 
&\leq&\frac{1}{N}\sum_{n=1}^N|(f-h)(x_n)|+\int_X |f-h|\ d\mu \nonumber\\&+&\Bigg|\frac{1}{N}\sum_{n=1}^Nh(x_n)-\int_X h\ d\mu\Bigg|\nonumber\\
&\leq& 2\|f-h\|_{\infty}+\Bigg|\frac{1}{N}\sum_{n=1}^Nh(x_n)-\int_X h\ d\mu\Bigg|\nonumber\\&<&3\varepsilon\nonumber
\end{eqnarray}
for sufficiently large $N$.\\
\endproof

Now, we can generalize to compact spaces the concept of discrepancy.

\begin{definition}[ $\cal G-$discrepancy]\label{Gdiscrepancy}\ \\
Let $\cal G$ be a discrepancy system in $X$ and $\omega_N=\{x_1,\ldots,x_N\}$ a finite set of points in $X$. Then the discrepancy with respect to $\cal G$ (or $\cal G-$discrepancy) is defined by
$$D_N^{\cal G}(\omega_N)=\sup_{M\in\cal G}\Bigg|\frac{1}{N}\sum_{n=1}^N\chi_M(x_n)-\mu(M)\Bigg|.
$$
\end{definition}

If $(x_n)$ is an infinite sequence of points in $X$, we associate to it the sequence of positive real numbers $D_N^{\cal G}(\{x_1, x_2, \dots x_N\})$. Often it is used the symbol $D_N^{\cal G}(x_n)$ to denote the quantity $D_N^{\cal G}(\{x_1, x_2, \dots x_N\})$.  

It follows from the definition that $(x_n)$ is u.d.\! if and only if $D_N^{\cal G}(x_n)$ tends to zero when $N$ tends to infinity.

%% file: Cap2.tex
\chapter{Uniform distribution on fractals}
In this chapter, we will be concerned with uniform distribution theory on a special class of fractals, namely those which are defined by an Iterated Function System (IFS) of similarities having the same ratio and satisfying the Open Set Condition (OSC). More precisely, we will give an explicit procedure to generate u.d.\! sequences of partitions and of points on this class of fractals and we will present some results about the elementary discrepancy of these sequences~\cite{Inf-Vol}.

\section{Fractals defined by Iterated Function Systems}
Let us introduce a general method of construction for some fractals which is based on their self-similarity, that is the property of many fractals to be made up of parts similar to the whole. For instance, the Cantor set is given by the union of two similar copies of itself and the von Koch curve consists of four similar copies. This property may actually be used to define these fractals, which are called IFS fractals because they are generated by an Iterated Function System. Before introducing this kind of construction, let us recall some basic notions (see \cite{Fal2, Fal1}).\\

Let $\|\cdot\|$ be the usual norm on the $d-$dimensional Euclidean space $\mathbb{R}^d$. By the diameter of a set $U\subset\mathbb{R}^d$ we mean the quantity diam$(U)=\sup\limits_{x,y\in U}\|x-y\|$.
\begin{definition}[$\delta-$covering]\ \\
Let $E\subset\mathbb{R}^d$. Fixed $\delta>0$, a countable family $(U_i)$ of sets of $\mathbb{R}^d$ is said to be a \emph{$\delta-$covering} of $E$ if the union of all $U_i$'s covers the set $E$ and for each $i$ we have $0<$\emph{diam}$(U_i)\leq\delta$.
\end{definition}

\begin{definition}[$s-$dimensional Hausdorff measure]\ \\
For $E\subset\mathbb{R}^d$ and $s\geq 0$, we define for each $\delta>0$
$$
\mathcal{H}_\delta^s(E)=\inf\Bigg\{\sum_{i=1}^\infty \emph{diam}(U_i)^s:\  (U_i)\ \text{is a}\   \delta-\text{covering of}\ E \Bigg\}
$$
and subsequently the $s-$dimensional Hausdorff measure of $E$ is given by
$$\mathcal{H}^s(E)=\lim_{\delta\rightarrow 0}\mathcal{H}_\delta^s(E).$$
\end{definition}
$\mathcal{H}^s$ is not a measure in the usual sense, but it is an outer measure. In fact, in general the countable additivity does not hold but it is possible to prove that $\mathcal{H}^s$ is a measure only when it is defined over the Borel sets of $\mathbb{R}^d$ (see \cite{Fal1}). The Hausdorff measure generalizes the concept of the Lebesgue measure on $\mathbb{R}^d$. Indeed, $\mathcal{H}^d$ is equal to the $d-$dimensional Lebesgue measure $\lambda_d$ up to a constant, i.e. $$\mathcal{H}^d(E)=c_d\lambda_d(E)$$
where $c_d=\frac{\pi^{\frac{1}{2}d}}{2^d\big(\frac{1}{2}d\big)!}$.
Moreover, the Hausdorff measure has a very useful scaling property.
\begin{prop}\label{scalingProp}\ \\
Let $E\subset\mathbb{R}^d$, $k>0$ and $s\geq 0$ then
\begin{equation*}
\mathcal{H}^s(kE)=k^s\mathcal{H}^s(E)
\end{equation*}
where $kF=\{kx\colon x\in E\}$, i.e. the set $kE$ is the set $E$ scaled of a factor $k$.
\end{prop}
$\mathcal{H}^s(E)$ is non-increasing with $s$ and there exists a unique value of $s$ where $\mathcal{H}^s(E)$ jumps from $\infty$ to $0$. This value is called \emph{Hausdorff dimension} of $E$ and it is given by
$$\dim_H(E)=\inf\{s\colon\mathcal{H}^s(E)=0\}=\sup\{s\colon\mathcal{H}^s(E)=\infty\}.$$

Let us denote by $\mathcal{K}(\mathbb{R}^d)$ the space of all the non-empty compact subsets of~$\mathbb{R}^d$ endowed with the Hausdorff distance, which makes it a complete metric space. Let us recall the definition of \emph{Hausdorff distance}. Let $K$ and $L$ be two non-empty subsets of~$\mathbb{R}^d$, then we define their Hausdorff distance $d_H(K,L)$ by 
$$
d_H(K,L)=\max\left\{\sup_{x\in K}\inf_{y\in L}\|x-y\|, \sup_{x\in L}\inf_{y\in K}\|x-y\| \right\}
$$
or equivalently
$$
d_H(K,L)=\min\left\{ \lambda\geq 0 \colon K\subset L_\lambda \ \text{and}\ L\subset K_\lambda \right\}
$$
where
$$K_\lambda=\left\{x\in\mathbb{R}^n \colon \|x-y\|\leq\lambda\ \text{for some}\ y\in K\right\}$$
and
$$L_\lambda=\left\{x\in\mathbb{R}^n \colon \|x-y\|\leq\lambda\ \text{for some}\ y\in L\right\}.$$

Let us give the following results, due to Hutchinson, which show how an IFS defines a unique non-empty self-similar compact set (see \cite{Hutch}).

\begin{theorem}\ \\
Let $\psi_1,\ldots,\psi_m$ be $m$ contractions defined on $\mathbb{R}^d$ so that
$\|\psi_i(x)-\psi_i(y)\|\leq c_i\|x-y\| \,$
for all $x,y\in\mathbb{R}^d $, with $0<c_i<1$ for each $i$. Then the mapping $\psi(E) \mapsto \bigcup\limits_{i=1}^m\psi_i(E)$ is a contraction on $\mathcal{K}(\mathbb{R}^d)$ and its unique fixed point is a non-empty compact set $F$, called the {\it attractor} of the IFS. The set $F$ is said to be a \emph{self-similar} set and we have $$F=\bigcup_{i=1}^m\psi_i(F).$$ 
Moreover, if $F_0\in\mathcal{K}(\mathbb{R}^d)$ is such that $\psi_i(F_0)\subset F_0$ for $1\leq i\leq m$, then the sequence of iterates $\left(\psi^n(F_0)\right)$ is decreasing and convergent to $F$ in the Hausdorff metric as $n\rightarrow\infty$, with
$$F=\bigcap_{n=0}^\infty \psi^n(F_0)$$
(where $\psi^0(F_0)=F_0 \ \textrm{and} \ \psi^{n+1}(F_0)=\psi(\psi^{n}(F_0))\ \textrm{for} \ n\geq 0$).
\end{theorem}

\noindent The set $F_0=\psi^0(F_0)$ is called \emph{initial set} and the iterates $\psi^n(F_0)$ are called \emph{pre-fractals} for $F$. 

One advantage of dealing with fractals generated by an IFS is that their Hausdorff dimension is often easy to calculate. In particular, the evaluation of the Hausdorff dimension is very simple when we consider $m$ similarities $\psi_1,\ldots,\psi_m$ on $\mathbb{R}^d$ with ratios $0<c_i<1$ for each $i$, i.e $\|\psi_i(x)-\psi_i(y)\|= c_i\|x-y\| \,,$
for all $x,y\in\mathbb{R}^d $ and assume that the following condition holds.
\begin{definition}[OSC]\ \\
A class of similarities $\psi_1,\ldots,\psi_m\colon\mathbb{R}^d\rightarrow\mathbb{R}^d$ satisfies the \emph{open set condition} if there exists a non-empty bounded open set $V$ such that
\begin{equation*}
  V\supset\bigcup_{i=1}^m\psi_i(V) 
\end{equation*}
where $\psi_i(V)$ are pairwise disjoint.
\end{definition}

Then for this special class of IFS fractals we have the following theorem due to Moran (see \cite{Moran}, \cite{Fal2}).

\begin{theorem}\label{MoranThm}\ \\
Assume that $m$ similarities $\psi_1,\ldots,\psi_m$ defined on $\mathbb{R}^d$ with ratios $0<c_i<1$ (for $i=1,\ldots, m$) satisfy the OSC. Let $F$ be the attractor of the $\psi_i$'s then the Hasudorff dimension $s$ of $F$ is the solution of the equation
\begin{equation}\label{dimMoran}
    \sum_{i=1}^mc_i^s=1.
\end{equation} 
Moreover, we have that the $s-$dimensional Hausdorff measure $\mathcal{H}^s(F)$ is positive and finite, i.e. $0<\mathcal{H}^s(F)<\infty$.
\end{theorem}
Note that, the OSC ensures that the components $\psi_i(F)$ of the invariant set $F$ cannot overlap too much and this property is expressed by the following corollary of the previous theorem.

\begin{corollary}\label{imp}\ \\
Let $\psi_1,\ldots,\psi_m$ be $m$ similarities on $\mathbb{R}^d$ with ratios $0<c_i<1$ for each $i$ and let $F$ be their attractor. If the OSC holds, then $\mathcal{H}^s(\psi_i(F)\cap \psi_j(F))=0 \ \textrm{for}\ i\neq j$.
\end{corollary}
\proof \ \\
Using the assumption that the $\psi_i$'s are similarities and Proposition \ref{scalingProp}, we have  
$$
\sum_{i=1}^m\mathcal{H}^s\left(\psi_i(F)\right)=\sum_{i=1}^mc_i^s\mathcal{H}^s(F).
$$
Since the OSC holds, we can use the relation (\ref{dimMoran}) and so we have 
$$
\sum_{i=1}^m\mathcal{H}^s\left(\psi_i(F)\right)=\mathcal{H}^s(F)\sum_{i=1}^mc_i^s=\mathcal{H}^s(F)
=\mathcal{H}^s\left(\bigcup\limits_{i=1}^m\psi_i(F)\right).
$$
By Theorem \ref{MoranThm} we have $0<\mathcal{H}^s(F)<\infty$, so the previous relation can only happen if $\mathcal{H}^s\left(\psi_i(F)\cap\psi_j(F)\right)=0$ for $i\neq j$.\\
\endproof

Let us cite some of the most popular examples of fractals which are included in the class considered by Theorem \ref{MoranThm}.
\newpage
\begin{ex}\label{es}\ 
\begin{description}
\item[Cantor set]\ \\ The Cantor set is constructed starting from the unit interval through a sequence of deletion operations. Put $C_0=[0,1]$. At the first step we remove the open middle third of $C_0$ and so we obtain the set $C_1=\big[0,\frac{1}{3}\big]\cup\big[\frac{2}{3},1\big]$. Deleting the open middle third of the intervals $\big[0,\frac{1}{3}\big]$ and $\big[\frac{2}{3},1\big]$, we obtain four intervals of lenght $\frac{1}{9}$. So at the end of the second step we have constructed the set $C_2=\big[0,\frac{1}{9}\big]\cup\big[\frac{2}{9},\frac{1}{3}\big]
\cup\big[\frac{2}{3},\frac{7}{9}\big]\cup\big[\frac{8}{9},1\big]$.
By repeating this procedure, at the $k-$th step we have that $C_k$ consists of $2^k$ intervals of lenght $3^{-k}$, generated by removing the open middle third of each interval in $C_{k-1}$. In Figure \ref{figCantor} the first three steps of this construction are illustrated. The Cantor set $C$ is given by 
$$C=\bigcap_{k=0}^\infty C_k.$$
\begin{figure}[!h]
\begin{center}
\includegraphics[width=14cm]{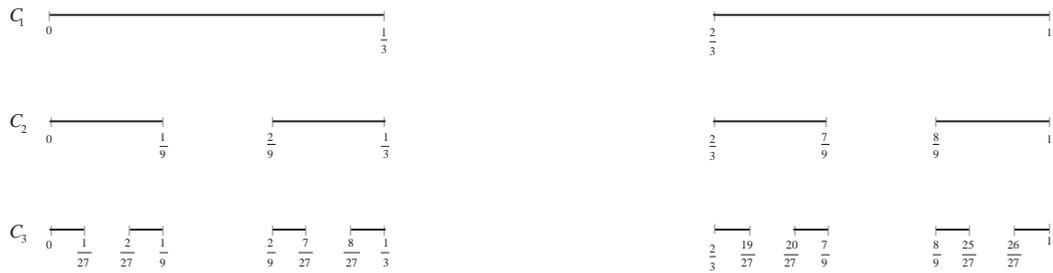}
\end{center}
\caption{Construction of the Cantor set}
\label{figCantor}
\end{figure}\\
So the Cantor set $C$ is the attractor of the two following similarities
\begin{equation*}
 \left\{ \begin{array}{ll}
S_1(x)=\frac{1}{3}x & \textrm{}\\ \\
S_2(x)=\frac{1}{3}x+\frac{2}{3} & \textrm{}
\end{array}\right.. 
\end{equation*}
The set $C$ satisfies the OSC taking $V=]0,1[$. Then by Theorem \ref{MoranThm} we have that the Hausdorff dimension $s$ of $C$ is given by
$$\sum_{i=1}^2\bigg(\frac{1}{3}\bigg)^s=1\ \Rightarrow\ 2=3^s\ \Rightarrow\ s=\frac{\log 2}{\log 3}.$$

\item[Sierpi\' nski Triangle]\ \\ The Sierpi\' nski triangle $T$ is constructed starting from an equilateral triangle by repeatedly removing inverted equilateral triangles. In fact, let $T_0$ be an equilateral triangle in $\mathbb{R}^2$ and take the three middle points of its sides. These three points and the vertices of $T_0$ define four equilateral congruent triangles and we remove the central open one.  At the end of the first step we have obtained three congruent closed triangles and we denote their union by $T_1$. At the second step we repeat this procedure on each triangle of $T_1$, so we get nine triangles whose union is $T_2$. In Figure \ref{figSierp} the first four steps of this construction are illustrated. The Sierpi\' nski triangle $T$ is given by 
$$T=\bigcap_{k=0}^\infty T_k.$$ 
\begin{figure}[!ht]
\begin{center}
\includegraphics[width=14cm]{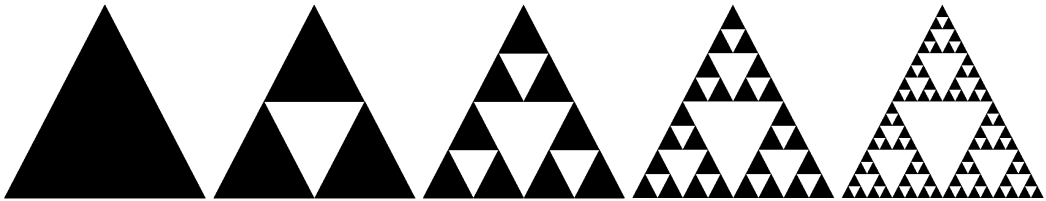}
\end{center}
\caption{Construction of the Sierpi\' nski triangle}
\label{figSierp}
\end{figure}

So the Sierpi\' nski triangle $T$ is the attractor of the following similarities
\begin{equation*}
 \left\{ \begin{array}{ll}
S_1(x,y)=\bigg(\frac{1}{2}x,\frac{1}{2}y\bigg) & \textrm{}\\ \\
S_2(x,y)=\bigg(\frac{1}{2}x+\frac 14,\frac{1}{2}y+\frac{\sqrt{3}}{4}\bigg) & \textrm{}\\ \\
S_3(x,y)=\bigg(\frac{1}{2}x+\frac{1}{2},\frac{1}{2}y\bigg) & \textrm{}
\end{array}\right. 
\end{equation*}
where the origin is taken in the left down vertex of the initial triangle. The set $T$ satisfies the OSC taking $V$ as the interior of the initial triangle $T_0$. Consequently, by Theorem~\ref{MoranThm}, we have that the Hausdorff dimension $s$ of $T$ is given by
$$\sum_{i=1}^3\bigg(\frac{1}{2}\bigg)^s=1\ \Rightarrow\ 3=2^s\ \Rightarrow\ s=\frac{\log 3}{\log 2}.$$

\item[von Koch Curve]\ \\ The von Koch curve $K$ is constructed starting from the unit interval $K_0=[0,1]$. At the first step we remove the open middle third of $K_0$ and replace it by the other two sides of the equilateral triangle based on the removed segment. The union of these four segments is denoted by $K_1$. We construct $K_2$ applying this procedure to each segment in $K_1$, and so on. The sequence of polygonal curves $K_j$ tends to a limiting curve $K$, called von Koch curve. 
In Figure \ref{figKoch} the first four steps of this construction are illustrated.\\
\begin{figure}[!h]
\begin{center}
\includegraphics[height=7.5cm, width=8cm]{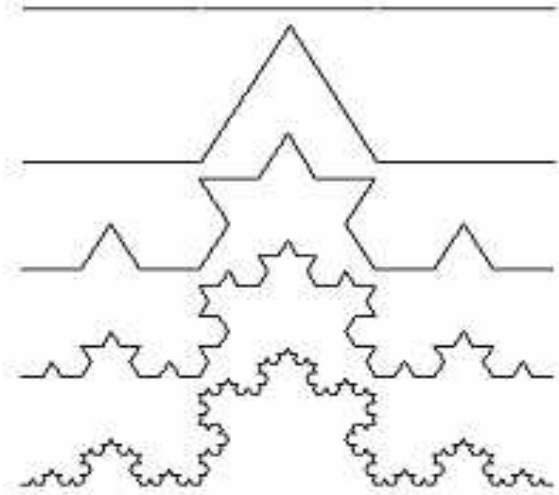}
\end{center}
\caption{Construction of the von Koch curve}
\label{figKoch}
\end{figure}\\
The von Koch curve $K$ is the attractor of the following similarities.
\begin{equation*}
 \left\{ \begin{array}{ll}
S_1(x,y)=\bigg(\frac{1}{3}x,\frac{1}{3}y\bigg) & \textrm{}\\ \\
S_2(x,y)=\bigg(\frac{1}{6}x-\frac{\sqrt{3}}{6}y+\frac{1}{3},\frac{\sqrt{3}}{6}x+\frac{1}{6}y\bigg) & \textrm{}\\ \\
S_3(x,y)=\bigg(-\frac{1}{6}x+\frac{\sqrt{3}}{6}y+\frac{2}{3},\frac{\sqrt{3}}{6}x+\frac{1}{6}y \bigg) & \textrm{}\\ \\
S_4(x,y)=\bigg(\frac{1}{3}x+\frac{2}{3},\frac{1}{3}y\bigg) & \textrm{}
\end{array}\right.. 
\end{equation*}
The curve $K$ satisfies the OSC taking $V$ as the interior of the isosceles triangle of height equal to $\frac{\sqrt{3}}{6}$ with basis the segment $K_0$. So by Theorem \ref{MoranThm} we have that the Hausdorff dimension $s$ of $K$ is given by
$$\sum_{i=1}^4\bigg(\frac{1}{3}\bigg)^s=1\ \Rightarrow\ 4=3^s\ \Rightarrow\ s=\frac{\log 4}{\log 3}.$$
\end{description}
\end{ex}

\section{Van der Corput sequences on fractals}\label{UDFractals}
In this section, we extend to certain fractals the concept of u.d.\! sequences of partitions defined for the interval~$[0,1]$ in Subsection \ref{sec:partitions}. In particular, we introduce our recent results in this setting \cite{Inf-Vol}. We present a general algorithm to produce u.d.\! sequences of partitions and of points on the class of fractals generated by a system of similarities on $\mathbb{R}^d$ having the same ratio and verifying the open set condition. Moreover, we provide an estimate for the elementary discrepancy of van der Corput type sequences constructed on this class of fractals. 

\subsection{Algorithm to construct u.d.\! sequences of points and of partitions on a class of fractals}\label{subsec: Algo}
The classical concept of u.d.\! sequences of points is more natural when we deal with the interval $[0,1]$ and with manifolds. On the other hand when we work on fractals, in particular with fractals generated by iterated function systems, partitions become a convenient tool for introducing a uniform distribution theory.

The advantage of considering partitions was implicitely used by Grabner and Tichy in \cite{Grab-Tichy} and by Cristea and Tichy in \cite{CT}, even if they treated u.d.\! sequences of points. In these papers various concepts of discrepancy were introduced on the planar Sierpi\' nski gasket and on the multidimensional Sierpi\' nski carpet respectively, by using different kinds of partitions on these two fractals. In \cite{Grab-Tichy} an analogue of the classical van der Corput sequence has been constructed on the planar Sierpi\' nski gasket. Similarly, in a succesive paper of Cristea, Pillichshammer, Pirsic and Scheicher~\cite{CT2} a sequence of van der Corput type has been defined on the $s$-dimensional Sierpi\' nski carpet by exploiting the IFS-addresses of the carpet points. In all these papers the order of convergence of the several notions of discrepancy is determined for the van der Corput type sequences constructed on these fractals.

The idea to study this special kind of sequences in relation to uniform distribution on IFS fractals is also our starting point. In fact, the algorithm we are going to introduce generalizes the results cited above and allows to construct van der Corput sequences on a whole class of fractals including the ones considered in \cite{Grab-Tichy}, \cite{CT} and~\cite{CT2}.
 
From now on, we consider $m$ similarities $\psi_1,\ldots,\psi_m$ defined on $\mathbb{R}^d$ having all the same ratio, i.e. for each $i$ we  have
$\|\psi_i(x)-\psi_i(y)\|= c\|x-y\| \,$
for all $x,y\in\mathbb{R}^d $ with $0<c<1$.
Moreover, we assume that our system of similarities satisfies the OSC. According to Theorem~\ref{MoranThm}, the Hausdorff dimension of the attractor $F$ of this IFS is $s=-\frac{\log m}{\log c}$ and its $s$-dimensional Hausdorff measure ${\cal H}^s$ is positive and finite.

Our class of fractals includes the most popular fractals as for instance the ones cited in Examples \ref{es}. But also $[0,1]$ can be seen as the attractor of an IFS, in fact of infinitely many IFS's. Indeed, fix a positive integer $m>1$ and consider the mappings $\varphi_1, \dots, \varphi_m$ from $\mathbb{R}$ to~$\mathbb{R}$, where
\begin{equation}\label{vdCIfs}
\varphi_k (x)=  \frac{k-1}{m} + \frac{x}{m}\ , \textrm{\ for\ } 1\le k \le m.
\end{equation}
Then $[0,1]$ is the attractor of this IFS. This observation goes back to Mandelbrot~(see~\cite{M2}) and suggests how to define on the kind of fractals we are considering (and also on $[0,1]$) the van der Corput sequences.

Let $\psi=\{\psi_1,\ldots,\psi_m\}$ be our IFS and $F$ its attractor. Assume that $F_0$ is the initial set such that
$$
\psi_i(F_0)\subset F_0 \textrm{ for} \ i=1,\ldots,m.
$$
Fix a point $x_0\in F$ and apply $\psi_1, \ldots, \psi_m$ in this order to $x_0$ getting so the points $x_1,\ldots,x_m$. At the second step, we apply the $m$ mappings first to $x_1$, then to $x_2$ and so on, getting finally $m^2$ points ordered in a precise manner. Now we keep going, applying the functions of the IFS first to $x_1$, then to $x_2$ and continue so until we reach the point $x_{m^2}$, getting so $m^3$ points in the order determined by the construction. Iterating this procedure we get a sequence $(x_n)$ of points in $F$ which will be called the van der Corput sequence generated by $\psi$.

\begin{figure}[!h]
\begin{center}
\includegraphics[width=8cm]{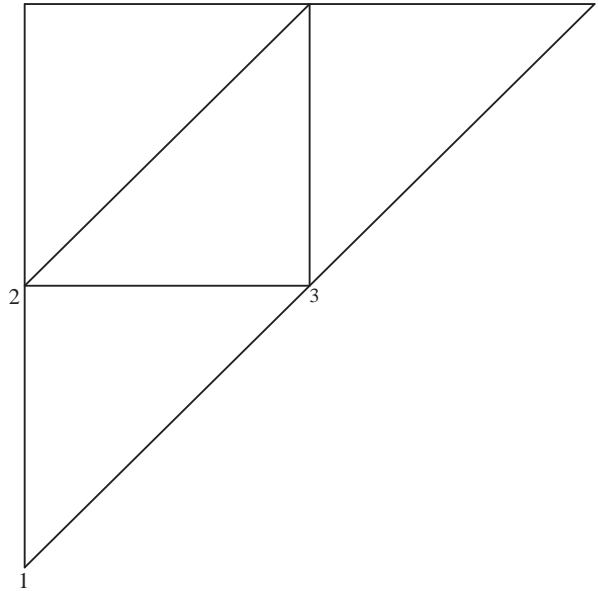}
\end{center} 
\caption[U.d.\! sequence of points on the Sierpi\' nski triangle (I step)]{Construction of a u.d.\! sequence of points on $T$ (I step)}
\label{AlgSierp1}
\end{figure}

\begin{figure}[!h]
\begin{center}
\includegraphics[width=8cm]{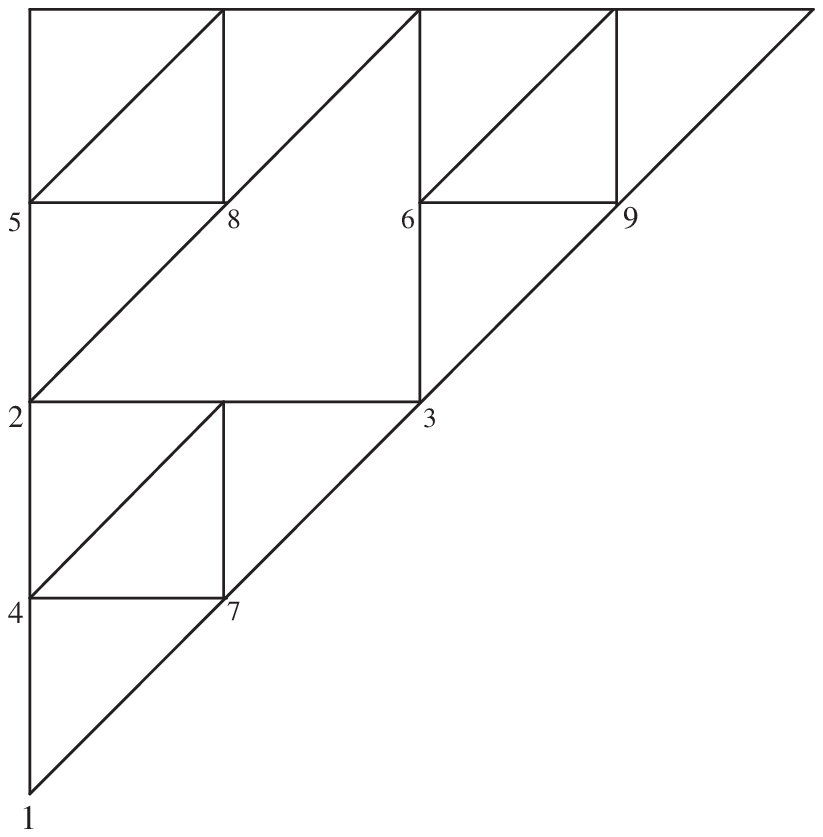}
\end{center}
\caption[U.d.\! sequence of points on the Sierpi\' nski triangle (II step)]{Construction of a u.d.\! sequence of points on $T$ (II step)}
\label{AlgSierp2}
\end{figure}

\begin{figure}[!ht]
\begin{center}
\includegraphics[width=8cm]{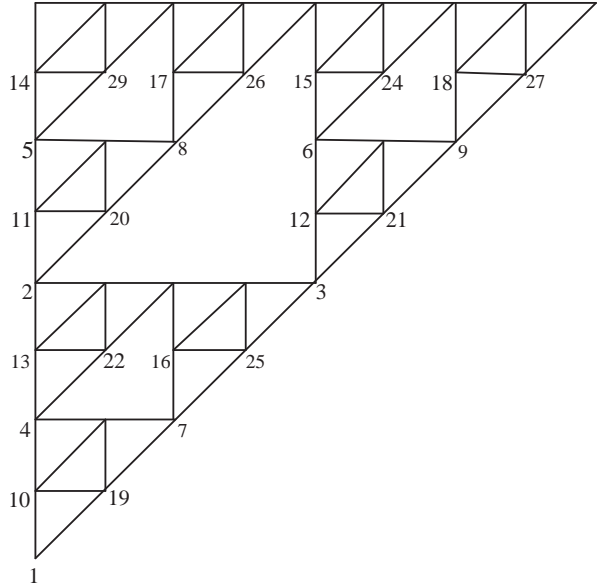}
\end{center}
\caption[U.d.\! sequence of points on the Sierpi\' nski triangle (III step)]{Construction of a u.d.\! sequence of points on $T$ (III step)}
\label{AlgSierp3}
\end{figure}

\paragraph{Example}\ \\
Consider the triangle $T_0\subset\mathbb{R}^2$ of vertices $(0,0)$, $(0,1)$ e $(1,1)$ and the similarities $S_1$, $S_2$ e $S_3$, defined as follows
\begin{equation*}
 \left\{ \begin{array}{ll}
S_1(x,y)=\bigg(\frac{1}{2}x,\frac{1}{2}y\bigg) & \textrm{}\\ \\
S_2(x,y)=\bigg(\frac{1}{2}x,\frac{1}{2}y+\frac{1}{2}\bigg) & \textrm{}\\ \\
S_3(x,y)=\bigg(\frac{1}{2}x+\frac{1}{2},\frac{1}{2}y+\frac{1}{2}\bigg) & \textrm{}
\end{array}\right.. 
\end{equation*} 
Let $T$ be the Sierpi\' nski triangle generated by this IFS starting from the initial set~$T_0$. 
Fixed $x_0=(0,0)$ by applying the algorithm we have
\begin{description}
\item[I step]
\begin{equation*}
\begin{array}{llll}
x_1=(0,0) & x_2=\big(0,\frac{1}{2}\big) & x_3=\big(\frac{1}{2},\frac{1}{2}\big) \textrm{}
\end{array}
\end{equation*} 

\item[II step]
\begin{equation*}
\begin{array}{llll}
x_1=(0,0) & x_2=\big(0,\frac{1}{2}\big) & x_3=\big(\frac{1}{2},\frac{1}{2}\big) \textrm{}\\ \\
x_4=\big(0,\frac{1}{4}\big) & x_5=\big(0,\frac{3}{4}\big) & x_6=\big(\frac{1}{2},\frac{3}{4}\big) \textrm{}\\ \\
x_{7}=\big(\frac{1}{4},\frac{1}{4}\big) & x_{8}=\big(\frac{1}{4},\frac{3}{4}\big) & x_{9}=\big(\frac{3}{4},\frac{3}{4}\big) \textrm{}
\end{array}
\end{equation*}

\item[III step]
\begin{equation*}
\begin{array}{llll}
x_{1}=(0,0) & x_{2}=\big(0,\frac{1}{2}\big) & x_{3}=\big(\frac{1}{2},\frac{1}{2}\big) \textrm{}\\ \\
x_{4}=\big(0,\frac{1}{4}\big) & x_{5}=\big(0,\frac{3}{4}\big) & x_{6}=\big(\frac{1}{2},\frac{3}{4}\big) \textrm{}\\ \\
x_{7}=\big(\frac{1}{4},\frac{1}{4}\big) & x_{8}=\big(\frac{1}{4},\frac{3}{4}\big) & x_{9}=\big(\frac{3}{4},\frac{3}{4}\big) \textrm{}\\ \\
x_{10}=\big(0,\frac{1}{8}\big) & x_{11}=\big(0,\frac{5}{8}\big) & x_{12}=\big(\frac{1}{2},\frac{5}{8}\big) \textrm{}\\ \\
x_{13}=\big(0,\frac{3}{8}\big) & x_{14}=\big(0,\frac{7}{8}\big) & x_{15}=\big(\frac{1}{2},\frac{7}{8}\big) \textrm{}\\ \\
x_{16}=\big(\frac{1}{4},\frac{3}{8}\big) & x_{17}=\big(\frac{1}{4},\frac{7}{8}\big) & x_{18}=\big(\frac{3}{4},\frac{7}{8}\big) \textrm{}\\ \\
x_{19}=\big(\frac{1}{8},\frac{1}{8}\big) & x_{20}=\big(\frac{1}{8},\frac{5}{8}\big) & x_{21}=\big(\frac{5}{8},\frac{5}{8}\big) \textrm{}\\ \\
x_{22}=\big(\frac{1}{8},\frac{3}{8}\big) & x_{23}=\big(\frac{1}{8},\frac{7}{8}\big) & x_{24}=\big(\frac{5}{8},\frac{7}{8}\big) \textrm{}\\ \\
x_{25}=\big(\frac{3}{8},\frac{3}{8}\big) & x_{26}=\big(\frac{3}{8},\frac{7}{8}\big) & x_{27}=\big(\frac{7}{8},\frac{7}{8}\big) \textrm{}
\end{array}
\end{equation*}
\end{description}
and so on. Figures \ref{AlgSierp1}, \ref{AlgSierp2} and \ref{AlgSierp3} show the first three steps of the algorithm.\\

Observe that if $[0,1]$ is seen as the attractor of the IFS described in (\ref{vdCIfs}) and if $x_0=0$, the sequence $(x_n)$ generated by this procedure is exactly the classical van der Corput sequence of base $m$ (see Subsection~\ref{sec:vdC}).
\newpage

\newpage

Let us come back to the general situation and show how a similar construction produces u.d.\! sequences of partitions.

Note that if we apply the $\psi_i$'s to $F$ in the same order as before, we construct a sequence $(\pi_k)$ of partitions of $F$ 
$$
\pi_k=\big\{ \psi_{j_k}\psi_{j_{k-1}}\cdots \ \psi_{j_1}(F) : {j_1},\ldots,{j_k}\in\{1,\ldots,m\}\big\}.
$$
Each of the $m^k$ sets $E_j^k$ of the partition $\pi_k$ contains exactly one point of the van der Corput sequence $(x_n)$ constructed above for $n=m^k$. We order the sets $E_j^k$ accordingly.

Let us denote by $\mathscr{E}_k$ the collection of the $m^k$ sets in $\pi_k$ and by
$\mathscr{E}$ the union of the families $\mathscr{E}_k$, for $k\in\mathbb{N}$. The sets of the class $\mathscr{E}$ are called \emph{elementary sets}.

As probability on $F$ we consider the normalized $s$-dimensional Hausdorff measure~$P$, i.e.
\begin{equation}\label{ProbFrac}
    P(A)=\frac{\mathcal{H}^s(A)}{\mathcal{H}^s(F)} \ \textrm{for any Borel set}\  A\subset F
\end{equation}
which is a regular probability (see Definition \ref{RegMis}).

The sequences of points and of partitions generated on $F$ by the algorithm are both u.d.\! with respect to the probability $P$. But before proving these results, we need to introduce some preliminary lemmas about the main properties of the family of the elementary sets.

\begin{lemma}\label{ElemContSets}\ \\
The elementary sets are $P$-continuity sets. 
\end{lemma}
\proof\ \\
Consider an elementary set $E_i=\psi_i(F)\in\mathscr{E}_1$. Let $x\in\partial E_i$. By definition, every neighbourhood $U$ of $x$ in the relative topology is such that $U\cap E_j\neq\emptyset$ for some $j\in\{1,2,\ldots,m\}$ and $j\neq i$. But each $E_j$ is closed, therefore $x\in E_j$. Hence $\partial E_i$ is contained in $\bigcup\limits_{\stackrel{j=1}{j\neq i}}^m(E_i\cap E_j).$
By Theorem \ref{imp}, we have
\begin{eqnarray}
0&\leq&\mathcal{H}^s(\partial E_i)\leq\mathcal{H}^s\left(\bigcup_{\stackrel{j=1}{j\neq i}}^m\big(E_i\cap E_j\big)\right)=\mathcal{H}^s\left(\bigcup_{\stackrel{j=1}{j\neq i}}^m\big(\psi_i(F)\cap \psi_j(F)\big)\right)\nonumber\\
&\leq&\sum_{\stackrel{j=1}{j\neq i}}^m\mathcal{H}^s\big(\psi_i(F)\cap \psi_j(F)\big)=0.\nonumber
\end{eqnarray}
Now, a generic elementary set $A\in\mathscr{E}_k$ with $k\geq 2$ is a homothetic image of an elementary set in $\mathscr{E}_1$ and therefore $\mathcal{H}^s(\partial A)=0$, too.\\
\endproof
\begin{lemma}\label{det-elementary}\ \\
The class $\mathscr{E}$ consisting of all elementary sets is determining.
\end{lemma}
\proof\
\\ Let $\mathcal{M}$ be the class consisting of all characteristic functions of sets $E\in\mathscr{E}$ and $f\in\mathcal{C}(F)$. By uniform continuity, for every $\varepsilon>0$ there exists $\delta>0$ such that $|f(x')-f(x'')|<\varepsilon$ whenever $\|x'-x''\|<\delta$. Choose $n\in\mathbb{N}$ such that every $E_k^n\in \mathscr{E}_n$ has diameter smaller than $\delta$. Take for any $E_k^n\in\mathscr{E}_n$ a point $t_k$ and consider the function
$$g(y)=\sum_{k=1}^{m^n}f(t_k)\chi_{E_k^n}(y),\quad y\in F.$$
For some $k$ we have $y\in E_k^n$ and so $|g(y)-f(y)|=|f(t_k)-f(y)|<\varepsilon$.\\
Hence, $\mbox{span}(\mathcal{M})$ is uniformly dense in $\mathcal{C}(F)$ and the conclusion follows by Theorem~\ref{teoSpan}.\\
\endproof

Now, we are ready to introduce the main results of this section.

\begin{theorem}\label{MainThm1}\ \\
The sequence $(\pi_n)$ of partitions of $F$ generated by the algorithm is u.d.\! with respect to the probability $P$.
\end{theorem}
\proof\ \\
By Lemma \ref{det-elementary}, we have to show that
$$
\lim_{n\rightarrow\infty}\frac{1}{m^n}\sum_{j=1}^{m^n}\chi_{E_h^k}(t_j^n)=\int_F {\chi_{E_h^k}(t)\ dP(t)}  
$$
for every $ E_h^k\in\pi_k$ and for every choice of $t_j^n\in E_j^n$. Let us fix $E_h^k\in\pi_k$. The previous relation is equivalent to 
\begin{equation*}
\lim_{n\rightarrow\infty}\frac{1}{m^n}\sum_{j=1}^{m^n}\chi_{E_h^k}(t_j^n)=\frac{1}{m^k}
\end{equation*}
because $$\int_F{\chi_{E_h^k}(t)\ dP(t)}=P({E_h^k})=c^{sk}P(F)=c^{sk}=\frac{1}{m^k}.$$
Now, observe that for $n>k$, among the $m^n$ sets generated by the algorithm, exactly one set of $\pi_n$ is contained in the fixed set ${E_h^k}$. Since there are $m^{n-k}$ sets of $\pi_n$ which are contained in ${E_h^k}$, then
$$
\lim_{n\rightarrow\infty}\frac{1}{m^n}\sum_{j=1}^{m^n}\chi_{E_h^k}(t_j^n)=\frac{m^{n-k}}{m^n}=\frac{1}{m^k}.
$$\endproof

\begin{theorem}\label{MainThm2}\ \\
The sequence $(x_i)$ of points of $F$ generated by the algorithm is u.d.\! with respect to $P$.
\end{theorem}
\proof\ \\
By Lemma \ref{det-elementary}, the class $\mathscr{E}$ is determining. Hence, for a fixed set $E\in\mathscr{E}_k$, we have to prove that
\begin{equation}\label{tesi}
\lim_{N\rightarrow\infty}\frac{1}{N}\sum_{i=1}^N\chi_{E }(x_i)=\int_{F}\chi_{E}\ dP=\frac{1}{m^k}.
\end{equation}
Let $m^t\leq N<m^{t+1}$, then
\begin{eqnarray}\label{rel0}
\frac{1}{N}\sum_{i=1}^N\chi_{E}(x_i) &=& \frac{1}{N}\sum_{i=1}^{m+m^2+\ldots+m^{t-1}}\chi_{E}(x_i)+
\frac{1}{N}\sum_{i={\frac{m^t-m}{m-1}}}^N\chi_{E}(x_i)\nonumber \\
&=&\frac{\big(\frac{m^t-m}{m-1}\big)}{N}\cdot\frac{1}{\big(\frac{m^t-m}{m-1}\big)}
\sum_{i=1}^{m+m^2+\ldots+m^{t-1}}\chi_{E}(x_i)\nonumber\\&+&\frac{N-\big(\frac{m^t-m}{m-1}\big)}{N}
\cdot\frac{1}{N-\big(\frac{m^t-m}{m-1}\big)}\sum_{i={\frac{m^t-m}{m-1}}}^N\chi_{E}(x_i)
\end{eqnarray}
since $1+m+m^2+\ldots+m^{t-1}=\frac{m^t-1}{m-1}$.\\
Observe that for $i>\frac{m^t-1}{m-1}$, because of the order of the points $x_i$, among the first $m^k$ points exactly one point of the sequence $(x_i)$ is contained in the fixed set $E$. Hence, for $t\rightarrow\infty$ we have
\begin{equation}\label{rel1}
\frac{1}{\big(\frac{m^t-m}{m-1}\big)}\sum_{i=1}^{m+m^2+\ldots+m^{t-1}}\chi_{E}(x_i)\rightarrow\frac{1}{m^k}.
\end{equation}
Writing $N$ as $N=\big(\frac{m^t-m}{m-1}\big)+Mm^k+r$ with $0\leq r<m^k$, we have
\begin{eqnarray}\label{rel2}
\frac{1}{N-\big(\frac{m^t-m}{m-1}\big)}\sum_{i=\frac{m^t-m}{m-1}}^N\chi_{E}(x_i)
&=&\frac{Mm^k}{{N-\big(\frac{m^t-m}{m-1}\big)}}\cdotp
\frac{1}{Mm^k}\sum_{i=\frac{m^t-m}{m-1}}^{N-r}\chi_{E}(x_i)\nonumber\\&+&
\frac{r}{{N-\big(\frac{m^t-m}{m-1}\big)}}\cdotp\frac{1}{r}\sum_{i={N-r+1}}^{N}\chi_{E}(x_i).
\end{eqnarray} 
By the previous remarks we get that
$$
\frac{1}{Mm^k}\sum_{i=\frac{m^t-m}{m-1}}^{N-r}\chi_{E}(x_i)=\frac{1}{m^k},
$$
while for $N\rightarrow\infty$ and hence for $t\rightarrow\infty$ we have
$$
\frac{r}{{N-\big(\frac{m^t-m}{m-1}\big)}}\sum_{i={N-r+1}}^{N}\chi_{E}(x_i)\rightarrow 0
$$
because $0\leq\frac{r}{N}<\frac{m^k}{m^t}$.\\
Using the last two relations in (\ref{rel2}) and taking the limit for $N\rightarrow\infty$ (and hence for $t\rightarrow\infty$) we have
\begin{equation}\label{rel3}
\frac{1}{N-\big(\frac{m^t-m}{m-1}\big)}\sum_{i=1}^N\chi_{E}(x_i)\rightarrow\frac{1}{m^k}.
\end{equation}
Finally, (\ref{rel0}) is a convex combination of two terms which both tend to $\frac{1}{m^k}$ for $N\rightarrow\infty$ because of (\ref{rel1}) and (\ref{rel3}). Therefore, the conclusion (\ref{tesi}) holds. \\
\endproof

\subsection{Order of convergence of the elementary discrepancy}\label{subsec: ElemDiscr}
In the following, we will provide an estimate for the elementary discrepancy of the sequences of van der Corput type generated by our explicit algorithm. 

Note that Lemma \ref{det-elementary} implies the family $\mathscr{E}$ of all elementary sets of $F$ is a discrepancy system (see Definition \ref{DisSystem}). So, according to Definition \ref{Gdiscrepancy}, the \emph{elementary discrepancy} of a sequence $\omega$ of points in $F$ is given by
\begin{equation*}
D_N^{\mathscr{E}}(\omega)=\sup_{E\in\mathscr{E}}\Bigg|\frac{1}{N}\sum_{i=1}^N\chi_E(x_i)-P(E)\Bigg|.
\end{equation*}
 
The choice of the elementary discrepancy is convenient because the family of elementary sets is obtained in the most natural way by the construction of the fractal $F$ and because the elementary sets can be constructed for every IFS fractal regardless of the complexity of its geometric structure.

In the papers \cite{Grab-Tichy}, \cite{CT} and \cite{CT2} the authors also gave estimates for the elementary discrepancy of the van der Corput type sequences produced on the particular fractals considered, finding that is of the order $\mathcal{O}\big(\frac{1}{N}\big)$. 
Our results include theirs, giving a more trasparent proof and taking in consideration the whole class of fractals described in the previous subsection. 

The next theorem evaluates the elementary discrepancy for the sequences of points generated by our algorithm.

\begin{theorem}\label{Thm: ElemDiscr}\ \\
Let $(x_i)$ be the sequence of points generated on $F$ by the algorithm described in the previous subsection and let $N\geq 1$. Then for the elementary discrepancy we have
$$
   \lim_{N\to\infty}ND_{N}^{\mathscr{E}}(\omega_N)=1.
$$
where $\omega_N=(x_1,\ldots, x_N)$.
\end{theorem}
\proof \ \\
The lower bound is trivial. In fact, for any $k\in\mathbb{N}$ we have
$$D_N^{\mathscr{E}}(\omega_N)\geq\frac{1}{N}-\frac{1}{m^k}.$$
In order to find an upper bound for $D_N^{\mathscr{E}}(\omega_N)$, let us consider $D_N^{\mathscr{E}_k}(\omega_N)$ for any $k\in\mathbb{N}$: 
$$
D_N^{\mathscr{E}_k}(\omega_N)=
\sup_{E\in\mathscr{E}_k}\Bigg|\frac{1}{N}\sum_{i=1}^N\chi_{E}(x_i)-\frac{1}{m^k}\Bigg|.
$$
Fix $k\in\mathbb{N}$ and let $E\in\mathscr{E}_k$. Among the first $m^k$ points of the sequence $(x_i)$ exactly one point is contained in the fixed set $E$ because of the special order induced by the algorithm.
 
Let us distinguish two different cases:
\begin{enumerate}
\item For $N\leq m^k$, the set $E$ contains at most one point of $\omega_N$. Hence
\begin{equation*}
\Bigg|\frac{1}{N}\sum_{i=1}^N\chi_{E}(x_i)-\frac{1}{m^k}\Bigg|=\max\Bigg\{\bigg|\frac{1}{N}-\frac{1}{m^k}\bigg|,
\bigg|0-\frac{1}{m^k}\bigg|\Bigg\}\leq\frac{1}{N}.
\end{equation*}
\item If $N> m^k$, we can write $N$ as follows
$$N=Q\cdot m^k+r\ \textrm{with}\ \ 0\leq r < m^k \ \textrm{and}\ Q\geq 1.$$
Therefore, every $E\in\mathscr{E}_k$ contains either $Q$ points or $Q+1$ points and hence
\begin{equation*}
\Bigg|\frac{1}{N}\sum_{i=1}^N\chi_{E}(x_i)-\frac{1}{m^k}\Bigg|\leq\max\Bigg\{\bigg|\frac{Q}{N}-\frac{1}{m^k}\bigg|,
\bigg|\frac{Q+1}{N}-\frac{1}{m^k}\bigg|\Bigg\}.
\end{equation*}

Note that 
\begin{equation*}
    \bigg|\frac{Q}{N}-\frac{1}{m^k}\bigg|=\bigg|\frac{Qm^k-N}{Nm^k}\bigg|
                                       =\bigg|\frac{-r}{Nm^k}\bigg|
                                       <\frac{m^k}{Nm^k}
                                       =\frac{1}{N},              
\end{equation*}
while
\begin{equation*}
    \bigg|\frac{Q+1}{N}-\frac{1}{m^k}\bigg|
                                =\bigg|\frac{Qm^k+m^k-N}{Nm^k}\bigg|
                            = \bigg|\frac{m^k-r}{Nm^k}\bigg|
                             = \bigg|\frac{1}{N}-\frac{r}{Nm^k}\bigg|
                             < \frac{1}{N}.        
\end{equation*}
So we have that 
$$\Bigg|\frac{1}{N}\sum_{i=1}^N\chi_{E}(x_i)-\frac{1}{m^k}\Bigg|<\frac{1}{N}.$$
\end{enumerate}
It follows that for any $k\in\mathbb{N}$ we have $
    D_{N}^{\mathscr{E}_k}(\omega_N)<\frac{1}{N}.
$
This implies that  
$
    D_{N}^{\mathscr{E}}(\omega_N)\leq\frac{1}{N},
$ as we wanted to prove.\\
\endproof

Note that $D_{N}^{\mathscr{E}}(\omega_N)$ is equal to zero for infinitely many $N$ and precisely when $N=\sum\limits_{i=1}^nm^i$ for any $n\in\mathbb{N}$\ . This is due to the fact that the elementary discrepancy of the sequence of partitions $(\pi_n)$ generated by the algorithm is exactly zero.

%% file: Cap3.tex
\chapter{Generalized Kakutani's sequences of partitions}
A first extension of the results introduced in Section \ref{UDFractals} can be obtained by eliminating the restriction that the similarities defining the fractal have all the same ratio. The study of this problem on the simplest case of attractor of an IFS, namely $[0,1]$, has taken us to consider Kakutani's sequences of partitions and their recent generalization: the $\rho-$refinements \cite{Vol}. In this chapter, we firstly introduce the technique of successive $\rho-$refinements which generalizes Kakutani's splitting procedure and allows to construct new families of u.d.\! sequences of partitions. Successively, we analyze the problem of deriving bounds for the discrepancy of these sequences. The approach that we use is based on a tree representation of any sequence of partitions constructed by successive $\rho-$refinements, which is precisely the parsing tree generated by Khodak's coding algorithm. Finally, with the help of this technique, we present an application of these results to a class of fractals which includes the one considered in Section \ref{UDFractals}.
\section{A generalization of Kakutani's splitting procedure}
In Subsection \ref{KakProc} we introduced Kakutani's splitting procedure, which works through successive $\alpha-$refinements of the unit interval. In a recent paper \cite{Vol}, this concept has been generalized through the new notion of $\rho-$refinement and it has been proved that the sequence of partitions generated by successive $\rho-$refinements of the trivial partiton is u.d..
\subsection{$\rho-$refinements}
Consider any non-trivial finite partition $\rho$ of $[0,1]$ and from now on we keep it fixed. 
\begin{definition}\ \\
Let $\pi$ be any partition of $[0,1]$. The \emph{$\rho$-refinement} of $\pi$ (which will be denoted by $\rho\pi$) is obtained by splitting all the intervals of $\pi$ having maximal lenght into a finite number of parts positively homothetically to $\rho$.
\end{definition}
Note that, if $\rho=\left\{[0,\alpha], [\alpha,1]\right\}$ then the $\rho-$refinement is just Kakutani's $\alpha-$refine\-ment. As in Kakutani's case, we can iterate the splitting procedure. We will denote by $\rho^2\pi$ the $\rho$-refinement of $\rho\pi$ and, in general, by $\rho^n\pi$ the $\rho$-refinement of $\rho^{n-1}\pi$.

In the following we will consider the sequence $(\rho^n\omega)$, where $\omega$ is the trivial partition of $[0,1]$, and we will prove that $(\rho^n\omega)$ is u.d..
 
\begin{rem}\ \\
It is important to note that in general $(\rho^n\pi)$ is not u.d.\! for any partition $\pi$. For instance, let $\pi=\left\{\left[0, \frac 25\right], \left[\frac 25, 1\right]\right\}$ and $\rho=\left\{\left[0, \frac 12\right], \left[\frac 12, 1\right]\right\}$. It is clear that the $\rho-$refinement operates alternatively on $\left[\frac 25, 1\right]$ and $\left[0, \frac 25\right]$. So, if we consider the sequence of measures $(\nu_n)$ associated to $(\rho^n\pi)$ (see (\ref{AssMeas}) for the definition), then the subsequence $(\nu_{2n})$ converges to $\mu_1$ while the subsequence $(\nu_{2n+1})$ converges to $\mu_2$ where
$$\mu_1(E)=\frac 54\cdot\lambda\left(E\cap\left[0, \frac 25\right]\right)
+\frac 56\cdot\lambda\left(E\cap\left[\frac 25, 1\right]\right) $$

$$\mu_2(E)=\frac 56\cdot\lambda\left(E\cap\left[0, \frac 25\right]\right)
+\frac {10}{9}\cdot\lambda\left(E\cap\left[\frac 25, 1\right]\right)$$
for any measurable set $E\subset[0,1]$. Hence, $(\nu_n)$ does not converge and consequently $(\rho^n\pi)$ is not u.d..

We can find the problem showed by this example also in the simplest case of Kakutani's splitting procedure. So, it could be interesting to find significant sufficient conditions on $\pi$ in order to guarantee the uniform distribution of $(\alpha^n\pi)$ or more in general of $(\rho^n\pi)$.
\end{rem}

Before introducing the analogous of Kakutani's theorem for these new sequences of partitions (Theorem \ref{VolThm}), let us fix some notations and recall some preliminary results. 
\\

Firstly, we need some classical definitions from ergodic theory (see \cite[Chapter~29]{P} due to F.Blume).

\begin{definition}\ \\
A measurable function $\varphi:[0,1]\rightarrow [0,1]$ is said to be \emph{measure preserving} if
\begin{itemize}
\item $\varphi$ is bijective,
\item $\varphi(A), \varphi^{-1}(A)$ are measurable when $A$ is measurable,
\item $\lambda(\varphi^{-1}(A))=\lambda(A)$.
\end{itemize}
A countable family $\mathscr{F}$ of measurable functions is said to be measure preserving if any $\varphi\in\mathscr{F}$ is measure preserving.
\end{definition}

\begin{definition}\ \\
Given a measurable function $\varphi:[0,1]\rightarrow [0,1]$, a measurable set $A$ is called \emph{$\varphi-$invariant} if
$$\lambda(\varphi^{-1}(A)\ \Delta\  A)=0,$$
where $\Delta$ is the symmetric difference. If $\mathscr{F}$ is a countable family of measurable functions, a measurable set $A$ is said to be $\mathscr{F}-$invariant if it is $\varphi-$invariant for any $\varphi\in\mathscr{F}$.
\end{definition}

\begin{definition}\ \\
A measurable function $\varphi:[0,1]\rightarrow [0,1]$ is called \emph{ergodic} if it is measure preserving and if for each set $A$ $\varphi$-invariant we have
$$\lambda(A)=0 \ \ \text{or} \ \ \lambda(A)=1.$$
A countable family $\mathscr{F}$ of measurable functions is said to be ergodic if any $\varphi\in\mathscr{F}$ is ergodic.
\end{definition}

\begin{definition}\ \\
Given a measurable function $\varphi:[0,1]\rightarrow [0,1]$, a real-valued function $f$ on $[0,1]$ is said to be \emph{$\varphi$-invariant} if
$$f(\varphi(x))=f(x)$$
holds $\lambda-$almost everywhere. If $\mathscr{F}$ is a countable family of measurable functions, $f$ is called $\mathscr{F}-$invariant if it is $\varphi$-invariant for any $\varphi\in\mathscr{F}$.
\end{definition}

It is easy to prove that
\begin{theorem}\label{ThmErgod}\ \\
Let $\mathscr{F}$ be a countable ergodic family and suppose that $f$ is a measurable real-valued function on $[0,1]$. If $f$ is $\mathscr{F}-$invariant, then $f$ is constant almost everywhere. 
\end{theorem}
\proof\ \\
Consider for every $\alpha\in\mathbb{R}$ the set
$$A_{\alpha}=\left\{x\in [0,1]: f(x)>\alpha\right\}.$$
$A_{\alpha}$ is measurable, since the function $f$ is measurable. Moreover, $A_{\alpha}$ is $\mathscr{F}-$invariant because $f$ is $\mathscr{F}-$invariant. Therefore we have
$$\lambda(A_{\alpha})=0 \ \ \text{or}\ \ \lambda(A_{\alpha})=1\ \ \text{for all}\ \alpha\in\mathbb{R}\,,$$
since $\mathscr{F}$ is an ergodic family. It follows that $f$ is constant almost everywhere. In fact if not, there would exist an $\alpha\in\mathbb{R}$ such that $0<\lambda(A_{\alpha})<1$ which would contradict the ergodicity of~$\mathscr{F}$.\\
\endproof

Let us recall an important theorem due to Hewitt and Savage \cite[Theorem~11.3]{HS}.
\begin{theorem}\label{HS}\ \\
A product measure on an infinite product of measure spaces can assume only the values $0$ and $1$ for sets which are invariant under all finite permutations of the coordinates.
\end{theorem}

Now, we can come back to the $\rho-$refinements and introduce some concepts and properties necessary to the proof of Theorem \ref{VolThm}.

Let $\rho=\left\{[u_{i-1}, u_i] : 1\leq i\leq k \right\}$ be the fixed partition of $[0,1]$ and let us denote by $\alpha_i=u_i-u_{i-1}$ for $1\leq i\leq k$, the lenghts of the $k$ intervals of $\rho$. Let $[\rho]^n$ be the $n-$th \emph{$\rho-$adic partition} of $[0,1]$, obtained from $[\rho]^{n-1}$ (where $[\rho]^1=\rho$) by subdividing all its $k^{n-1}$ intervals positively homothetically to $\rho$. If an interval belongs to $[\rho]^n$, we will say that it has rank $n$. 

The $k$ intervals of $\rho$ will be denoted by
$$I(\alpha_i)=\left[\sum\limits_{h=1}^{i-1}\alpha_h, \sum\limits_{h=1}^{i}\alpha_h\right]=[u_{i-1}, u_i].$$
If $[y_{j-1}, y_j]=I(\alpha_{i_1}\alpha_{i_2}\ldots\alpha_{i_{n-1}})$ with $1\leq j\leq k^{n-1}$ is a generic interval of rank $n-1$, then its subintervals belonging to $[\rho]^n$ are 
$$I(\alpha_{i_1}\alpha_{i_2}\ldots\alpha_{i_{n-1}}\alpha_{i_n})
=\left[y_{j-1}+(y_j-y_{j-1})\sum_{h=1}^{i_n-1}\alpha_h,\  y_{j-1}+(y_j-y_{j-1})\sum_{h=1}^{i_n}\alpha_h\right] $$
for $\alpha_{i_n}=\alpha_1,\ldots, \alpha_k$. Moreover, by varying $[y_{j-1}, y_j]\in[\rho]^{n-1}$, $1\leq j\leq k^{n-1}$, we obtain all the $k^n$ intervals of $[\rho]^n$. Note that 
\begin{equation}\label{measInt}
\lambda(I(\alpha_{i_1}\alpha_{i_2}\ldots\alpha_{i_{n-1}}\alpha_{i_n}))=\prod_{m=1}^n\alpha_{i_m}.
\end{equation}

\paragraph{Example}\ \\
Let $\rho=\left\{\left[0,\frac 14\right], \left[\frac 14, \frac 12\right], \left[\frac 12, 1\right]\right\}$. In this case we have $\alpha_1=\alpha_2=\frac 14$ and $\alpha_3=\frac 12$. We only want to construct $[\rho]^n$ for $n=1,2$. Then

$$[\rho]^1=\rho=\Bigg\{\underbrace{\left[0,\frac 14\right]}_{I(\alpha_1)}, \underbrace{\left[\frac 14, \frac 12\right]}_{I(\alpha_2)}, \underbrace{\left[\frac 12, 1\right]}_{I(\alpha_3)}\Bigg\}.$$
Now, $[\rho]^2$ can be obtained from $[\rho]^1$ by splitting all its intervals homothetically to $\rho$. Practically, we have to subdivide each interval in two equal parts, then take the first of these two halves and split it again in two equal parts. So we have
$$[\rho]^2=\Bigg\{\underbrace{\left[0,\frac 1{16}\right]}_{I(\alpha_1\alpha_1)}, \underbrace{\left[\frac 1{16}, \frac 18\right]}_{I(\alpha_1\alpha_2)}, \underbrace{\left[\frac 18, \frac 14\right]}_{I(\alpha_1\alpha_3)}, \underbrace{\left[\frac 14, \frac 5{16}\right]}_{I(\alpha_2\alpha_1)}, \underbrace{\left[\frac 5{16}, \frac 38\right]}_{I(\alpha_2\alpha_2)}, \underbrace{\left[\frac 38, \frac 12\right]}_{I(\alpha_2\alpha_3)}, \underbrace{\left[\frac 12, \frac 58\right]}_{I(\alpha_3\alpha_1)}, \underbrace{\left[\frac 58, \frac 34\right]}_{I(\alpha_3\alpha_2)}, \underbrace{\left[\frac 34, 1\right]}_{I(\alpha_3\alpha_3)} \Bigg\}.$$
\\

Let $X=\{\alpha_1,\ldots,\alpha_k\}$ and let $\sigma$ be the probability on $X$ such that $\sigma(\{\alpha_i\})=\alpha_i$ for $1\leq i\leq k$. Put $X_m=X$ and $\sigma_m=\sigma$ for any $m\in\mathbb{N}$. Denote by $$Y=\prod_{m=1}^\infty X_m$$ and consider on $Y$ the usual product probability $\mu$. 

If $C=C(\alpha_{i_1}\alpha_{i_2}\ldots\alpha_{i_n})=\prod\limits_{m=1}^\infty X'_m$, where $X'_m=\{\alpha_{i_m}\}$ for $m\leq n$ and $X'_m=X$ for $m>n$, is a cylinder set then
\begin{eqnarray}\label{measCyl}
\mu(C)&=&
\mu\left(\{\alpha_{i_1}\}\times\{\alpha_{i_2}\}\times\cdots\times\{\alpha_{i_n}\}\times X'_{n+1}\times X'_{n+2}\times\cdots\right)\nonumber\\
&=&\sigma(\{\alpha_{i_1}\})\cdot\sigma(\{\alpha_{i_2}\})\cdots\sigma(\{\alpha_{i_n}\})\cdot \sigma(X)\cdot \sigma(X)\cdots\nonumber\\
&=&\prod_{m=1}^n\alpha_{i_m}.
\end{eqnarray}
To every point $t\in[0,1]$ we can associate a sequence $(\alpha_{i_m})$ such that
$$t\in\bigcap_{m=1}^\infty I(\alpha_{i_1}\ldots\alpha_{i_m}),$$
that is
$$[t]_{\rho}=\alpha_{i_1}\ldots\alpha_{i_m}\ldots.$$
We called $[t]_{\rho}$ the $\rho-$adic representation of $t$. 

It is important to take care of an expected ambiguity of this representation. In fact, there are two such sequences $(\alpha_{i_m})$ associated to a $t$ in the countable set of points belonging to the endpoints of some $[\rho]^n$. In this case to solve the problem, we associate to $t$ the sequence for which definitively $\alpha_{i_m}=\alpha_1$. 

This defines a $1-1$ mapping $\phi:[0,1]\to Y'$, where $Y'$ is a subset of $Y$ obtained by removing from $Y$ the countable set of sequences $(\alpha_{i_m})$ such that definitively $\alpha_{i_m}=\alpha_{k}$, i.e.
$$Y'=Y\setminus\left\{(\alpha_{i_m}): \alpha_{i_m}=\alpha_{k}\quad \forall m\geq m_0\right\} .$$

Note that $\mu(Y\setminus Y')=0$. In fact, if $y\in Y\setminus Y'$ then $y$ is of the form
$$y=(\alpha_{i_1}\ldots\alpha_{i_{m_0-1}}\alpha_{k}\alpha_{k}\alpha_{k}\ldots)$$
so
$$Y\setminus Y'\subset C(\alpha_{i_1}\ldots\alpha_{i_{m_0-1}}\underbrace{\alpha_{k}\ldots\alpha_{k}}_{n\  \text{times}}).$$
Therefore we have that
$$0\leq\mu(Y\setminus Y')\leq\mu\left(C(\alpha_{i_1}\ldots\alpha_{i_{m_0-1}}\underbrace{\alpha_{k}\ldots\alpha_{k}}_{n\  \text{times}})\right)=\alpha_{i_1}\cdots\alpha_{i_{m_0-1}}\alpha_{k}^n\to 0 $$
as $n\to\infty$, since $\alpha_{k}<1$.
 
Moreover, $\phi$ is a measure preserving mapping if we take on $[0,1]$ the Lebesgue measure $\lambda$ and on $Y'$ the restriction of $\mu$. This follows immediately noting that $\rho-$adic intervals and cylinder sets with the same indices have the same measure (see (\ref{measInt}) and (\ref{measCyl})). Hence, $\phi$ is a measure isomorphism between $\left([0,1], \lambda|_{[0,1]}\right)$ and  $\left(Y', \mu|_{Y'}\right)$.

Let $I$ and $J$ be two disjoint subintervals of $[0,1]$ having the same lenght and let $J=I+c$ with $0<c<1$. Let us define the following function
\begin{equation}\label{funzIJ}
f_{I,J}(x)=\left\{
\begin{array}{ll}
x+c & \textrm{if $x\in I$} \\
x-c & \textrm{if $x\in J$} \\
x & \textrm{otherwise}  
\end{array} 
\right..
\end{equation}
It is important to observe that $f_{I,J}$ is measure preserving. Let us denote by $\mathcal{F}$ the family of all functions $f_{I,J}$ such that $I$ and $J$ are two $\rho-$adic intervals having the same lenght. The intervals $I$ and $J$ do not need to have necessarily the same rank.

\begin{lemma}\label{lem1}\ \\
The family $\mathcal{F}$ is ergodic.
\end{lemma}
\proof\ \\
Let us denote by $\mathcal{F}'$ the family of transformations on $Y'$ correspondent to $\mathcal{F}$ by using the isomorphism described above, i.e 
$$\mathcal{F}'=\{f':Y'\to Y' \ \text{s.t.}\  f'=\phi f\phi^{-1}\ \text{for some}\  f\in\mathcal{F}\}.$$ 
When $f=f_{I,J}$ and $I$ and $J$ have the same rank, the correspondent function $f'$ on~$Y'$ is a permutation of a finite number of coordinates and it preserves the product measure $\mu$. Let us denote by $\mathcal{G}'$ the family of such functions. By Theorem \ref{HS} the family $\mathcal{G}'$ is ergodic and consequently $\mathcal{F}'$ is also ergodic, since $\mathcal{G}'\subset\mathcal{F}'$. In conclusion, since the isomorphism $\phi$ is measure preserving the family $\mathcal{F}$ results to be ergodic, too.\\
\endproof

For the partition $\rho^n\omega$, let $A_n$ be the lenght of the longest interval and $a_n$ the lenght of the shortest interval. Moreover, let us denote by $k(n)$ the number of intervals of the partition $\rho^n\omega$. We have the following results.
\begin{lemma}\label{lem2}\ 
\begin{enumerate}
\item For any $n\in\mathbb{N}$ we have $a_1A_n\leq a_n.$
\item If $\pi_n=\rho^n\omega$ then $\lim\limits_{n\to\infty} diam(\pi_n)=0.$
\end{enumerate}
\end{lemma}
\proof\ 
\begin{enumerate}
\item Since $A_n<1$ for any $n\in\mathbb{N}$, the strict inequality holds for $n=1$. Now, proceed by induction. Suppose that $a_1A_{n-1}\leq a_{n-1}$ holds. There are two possibilities either $a_n=a_{n-1}$ or $a_n<a_{n-1}$. In the first case, since $A_n<A_{n-1}$ we have 
$$a_1A_n<a_1A_{n-1}\leq a_{n-1}=a_n.$$
In the second case, the shortest interval of the partition $\rho^n\omega$ is obtained by splitting the longest interval of $\rho^{n-1}\omega$, so $a_n=a_1A_{n-1}$. Hence, we have
$$a_1A_n<a_1A_{n-1}=a_n.$$

\item According to the notation introduced above, obviously $a_n<\frac 1n$ for any $n\in\mathbb{N}$. By applying the relation just proved we have
$$diam(\pi_n)=A_n\leq\frac{a_n}{a_1}<\frac{1}{a_1n}$$
and so the conclusion follows.\\ \endproof
\end{enumerate}

\begin{lemma}\label{lem3}\ \\
The family of the characteristic functions of all intervals belonging to the partitions $\rho^n\omega$ for $n\in\mathbb{N}$ is determining.
\end{lemma}
\proof\ \\
Let $\mathcal{M}$ be the class consisting of the characteristic functions of all intervals belonging to the partitions $\rho^n\omega$ for $n\in\mathbb{N}$ and let $f\in\mathcal{C}([0,1])$. 

By uniform continuity we have that  for any $\varepsilon>0$ there exists $\overline{\delta}>0$ such that $|f(x')-f(x'')|<\varepsilon$ whenever $|x'-x''|<\overline{\delta}$. Moreover, by the second part of Lemma~\ref{lem2}, we have that for any $\delta>0$ there exists $\overline{n}\in\mathbb{N}$ such that for every $n>\overline{n}$ we have $diam(\rho^n\omega)<\delta$. Hence for every $n>\overline{n}$ we have that each $E_j^n\in \rho^n\omega$ has diameter smaller than $\delta$. 

So we can choose $n\in\mathbb{N}$ such that each $E_j^n\in \rho^n\omega$ has diameter smaller than $\overline{\delta}$. Fixed a such $n$, take a point $t_j$ in any $E_j^n\in\rho^n\omega$ and consider the function
$$g(y)=\sum_{j=1}^{k(n)}f(t_j)\chi_{E_j^n}(y),\quad y\in [0,1].$$
Now, $y\in E_j^n$ for some $j\in\{1,\ldots,k(n)\}$ then $|g(y)-f(y)|=|f(t_j)-f(y)|<\varepsilon$. Hence, $\mbox{span}(\mathcal{M})$ is uniformly dense in $\mathcal{C}([0,1])$ with respect to the $||\cdot||_\infty$ and so the conclusion follows by Theorem \ref{teoSpan}.\\
\endproof

\subsection{A generalization of Kakutani's Theorem}
\begin{theorem}\label{VolThm}\ \\
The sequence $(\rho^n\omega)$ is u.d..
\end{theorem}
\proof\ \\
Let us denote by $(\nu_n)$ the sequence of measures associated to the sequence of partitions $(\rho^n\omega)$ defined accordingly to (\ref{AssMeas}). We have to prove that $(\nu_n)$  is weakly convergent to $\lambda$.

It is well known that the set of all Borel probability measures on $[0,1]$, with the topology associated to the weak convergence, is metrizable and compact (see \cite[Theorem 6.4]{Par}).
Then $(\nu_n)$ admits weakly convergent subsequences. So all we need to prove is that any such subsequence converges to $\lambda$.

First of all, let us prove that the family of the characteristic functions of all $\rho-$adic intervals is determining. In fact, every interval $J\in\rho^n\omega$ belongs to some $[\rho]^m\omega$. The viceversa is also true, namely every interval $I\in[\rho]^m\omega$ sooner or later belongs to some $\rho^n\omega$. This is due to the fact that
$$s=\sup\left\{r: J_r\in\rho^r\omega,\ I\subset J_r\right\}$$
is well defined, since the $\rho-$adic intervals are either disjoint or contained one in the other. Moreover, by the second part of Lemma \ref{lem2}, the diameter of $\rho^n\omega$ tends to zero and so $I=J_s$. Therefore, by Lemma \ref{lem3} we can conclude that the family of the characteristic functions of all $\rho-$adic intervals is determining, too.

Consequently, it is sufficient to prove that for each weakly convergent subsequence $(\nu_{n_k})$ we have
$$\nu_{n_k}(J)\rightharpoonup \lambda(J),\ \ \forall J\in[\rho]^m\omega, \ \forall m\in\mathbb{N}$$
where we denote by ``$\rightharpoonup$'' the weak convergence.

Let $J$ be any $\rho-$adic interval and suppose that $m\in\mathbb{N}$ is such that for any $n\geq m$ every $\rho^n\omega$ subdivides $J$. If $J$ is splitted in $k$ intervals by $\rho^n\omega$, then
$$
    ka_n\leq\lambda(J)\leq kA_n
$$
where $a_n$ and $A_n$ are the quantities considered in Lemma \ref{lem2}.

If $k(n)$ is the number of intervals in $\rho^n\omega$, we have that for all $n\geq m$
\begin{equation}\label{Irel}
    \frac{\lambda(J)}{k(n)A_n}\leq\nu_n(J)\leq\frac{\lambda(J)}{k(n)a_n}.
\end{equation}

In fact, let us denote by $K(J)$ the number of intervals in which $J$ is subdivided by $\rho^n\omega$, i.e
$$K(J):=\sum_{i=1}^{k(n)}\delta_{t_i^n}(J),$$
where the points $t_i^n$ are the points determining the partition $\rho^n\omega$.
\\So $\nu_n(J)=\frac{1}{k(n)}\sum\limits_{i=1}^{k(n)}\delta_{t_i^n}(J)=\frac{K(J)}{k(n)}$ and therefore
\begin{itemize}
\item $\lambda(J)\geq K(J) a_n=\nu_n(J) k(n) a_n$ $\Rightarrow$ $\nu_n(J)\leq\frac{\lambda(J)}{k(n)a_n},$
\item $\lambda(J)\leq K(J) A_n=\nu_n(J) k(n) A_n$ $\Rightarrow$ $\nu_n(J)\geq\frac{\lambda(J)}{k(n)A_n}$.
 \end{itemize}

\noindent By (\ref{Irel}) and by the first part of Lemma \ref{lem2} we have
\begin{equation}\label{IIrel}
a_1\lambda(J)\leq\frac{a_1\lambda(J)}{k(n)a_n}\leq\nu_n(J)
\leq\frac{\lambda(J)}{k(n)a_1A_n}\leq\frac{\lambda(J)}{a_1}.
\end{equation} 

Now, suppose that $(\nu_{n_k})$ is a subsequence weakly convergent to $\nu$. Then, by~(\ref{IIrel}) we have that for any $\rho-$adic interval $J$ the following holds
\begin{equation}\label{IIIrel}
    a_1\lambda(J)\leq\nu(J)\leq\frac{\lambda(J)}{a_1}.
\end{equation}
Since the family of the characteristic functions of all $\rho-$adic intervals is determining, the relation (\ref{IIIrel}) holds for any Borel set $B$ in $[0,1]$. Therefore $\lambda\ll\nu\ll\lambda$ and if we denote by $\frac{d\nu}{d\lambda}$ the Radon-Nikodym derivative of $\nu$ with respect to $\lambda$, then we have that
\begin{equation*}
    a_1\leq\frac{d\nu}{d\lambda}\leq\frac{1}{a_1}.
\end{equation*}
Note that if $I$ and $J$ are two intervals having the same lenght and belonging to some $\rho^n\omega$ (not necessarily having the same rank), then the splitting procedure behaves on them in the same way. This implies that $\nu(I)=\nu(J)$. Then by applying Radon-Nikodym's theorem we have
$$\int_I\frac{d\nu}{d\lambda} d\lambda=\int_J\frac{d\nu}{d\lambda} d\lambda$$
and so  
$$\frac{d\nu}{d\lambda}(I)=\frac{d\nu}{d\lambda}(J)\ \ \text{a.e.}.$$
But $J=f_{I,J}(I)$ for some $c\in]0,1[$ by (\ref{funzIJ}), so the previous relation becomes
$$\frac{d\nu}{d\lambda}(I)=\frac{d\nu}{d\lambda}(f_{I,J}(I))\ \ \text{a.e.}.$$
It follows that $\frac{d\nu}{d\lambda}$ is $\mathcal{F}-$invariant. Hence, since by Lemma~\ref{lem1} the family $\mathcal{F}$ is ergodic, Theorem~\ref{ThmErgod} implies that $\frac{d\nu}{d\lambda}$ is costant a.e.. In particular, we have
$$\frac{d\nu}{d\lambda}=1\ \ \text{a.e.}$$
because
$$1=\nu([0,1])=\int_{[0,1]}\frac{d\nu}{d\lambda}d\lambda=\frac{d\nu}{d\lambda}\cdot\lambda([0,1])
=\frac{d\nu}{d\lambda}.$$
Then $\nu=\lambda$. \\
\endproof

\section{Discrepancy of some generalized Kakutani's sequences}\label{sec:Discr}
A natural problem, posed in \cite{Vol}, is to estimate the behaviour of the discrepancy (\ref{discrPart}) of the sequence of partitions generated by successive $\rho-$refinements as $n$ tends to infinity. In particular, it is interesting to find partitions $\rho$ such that the speed of convergence of the discrepancy to zero is as high as possible. The only known discrepancy bounds for sequences of this kind have been obtained in \cite{Car} by Carbone, who considered the so-called $LS$-sequences that evolve from partitions $\rho$ with $L$ subintervals of $[0,1]$ of length $\alpha$ and $S$ subintervals of length $\alpha^2$ (where $\alpha$ is given by the equation $L\alpha + S\alpha^2 = 1$).

In this section, we analyze this problem with a new approach based on a parsing tree (related to Khodak's coding algorithm \cite{Kho}) which represents the successive $\rho$-refinements.
In particular, we use refinements of the results proved in \cite{Drmota} about Khodak's algorithm to give bounds of the discrepancy for a class of sequences of partitions constructed by successive $\rho-$re\-finements. Finally, we present some examples and applications including $LS$-sequences and u.d.\! sequences of partitions on a class of fractals. These results are collected in \cite{Drm-Inf}.

\subsection{Correspondence between $\rho-$refinements and Khodak's algorithm}
From now on, consider a partition $\rho$ of $[0,1]$ consisting of $m$ intervals of lengths $p_1,\ldots,p_m$ and the sequence of $\rho$-refinements of the trivial partition $\omega=\{[0,1]\}$. Our goal is to construct recursively an $m$-ary tree which represents the process of successive $\rho-$refinements of $\omega$. 
\begin{definition}\ \\
An $m$-ary tree is an ordered rooted tree, where each node has either $m$ ordered successors or it is a leaf with no successors. A node with $m$ successors is called {\it internal node}, while a leaf that has no successors is called also {\it external node}.
\end{definition}
The numbers $p_1,\ldots,p_m$ induce a natural labelling on the nodes. Suppose that the unique path from the root to a node $x$ at level $l$ is encoded by the sequence $(j_1,j_2,\ldots,j_l)$, with $j_i\in\{1,2,\ldots, m\}$, then we set $P(x)=p_{j_1}p_{j_2}\cdots p_{j_l}$. This can be also considered as the probability of reaching the node $x$ with a random walk that starts at the root and moves away from it according to the probabilities $p_1,\ldots,p_m$. For completeness the root $a$ is labelled with $P(a) = 1$. If $T$ is a finite $m$-ary tree then the labels of the external nodes sum up to $1$. Hence, the shape of an $m$-ary tree (together with $p_1,\ldots,p_m$) gives rise of a probability distribution.

The start of our iteration is a tree that only consists of the root which is then an external node with probability $1$. In the first step, the root is replaced by an internal node together with $m$ ordered successing leaves that are given the probability distribution $p_1,\ldots,p_m$. At each further iteration we select all leaves $y$ with the largest label $P(y)$ and grow $m$ children out of each of them. Actually, this construction corresponds to the procedure of successive $\rho$-refinements. The leaves of the tree correspond to the intervals of $\rho^n\omega$ and the labels of the leaves to the lengths of these intervals.

This procedure exactly leads to the same parsing tree of the \emph{Tunstall code} \cite{Drmota}. In fact, the $m$ initial leaves correspond to the symbols of an $m-$ary alphabet $\mathcal{A}$ and so the words $(j_1,j_2,\ldots,j_l)$ that encode the paths from the root to the leaves are the phrases of the dictionary. It is important to note that at each iteration we can have different leaves of the same highest probability, but Tunstall's algorithm selects (randomly) only one of these leaves and grow $m$ children out of it.
 
There is a second way to describe this tree evolution process, namely by \emph{Khodak's algorithm}~\cite{Kho}. Fix a real number $r\in]0,p_{min}[$, where $p_{min}=\min\{p_1,\ldots,p_m\}$, and consider all nodes $x$ among in an infinte $m$-ary tree with $P(x)\geq r$. Let us denote these nodes by $\mathcal{I}(r)$. Of course, if $P(x)\geq r$ then all nodes $x'$ on the path from the root to $x$ satisfy $P(x')\geq r$, too. Hence, these nodes of $\mathcal{I}(r)$ constitute a finite subtree. These nodes will be the $\emph{internal nodes}$ of Khodak's construction. Finally, we append to these internal nodes all successor nodes $y$. By construction all these nodes satisfy $p_{min}r\leq P(y) < r$ and we denote them by $\mathcal{E}(r)$. These nodes are the $\emph{external nodes}$ of Khodak's construction. 
We denote by $M_r = |\mathcal{E}(r)|$ the number of external nodes.
Obviously, we have got a finite $m$-ary tree
$\mathcal{T}(r)=\mathcal{I}(r)\cup\mathcal{E}(r)$
and it is clear that these trees grow when $r$ decreases. For certain values $r$, precisely the external nodes $y$ of largest value $P(y)=r$ turn into internal nodes and all their successors become new external nodes. Actually, the tree $\mathcal{T}(r)$ grows in correspondence to a decreasing sequence of values $(r_j)$. Indeed, when $r\in]r_j,r_{j-1}]$ the tree remains the same, i.e $\mathcal{T}(r_{j-1})=\mathcal{T}(r)$. 

The parsing tree resulting from Khodak's algorithm is exactly the same as the tree constructed by Tunstall's algorithm. However, we have to observe that in Khodak's construction all leaves with the same highest probability are selected to generate the children at once, while in Tunstall'algortithm they are selected one by one in an arbitrary order. Now, in the procedure of successive $\rho-$refinements, at each step we select the intervals having maximal lenght at once and we split them at the same time. So, Khodak's algorithm and $\rho-$refinements procedure not only are exactly represented by the same tree but they also have a common structure which allows to create a useful correspondence between them. 

In fact, if we fix a step $j$ in $\rho-$refinements procedure, then the tree associated to the partition $\rho^j\omega$ is exactly $\mathcal{T}(r_j)$. Therefore, we will only consider the values of the sequence $(r_j)$ for which the tree constructed by Khodak's algorithm actually grow. Note that we have to start with the value $r_1=1$, because we intend to consider the whole procedure since the first step which corresponds to $\rho\omega$. Hence, in our correspondence between Khodak's algorithm and the procedure of successive $\rho-$refinements the value $r_j\in]0,1]$ corresponds to the partition $\rho^{j}\omega$. Consequently, the number of external nodes in $\mathcal{E}(r_j)$ equals the number of points defining the partition $\rho^{j}\omega$, i.e. $M_{r_j}=k(j)$. Moreover, if $r\in]r_j,r_{j-1}]$ then $M_{r}=M_{r_{j-1}}=k(j-1)$.

From here on we denote by $\mathscr{E}_{r_j}$ the family of all intervals of the partition $\rho^{j}\omega$ corresponding to the leaves belonging to $\mathcal{E}(r_j)$ and the order of the intervals in $\mathscr{E}_{r_j}$ corresponds to the left-to-right order of the external nodes in $\mathcal{E}(r_j)$. We will call \emph{elementary intervals} all the intervals belonging to each $\mathscr{E}_r$ for $r\in]0,1]$.

In the following we denote by $H$ the entropy of the
probability distribution $p_1,\ldots,p_m$ , which is defined as
$$H=p_1\log\left(\frac{1}{p_1}\right)+\dots+p_m\log\left(\frac{1}{p_m}\right).$$

\begin{figure}
\centering
\subfloat[I step]{\label{t1}\includegraphics[width=13cm]{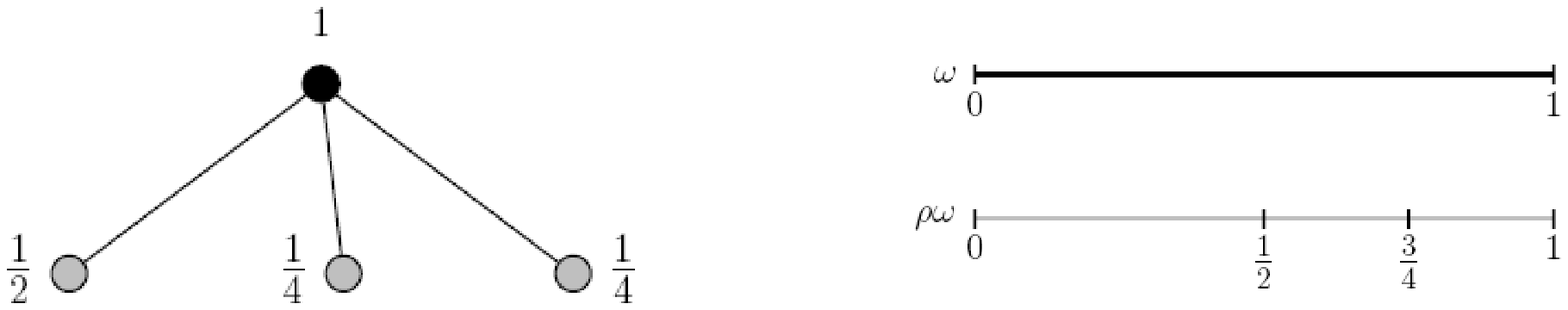}} 
\vspace{0.2cm}  
            
\subfloat[II step]{\label{t2}\includegraphics[width=13cm]{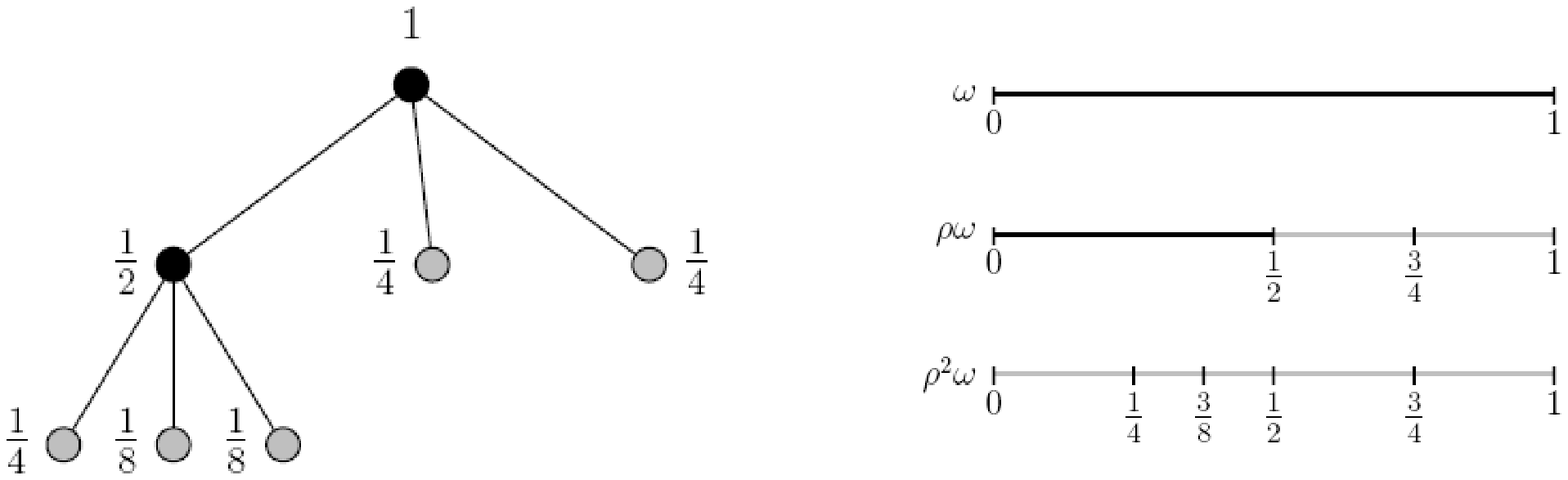}}
\vspace{0.2cm} 

\subfloat[III step]{\label{t3}
$
\begin{array}{ll}
\includegraphics[width=7cm]{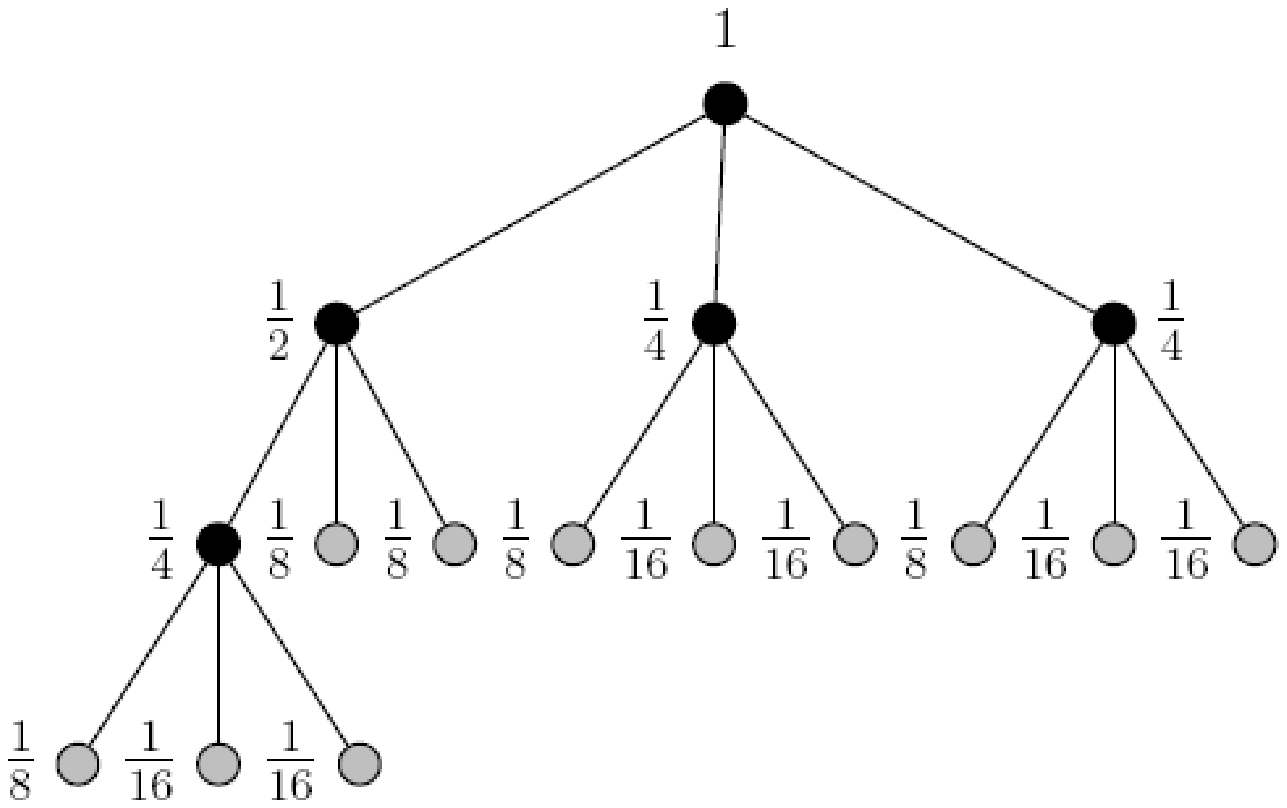}
&
\includegraphics[width=6cm]{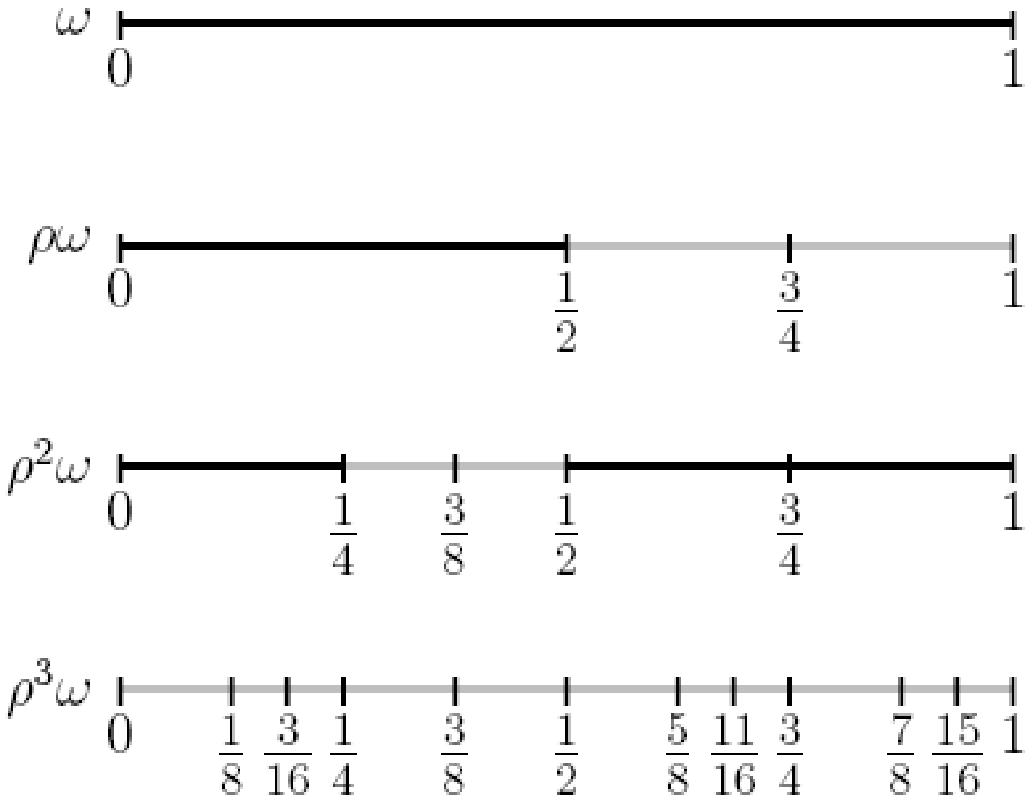}
\end{array}
$
}
\vspace{0.2cm}
 
\caption{Correspondence between $\rho-$refinements and Khodak's tree.}
\label{FigCorr}
\end{figure}

\paragraph{Example}\ \\
Let $\rho=\left\{\left[0,\frac 14\right], \left[\frac 14, \frac 12\right], \left[\frac 12, 1\right]\right\}$ and $\omega=\{[0,1]\}$. So in this case we have that $p_1=p_2=\frac 14$ and $p_3=\frac 12$. In Figure \ref{FigCorr} the correspondence between the tree constructed by Khodak's algorithm and the $\rho-$refinements of $\omega$ is illustrated in the first three steps of these procedures. In particular, the internal nodes are coloured in black and the external ones in grey. Moreover, it is easy to note that the label of each node is exactly the lenght of the corresponding interval.

The start of the procedure of $\rho-$refinements is the trivial partition $\omega$ which corresponds to the root node of probability $1$. At the first step we take $r_1=1$, so the root becomes an internal node and we grow three leaves out of it. The three children nodes have probability $p_1,p_2$ and $p_3$ respectively and each of them corresponds to an interval of $\rho\omega$ (see Figure \ref{t1}). The next value of $r$ for which we have an actual growth of the tree is $r_2=\frac 12$. Consequently, at the second step, we select all the nodes $x$ having $P(x)\geq r_2$ and grow three children out of each of them. The external nodes generated at the end of this step correspond exactly to the intervals of $\rho^2\omega$ (see Figure \ref{t2}). At the third step we choose $r_3=\frac 14$ and we proceed at the same way of the previous steps and so we get $11$ leaves which are associated to the intervals of $\rho^3\omega$ (see Figure \ref{t3}). By iterating this procedure for all the values of $(r_j)$, we will get the whole infinity tree corresponding to the sequence of partitions $(\rho^j\omega)$. \\

\subsection{Asymptotic results on Khodak's algorithm}
One of the main results from \cite{Drmota} provides asymptotic information on the number $M_r$ of external nodes in Khodak's construction. Actually, these relations can be used to prove Theorem~\ref{VolThm}. However, in order to derive bounds for the discrepancy of the sequence $(\rho^n\omega)$ we need more precise information on the error terms in the asymptotic results given in \cite{Drmota}. Therefore, we will extend the analysis of this paper in Theorem \ref{Thm1}. Before introducing it let us give a fundamental definition, which allows us to distinguish two main cases in our discussion.

\begin{definition}\label{defRR}\ \\
We say that $\log\left(\frac{1}{p_1}\right),\ldots,\log\left(\frac{1}{p_m}\right)$ are \emph{rationally related} if there exists a positive real number $\Lambda $ such that $\log\left(\frac{1}{p_1}\right),\ldots,\log\left(\frac{1}{p_m}\right)$ are integer multiples of $\Lambda $, that is
$$\log\left(\frac{1}{p_j}\right)=n_j\Lambda ,\quad \textrm{with}\ n_j\in\mathbb{Z}\ \textrm{for}\ j=1,\ldots,m.$$
Equivalently, all fractions $(\log p_i)/(\log p_j)$ are rational. Without loss of generality we can assume that $\Lambda $ is as large as possible which is equivalent to assume that $\gcd(n_1,\ldots,n_m)=1$. 

Similarly we say that $\log\left(\frac{1}{p_1}\right),\ldots,\log\left(\frac{1}{p_m}\right)$ are \emph{irrationally related} if they are not rationally related.
\end{definition}

\begin{theorem}\label{Thm1}\ \\
Let $M_r$ be the number of the external nodes generated at the step corresponding to the parameter $r$ in Khodak's construction, that is, the number of nodes in $\mathcal{E}(r)$. 
\begin{enumerate}
\item If $\log\left(\frac{1}{p_1}\right),\ldots,\log\left(\frac{1}{p_m}\right)$ are rationally related, let $\Lambda >0$ be the largest real number for
which $\log\left(\frac{1}{p_j}\right)$ is an integer multiple of $\Lambda $ (for $j=1,\ldots,m$). Then there exists a real number $\eta> 0$ and an integer
$d\ge 0$ such that
\begin{equation}\label{ratCase}
M_r=\frac{(m-1)}{rH}Q_1\left(\log\left(\frac{1}{r}\right)\right)+
\mathcal{O}\left((\log r)^d r^{-(1-\eta)}\right),
\end{equation}
where
$$Q_1(x)=\frac{\Lambda }{1-e^{-\Lambda }}e^{-\Lambda \left\{ \frac{x}{\Lambda }\right\} }$$
and $\{ y\}$ is the fractional part of the real number $y$.
Furthermore, the error term is optimal.
\item If $\log\left(\frac{1}{p_1}\right),\ldots,\log\left(\frac{1}{p_m}\right)$ are irrationally related, then 
\begin{equation}\label{irratCase}
M_r=\frac{(m-1)}{rH}+{o}\left(\frac{1}{r}\right).
\end{equation}
In particular, if $m=2$ and $\gamma = (\log p_1)/(\log p_2)$ is
badly approximable then
\begin{equation*}
M_r=\frac{(m-1)}{rH}\left( 1 + \mathcal{O}\left( \frac{ (\log\log (1/r))^{1/4} }
{ (\log (1/r))^{1/4} }\right) \right).
\end{equation*}
Moreover, if $p_1$ and $p_2$ are algebraic then there exists an
effectively computable constant $\kappa> 0$
with
\begin{equation*}
M_r=\frac{(m-1)}{rH}\left( 1 + \mathcal{O}\left( \frac{ (\log\log(1/r))^{\kappa} }
{ (\log (1/r))^{\kappa} } \right)\right).
\end{equation*}
\end{enumerate}
\end{theorem}

The proof of this theorem requires some auxiliary results, so let us introduce them preliminarly.

\subsubsection{Trigonometric Sums}

\begin{lemma}\label{LeApp1}\ \\
Let $f(n)=\sum\limits_{i=1}^kc_i\cos(2\pi\theta_in+\alpha_i)$, $c_i,\alpha_i,\theta_i\in\mathbb{R}$
 be defined for non-negative integers
$n$ and suppose that $f$ is not identically zero. Then
there exists $\delta>0$ such that $|f(n)|\geq\delta$ for infinitely many 
non-negative integers $n$.
\end{lemma}
\proof\ \\
We have to distinguish two cases:
\begin{description}
\item[Case 1] $\theta_1,\ldots,\theta_k$ are rationally related. \\
There exist $\Lambda \in\mathbb{R}\setminus\{0\}$ and $k_i\in\mathbb{Z}$ such that $\theta_i=\Lambda  k_i$.
In this case, we can rewrite the function $f$ as follows
$$f(n)=\sum_{i=1}^kc_i\cos(2\pi\Lambda  n k_i+\alpha_i)=\sum_{i=1}^kc_i\cos(2\pi\{\Lambda  n\} k_i+\alpha_i).$$

Hence, $f(n)=g(\{\Lambda  n\})$ where $g(x)=\sum\limits_{i=1}^kc_i\cos(2\pi k_ix+\alpha_i)$
is a periodic non-zero function of period 1.
\begin{description}
\item[Case 1.1] If $\Lambda \in\mathbb{Q}$, then $\Lambda =\frac{p}{q}$ for some 
coprime integers $p,q\in\mathbb{Z}$ and the sequence $(f(n))$ attains
periodically the set of values
\[
g\left(\left\{\frac{pn}{q}\right\}\right), \quad n=0,\ldots, q-1.
\]
Since they are not all equal to zero there exists $\delta>0$ such that 
$$|f(n)|=|g(\{\Lambda  n\})|\geq\delta$$ for infinitely many $n$.
In particular, we can use a linear subsequence
$qn+r$ for which $|f(qn+r)|\ge \delta$.
\item[Case 1.2] If $\Lambda \notin\mathbb{Q}$, then the sequence $(\{\Lambda  n\})$ is u.d.\! on $[0,1]$ and consequently dense in $[0,1]$. Hence, there again exists $\delta>0$ such that $$|f(n)|=|g(\{\Lambda  n\})|\geq\delta$$ for infinitely many $n$. 
\end{description}
\item[Case 2] $\theta_1,\ldots,\theta_k$ are irrationally related. \\
Here we divide the $\theta_i$ in groups which are rationally related. Assume that we have $s$ groups $\{ \theta_i : i\in I_j\}$ for $j=1,\ldots,s$, 
and in each group we write 
$$\theta_{i}=\Lambda _jk_i, \quad i\in I_j$$
with $k_i\in\mathbb{Z}$ and some $\Lambda_j\in \mathbb{R}\setminus\{0\}$.

In this case, we distinguish between three different subcases:
\begin{description}
\item[Case 2.1] $1,\Lambda _1,\ldots,\Lambda _s$ are linearly independent over $\mathbb{Q}$ (so~$\Lambda _1,\ldots,\Lambda _s\notin\mathbb{Q}$\ ).\\ 
We set $f_j(x)=\sum\limits_{i\in I_j}c_i\cos(2\pi xk_i+\alpha_i)$ (where we assume
w.l.o.g.\ that $f_j$ is non-zero) and $g(x_1,\ldots, x_s)=\sum\limits_{j=1}^sf_j(x_j)$.
Then
$$
f(n)=\sum_{j=1}^sf_j(\{n\Lambda _j\}) = g\left(\{n\Lambda _1\},\ldots,\{n\Lambda _s\}\right).
$$
By Kronecker's Theorem (Theorem \ref{KronMultidim}), the sequence $\left(\{n\Lambda _1\},\ldots,\{n\Lambda _s\}\right)$ is dense in the cube $[0,1]^s$. 
Thus, it follows (as above) that there exists $\delta>0$ such that $|f(n)|\geq\delta$
for infinitely many $n$.

Note that by same reasoning it follows that for every $\varepsilon> 0$ we
have $|f(n)| \le \varepsilon$ for infinitely many $n$. 
(Here we also use the fact that $f$ has zero mean.)
This observation will be used in Case 2.3.

\item[Case 2.2] $1,\Lambda _1,\ldots,\Lambda _s$ are linearly dependent over $\mathbb{Q}$ and $\Lambda _1,\ldots,\Lambda_s\notin\mathbb{Q}$\ .

In this case there exist $q,p_1,\ldots,p_s\in\mathbb{Z}$ such that $q=p_1\Lambda _1+\ldots+p_s\Lambda _s$.\\
Suppose (w.l.o.g.)\! that $p_1 > 0$ and consider the subsequence of integers~$(p_1n)$, then
\begin{eqnarray*}
f(p_1n)&=&\sum_{j=1}^s f_j(n \Lambda_j p_1) \\
		&=& f_1(n(q-\Lambda_2 p_2 - \cdots - \Lambda_s p_s)) 
		+ \sum_{j=2}^s f_j(n \Lambda_j p_1).
\end{eqnarray*}
By using the addition theorem for cosine and rewriting the sum accordingly, we
obtain a representation of the form 
\[
f(p_1n) = \sum_{j=2}^s \tilde f_j (n \Lambda_j p_j),
\]
where $\tilde f_j$ are certain trigonometric polynomials.
This means that we have eliminated~$\Lambda_1$. 

In this way we can proceed further. If 
$1,p_2\Lambda_2,\ldots,p_s\Lambda_s$ are linearly independent over~$\mathbb{Q}$ 
then we argue as in Case 2.1. However, if 
$1,p_2\Lambda_2,\ldots,p_s\Lambda_s$ are linearly dependent over~$\mathbb{Q}$ 
then we repeat the elimination procedure etc. Note that this elimination
procedure terminates, since we assume that 
$\Lambda _1,\ldots,\Lambda_s\notin\mathbb{Q}$. 
Hence, we always end up in Case~2.1.

\item[Case 2.3] $\Lambda _1,\ldots,\Lambda _s$ are not all irrationals. \\
Here we represent $f(n) = h_1(n)+h_2(n)$, where
\[
h_1(n) = \sum_{j\in\{j:\Lambda _j\in\mathbb{Q}\}} f_j(n)
\quad \mbox{and}\quad
h_2(n) = \sum_{j\in\{j:\Lambda _j\not\in\mathbb{Q}\}} f_j(n).
\]
If $h_1$ is non-zero then we can argue as in Case~1.1. All appearing
$\theta_i$ are rational and consequently there exists a linear subsequence
$qn+r$ such that $|h_1(qn+r)|\ge\frac {3\delta} 2$ for some $\delta > 0$.
Next we reduce the sum $h_2(qn+r)$ to a sum of the form that is discussed
in Case~2.1 (possibly we have to eliminate several terms as discussed
in Case~2.2). Consequently, it follows that there exist infinitely many
$n$ such that $|h_2(qn+r)| \le \delta/2$. Hence we have
$|f(n)| \ge \delta$ for infinitely many $n$.

If $h_1$ is zero, for all non-negative integers we just have to 
consider $h_2$. But this case is precisely that of Case 2.2.
\end{description}
\end{description}
\endproof

\subsubsection{Zerofree Regions}
The purpose of this paragraph is to discuss zerofree regions of the equation 
\begin{equation}\label{eq}
1-p^{-s}-q^{-s}=0
\end{equation}
where $p,q$ are positive numbers with $p+q = 1$. 

It is clear that $s=-1$ is a solution of (\ref{eq}) and that all its solutions have to satisfy $\Re(s) \ge -1$. 
Otherwise, we would have 
$$|p^{-s}+q^{-s}|\leq|p^{-s}| + |q^{-s}|=e^{-\Re(s)\log(p)}+e^{-\Re(s)\log(q)}<
e^{\log(p)}+e^{\log(q)}=p+q= 1.$$

Furthermore, it is easy to verify the following property. 
\begin{prop}\ \\
There are no solutions of (\ref{eq}), other than $s= -1$, on the line $\Re(s) = -1$ if and only if the ratio $\gamma = (\log p)/(\log q)$ is irrational.
\end{prop}
\proof\ \\
Let $s$ be of the form $s=-1+it$ with $t\neq 0$. Then we have
$$p^{-s}+q^{-s}=p^{1-it}+q^{1-it}=pp^{-it}+qq^{-it}
			=pe^{it\log\left(\frac{1}{p}\right)}+qe^{it\log\left(\frac{1}{q}\right)}.$$
Therefore, $s$ is a solution of (\ref{eq}) if and only if
$$pe^{it\log\left(\frac{1}{p}\right)}+qe^{it\log\left(\frac{1}{q}\right)}=1.$$
Since $p+q=1$, then
$$e^{it\log\left(\frac{1}{p}\right)}=1=e^{it\log\left(\frac{1}{q}\right)}$$
necessarily holds. But the last relations imply that there exist $k,l\in\mathbb{Z}$ such that
$$
t\log\left(\frac{1}{p}\right)=2\pi k\quad\textrm{and}\quad t\log\left(\frac{1}{q}\right)=2\pi l.
$$
It follows that $\gamma=\frac{\log\left(\frac{1}{p}\right)}{\log\left(\frac{1}{q}\right)}=\frac{k}{l}$, which is a rational number.\\
\endproof

Another important property about the structure of the set of the solutions of~(\ref{eq}) is the following (see \cite{Schachinger}).
\begin{prop}\label{SchachProp}\ \\
There exist two real numbers $\sigma_0>-1$ and $\tau> 0$
such that in each box of the form
\[
B_k = \{s\in \mathbb{C} : -1 \le \Re(s) \le \sigma_0,\, 
(2k-1)\tau \le \Im(s) < (2k+1)\tau \}, \quad k \in 
\mathbb{Z}\setminus \{ 0 \},
\]
there is precisely one zero of (\ref{eq}) and there are no other zeros.
\end{prop} 
However, the position of the zeros in the $B_k$'s is by no means clear. Nevertheless, with the help of the continued fractional expansion of $\gamma$ it is possible to construct infinitely many zeros $s$ of (\ref{eq}) with $\Re(s) < -1 + \varepsilon$ (for every
$\varepsilon> 0$). Therefore, it is natural to ask for zerofree regions
of this equation. Actually, one has to assume some Diophantine conditions on $\gamma$ to get precise information.

\begin{lemma}\label{freeZeroReg}\ \\
If $\gamma$ is badly approximable then for every solution $s\neq-1$ of the equation $$1-p^{-s}-q^{-s}=0$$ we have that
$$
\Re(s)>\frac{c}{(\Im(s))^2}-1
$$
for some positive constant $c$.
\end{lemma}

Before proving the lemma, let us recall some basic notions of the theory of continued fractions \cite{Khi}.  
\begin{definition}\ \\
The continued fractional expansion of a real number $x$ is given by
$$x=a_0+\frac{1}{a_1+\frac{1}{a_2+\frac{1}{a_3+\ldots}}}$$
where $a_0,a_1,\ldots$ are positive integers. In a compact notation we can write
$$x=[a_0; a_1; a_2;a_3 \ldots].$$
\end{definition}
\begin{definition}\ \\
An irrational number $\gamma$ is said to be badly approximable if its continued fractional expansion $\gamma =[a_0; a_1; \ldots]$ is bounded, i.e. there exists a positive constant $D$ such that $\max\limits_{j\geq 1}(a_j)\leq D$. 
\end{definition}
Equivalently, we have the following property.
\begin{prop}\ \\
An irrational number $\gamma$ is badly approximable if there exists a constant $d>0$ such that 
\begin{equation}\label{badAppr}
\left|\gamma-\frac{k}{l}\right|\geq\frac{d}{l^2}
\end{equation}
for all non-zero integers $k,l$.
\end{prop}

\proof[Proof of Lemma \ref{freeZeroReg}]\ \\
In order to make the presentation of the proof more transparent
we make a shift by~$1$ and consider the equation
\begin{equation}\label{eqshift}
p^{1-s} + q^{1-s} = 1
\end{equation}
and show that all its non-zero solutions satisfy 
$\Re(s) > c/\Im(s)^2$  
for some positive constant~$c$ that depends on~$\gamma$. 

Suppose that $s = \sigma + i\tau$ is a zero of (\ref{eqshift})
with $\sigma > 0$. Furthermore, we assume that $\sigma \le \varepsilon$,
where $\varepsilon$ is a sufficiently small constant. Since $p+q = 1$ and $|p^{1-s}| = p^{1-\sigma} = p(1+\mathcal{O}(\varepsilon))> p$
and $|q^{1-s}| = q^{1-\sigma} = q(1+\mathcal{O}(\varepsilon))> q$ we can only have
a solution if the arguments of $p^{1-s}$ and $q^{1-s}$ are small. (Actually they have to be of order $\mathcal{O}(\sqrt{\varepsilon})$ if
$\varepsilon$ is chosen sufficiently small). W.l.o.g.\ we write
\[
\arg(p^{1-s}) = \tau \log (1/p) = 2\pi k + \eta_1 \quad\mbox{and}
\quad \arg(q^{1-s}) = \tau \log (1/q) = 2\pi l - \eta_2
\]
for some integers $k,l$ and certain positive numbers $\eta_1$, $\eta_2$ 
(which are of order $\mathcal{O}(\sqrt{\varepsilon})$).
More precisely, by doing a local expansion in (\ref{eqshift})
we obtain 
\[
\eta_2 = \frac pq \eta_1 + \mathcal{O}(\eta_1^2)\quad\mbox{and}
\quad \sigma = \frac{p}{2qH} \eta_1^2 + \mathcal{O}(\eta_1^4).
\]

In fact (\ref{eqshift}) is equivalent to
$$
|p^{1-s}|(\cos(\arg(p^{1-s}))+i\sin(\arg(p^{1-s})))
+|q^{1-s}|(\cos(\arg(q^{1-s}))+i\sin(\arg(q^{1-s})))=1
$$
and so
\begin{equation}\label{trigFrom}
p^{1-\sigma}(\cos(\eta_1)+i\sin(\eta_1))+q^{1-\sigma}(\cos(\eta_2)-i\sin(\eta_2))=1.
\end{equation}
Therefore we have
$$p^{1-\sigma}\sin(\eta_1)-q^{1-\sigma}\sin(\eta_2)=0$$
and by doing a local expansion it follows that
\begin{eqnarray*}
&&p\left(1+\log\left(\frac 1p\right)\sigma+\mathcal{O}(\sigma^2)\right)\cdot(\eta_1+\mathcal{O}(\eta_1^2))\\
&+&q\left(1+\log\left(\frac 1q\right)\sigma+\mathcal{O}(\sigma^2)\right)\cdot(-\eta_2+\mathcal{O}(\eta_2^2))=0.
\end{eqnarray*}
Now, by taking into account that $\sigma\leq\epsilon$ and $\eta_1, \eta_2$ are both of the order $\mathcal{O}(\sqrt{\varepsilon})$, we have
$$p\eta_1-q\eta_2+\mathcal{O}(\eta_1^2)=0$$
and so
\begin{equation}\label{EtaDue}
\eta_2 = \frac pq \eta_1 + \mathcal{O}(\eta_1^2).
\end{equation}
On the other hand, from (\ref{trigFrom}) it follows also that 
$$p^{1-\sigma}\cos(\eta_1)+q^{1-\sigma}\cos(\eta_2)=0$$
and by doing a local expansion we have
\begin{eqnarray*}
&&p\left(1+\log\left(\frac 1p\right)\sigma+\mathcal{O}(\sigma^2)\right)\cdot\left(1-\frac{\eta_1^2}{2}+\mathcal{O}(\eta_1^4)\right)\\
&+&q\left(1+\log\left(\frac 1q\right)\sigma+\mathcal{O}(\sigma^2)\right)\cdot\left(1-\frac{\eta_2^2}{2}+\mathcal{O}(\eta_2^4)\right)=1.
\end{eqnarray*}
Now, by using the same argumentations of above, we have
$$\sigma\left(p\log\left(\frac 1p\right)+q\log\left(\frac 1q\right)\right)= \frac{p\eta_1^2}{2}+\frac{q\eta_2^2}{2}+\mathcal{O}(\eta_1^4)$$
and so by (\ref{EtaDue}) we get
$$\sigma=\frac{p}{2H} \eta_1^2+ \frac{p^2}{2Hq} \eta_1^2+ \mathcal{O}(\eta_1^4)=\frac{p}{2qH} \eta_1^2 + \mathcal{O}(\eta_1^4).$$

Furthermore, we have
\begin{align*}
\gamma &= \frac {\tau \log \frac 1p} {\tau \log \frac 1q} \\
& = \frac{2\pi k + \eta_1}{2\pi l - \eta_2} \\
&= \frac kl\left(1+\frac{\eta_1}{2\pi k}\right)\left(1+\frac{p\eta_1}{q 2\pi l}+\mathcal{O}(\eta_1/l^2)\right)\\
&= \frac kl + \frac 1{2\pi}\left( \frac 1l + \frac {kp}{l^2 q} \right) \eta_1
(1+ \mathcal{O}(\eta_1/l)).
\end{align*}
This means that $k/l$ is close to $\gamma$ and by applying (\ref{badAppr})
it follows that 
\[
\eta_1 \ge \frac {d'}{|l|}
\]
for some constant $d'> 0$. 
Consequently, we obtain $\sigma \ge d''/l^2$ (for some constant $d''>0$)
which translates directly to 
$\sigma > c/\tau^2$ for some positive constant $c$.\\
\endproof

Next we consider the case of algebraic number $p$ and $q$
such that $\log(p)/\log(q)$ is irrational.

\begin{lemma}\label{freeZeroRegBak}\ \\
If $p,q\in]0,1[$ are positive algebraic numbers with $p+q =1$
and with the property that $\log(p)/\log(q)$ is irrational.
Then for every solution $s\neq-1$ of the equation 
$$1-p^{-s}-q^{-s}=0$$ we have
$$
\Re(s)>\frac{D}{(\Im(s))^{2C}}-1
$$
with effectively computable positive constants $C,D$.
\end{lemma}

The classical theorem of Gelfond-Schneider says that if
$\gamma = \log(p)/\log(q)$ is irrational for algebraic numbers $p$ and $q$ 
then $\gamma$ is transcendental. Baker's theorem (see~\cite{Bak}) gives also effective bounds for Diophantine approximation of $\gamma$ that will be used 
in the subsequent proof of Lemma~\ref{freeZeroRegBak}.
Before introducing Baker's theorem, let us recall that the \emph{height} of an algebraic number is the maximum of the absolute values of the relatively prime integer coefficients in its minimal defining polynomial, while its \emph{degree} is the degree of this polynomial.
\begin{theorem}[Baker's Theorem]\label{ThmBak}\ \\ 
Let $\gamma_1,\ldots,\gamma_n$ be non-zero algebraic numbers with degrees at most $d$ and heights at most $A$. Further, $\beta_0,\beta_1,\ldots,\beta_n$ are algebraic numbers with degree at most $d$ and heights at most $B$ ($\geq 2)$. Then for
$$\Lambda =\beta_0+\beta_1\log\gamma_1+\ldots+\beta_n\log\gamma_n$$
we have either $\Lambda =0$ or $|\Lambda |\geq B^{-C}$, where $C$ is an effectively computable number depending only on $n,d$, and $A$.
\end{theorem}

\proof[Proof of Lemma~\ref{freeZeroRegBak}]\ \\
We apply Theorem~\ref{ThmBak} to the algebraic numbers $\gamma_1=p$ and 
$\gamma_2 = q$ and to the integers $\beta_0 =0$, $\beta_1 = l$, and
$\beta_2 = -k$. Then $B = \max\{|k|,|l|\}$.  W.l.o.g. we may assume that
$p> q$ which assures that we only have to consider cases with $|k|\le |l|$.
Thus 
$$\left|l\log p-k\log q\right|>B^{-C}$$
and consequently
\begin{equation}\label{caseBak}
\left|\frac{\log p}{\log q}-\frac{k}{l}\right|>\left(\frac{1}{\log q}\right)\frac{B^{-C}}{l}>
\left(\frac{1}{\log q}\right)\frac{1}{l^{1+C}},
\end{equation}
where $C$ is effectively computable.
 
By using (\ref{caseBak}) instead of (\ref{badAppr}) in the proof of Lemma~\ref{freeZeroReg} we easily complete the proof of Lemma~\ref{freeZeroRegBak}.\\
\endproof

\subsubsection{Differentiating Asymptotic Expansions}
For our analysis, we need a Tauberian theorem for the Mellin transform. A classical result in this direction is the following (see \cite{Korevaar}, \cite{Hen}).
\begin{theorem}\ \\
Suppose that $f(v)$ is a monotone function for $v\geq 0$ such that
$$F(v)=\int_0^v f(w)dw$$
is asymptotically given by
$$F(v)\sim\frac{v^{\lambda+1}}{(\lambda+1)}\quad \textrm{as} \quad v\to\infty,$$
for some $\lambda>-1$. Then
$$f(v)\sim v^\lambda\quad \textrm{as} \quad v\to\infty. $$
\end{theorem}
We make this result more precise in the next lemma.
\begin{lemma}\label{Taub}\ \\
Suppose that $f(v)$ is a non-negative increasing function for $v\geq 0$. Assume that  
$$F(v)=\int_0^v f(w)dw$$
has the asymptotic expansion
$$F(v)=\frac{v^{\lambda+1}}{\lambda+1}
\left(1+\mathcal{O}\left(g(v)\right)\right)\quad \textrm{as} \quad v\to\infty,$$
where $\lambda>-1$ and $g(v)$ is a decreasing function that tends to 
zero as $v\to\infty$. Then
$$f(v)= v^\lambda\left(1+\mathcal{O}\left(g(v)^{\frac{1}{2}}\right)\right)\quad \textrm{as} \quad v\to\infty. $$
\end{lemma}
\proof\ \\
By the assumption we have that there exist $v_0,c>0$ such that for all $v\geq v_0$ we have
$$\left|F(v)-\frac{v^{\lambda+1}}{(\lambda+1)}\right|\leq 
c|g(v)|\frac{v^{\lambda+1}}{(\lambda+1)}.$$
Now, set $h=|g(v)|^{\frac{1}{2}}v$. By monotonicity, for $v\geq v_0$ we get
$$
\frac{F(v+h)-F(v)}{h}=\frac{1}{h}\int_v^{v+h}f(w)dw\geq
\frac{1}{h}\int_v^{v+h}f(v)dw=f(v).
$$
Consequently we have
\begin{eqnarray*}
f(v)&\leq&\frac{F(v+h)-F(v)}{h}\\
&\leq&\frac{1}{h}\left(\frac{(v+h)^{\lambda+1}}{\lambda+1}-\frac{v^{\lambda+1}}{\lambda+1}\right)
+\frac{1}{h}\left(c|g(v+h)|\frac{(v+h)^{\lambda+1}}{(\lambda+1)}+c|g(v)|\frac{v^{\lambda+1}}{(\lambda+1)}\right)\\
&\leq&\frac{1}{h(\lambda+1)}\left(v^{\lambda+1}+(\lambda+1)v^\lambda h
+\mathcal{O}(v^{\lambda-1}h^2)-v^{\lambda+1}\right)+\mathcal{O}\left( |g(v)|\frac{v^{\lambda+1}}{h}\right)\\
&=&v^{\lambda}+\mathcal{O}\left({v^{\lambda-1}h}\right)+\mathcal{O}\left( |g(v)|\frac{v^{\lambda+1}}{h}\right)\\
&=&v^{\lambda}+\mathcal{O}\left({v^{\lambda-1}|g(v)|^{\frac{1}{2}}v}\right)+\mathcal{O}\left( |g(v)|\frac{v^{\lambda+1}}{|g(v)|^{\frac{1}{2}}v}\right)\\
&=&v^{\lambda}+\mathcal{O}\left(v^{\lambda}|g(v)|^{\frac{1}{2}}\right).
\end{eqnarray*}
\endproof

\subsubsection{Proof of Theorem \ref{Thm1}}
Set $v=\frac{1}{r}$ and denote by $A(v)$ the number of internal nodes (root node included) in Khodak's construction with parameter $r = 1/v$, i.e.
$$A(v)=\sum_{x: P(x)\geq\frac{1}{v}}1.$$
Hence, the number of external nodes generated at the step corresponding to the parameter $r$ is
\begin{equation}\label{eqMrrel}
M_r=(m-1)A(v)+1.
\end{equation}
The key relation is that that $A(v)$ satisfies the following recurrence (see \cite[Lemma~2]{Drmota})
\begin{equation} \label{eqrecA}
A(v)=\left\{
\begin{array}{ll}
0&\ v<1\\
1+\sum\limits_{j=1}^m A(p_jv)&\ v\geq 1
\end{array}\right..
\end{equation}
For the asymptotic analysis of $A(v)$ and consequently that of $M_r$ we distinguish between the rational and the irrational case.
\paragraph{\emph{Rational case}}\ \\
If the $\log(1/p_j)$'s are rationally related then $A(v)$ is constant for $v\in [e^{\Lambda n},e^{\Lambda (n+1)}[$ (for every integer $n$). Hence, it suffices to study the behaviour of the sequence $G(n) = A(e^{\Lambda n})$ which verifies the recurrence
$$
G(n)=\left\{
\begin{array}{ll}
0&\ n<0\\
1 + \sum\limits_{j=1}^m G(n-n_j)&\ n\geq 0
\end{array}\right.
$$
where $n_j=\frac{\log\left(\frac{1}{p_j}\right)}{\Lambda}$ for $j=1,\ldots,m$. The generating function
$g(z) = \sum\limits_{n\ge 0} G(n) z^n$ is then given by
\begin{eqnarray*}
g(z)&=&
\sum_{n\ge 0} G(n)z^n=\sum_{n\ge 0}\left(1+\sum_{j=1}^mG(n-n_j)\right)z^n\\
&=&\sum_{n\ge 0} z^n+\sum_{n\ge 0}\sum_{j=1}^mG(n-n_j)z^n\\
&=&\frac{1}{1-z}+\sum_{j=1}^m\sum_{n=-n_j}^\infty G(n)z^{n+n_j}\\
&=&\frac{1}{1-z}+\sum_{j=1}^m z^{n_j}\left(\sum_{n\ge 0} G(n)z^n\right)\\
&=&\frac{1}{1-z}+g(z)\sum_{j=1}^m z^{n_j}
\end{eqnarray*}
and so
\[
g(z) = \frac 1{(1-z)f(z)},
\]
where $f(z)=1-z^{n_1}+\cdots-z^{n_m}$. 
By Definition \ref{defRR}, it follows that $e^{-\Lambda }$ is a positive real root of $f$. Moreover, it is proved in \cite{CG} that if we denote by $\omega_1,\ldots,\omega_h$ all the other different roots with multiplicities $\mu_i$ of $f$, then $|\omega_i|>e^{-\Lambda }$ for $i=1,\ldots,h$. (Here we use the assumption that $n_1,\ldots,n_m$ are coprime.) Therefore, we have by partial fraction decomposition
\begin{eqnarray*}
g(z)&=&\frac{a}{1-z}+\frac{b}{z-e^{-\Lambda}}+
\frac{c^{(1)}_1}{z-\omega_1}+\cdots+\frac{c^{(1)}_{\mu_1}}{(z-\omega_1)^{\mu_1}}
+\cdots\\
&+&\frac{c^{(h)}_1}{z-\omega_h}+\cdots+\frac{c^{(h)}_{\mu_h}}{(z-\omega_h)^{\mu_h}}\\
&=&\frac{a}{1-z}+\frac{-be^\Lambda}{1-ze^{\Lambda}}+
\frac{-c^{(1)}_1\omega_1^{-1}}{1-z\omega_1^{-1}}+\cdots+\frac{-c^{(1)}_{\mu_1}\omega_1^{-{\mu_1}}}{(1-z\omega_1^{-1})^{\mu_1}}
+\cdots\\
&+&\frac{-c^{(h)}_1\omega_h^{-1}}{1-z\omega_h^{-1}}+\cdots+\frac{-c^{(h)}_{\mu_h}\omega_1^{-{\mu_h}}}{(1-z\omega_h^{-1})^{\mu_h}}
\end{eqnarray*}
for certain constants $a,b,c^{(i)}_1,\ldots,c^{(i)}_{\mu_i}$ with $i=1,\ldots,h$. By expanding in series we have
\[
g(z) =\sum_{n=0}^\infty\left(a -(be^\Lambda) e^{\Lambda n} + 
\sum_{i=1}^h P_i(n) \omega_i^{-n}\right) z^n ,
\]
where $P_i$ are polynomials of degree smaller than $\mu_i$. Consequently, the following holds
\[
G(n) =a - (be^\Lambda) e^{\Lambda n} + 
\sum_{i=1}^h P_i(n) \omega_i^{-n}.
\]
It remains to determine the constants $a,b$.
\begin{eqnarray*}
a&=&\lim_{z\to 1}(z-1)\cdot g(z)\\
 &=&\lim_{z\to 1}\frac{(z-1)}{(1-z)(1-z^{n_1}+\cdots-z^{n_m})}\\
 &=&-\frac{1}{1-m},
\end{eqnarray*} 
and
\begin{eqnarray*}
b&=&\lim_{z\to e^{-\Lambda}}(z-e^{-\Lambda})\cdot g(z)\\
 &=&\lim_{z\to e^{-\Lambda}}\frac{(z-e^{-\Lambda})}{(1-z)(1-z^{n_1}+\cdots-z^{n_m})}\\
 &=&\frac{1}{1-e^{-\Lambda}}\lim_{z\to e^{-\Lambda}}\left(\frac{1}{-n_1z^{n_1-1}+\cdots -n_mz^{n_m-1}}\right)\\
&=&\frac{-e^{-\Lambda}}{(1-e^{-\Lambda})(n_1e^{-n_1\Lambda}+\cdots+n_me^{-n_m\Lambda})}\\
&=&\frac{-\Lambda e^{-\Lambda}}{H(1-e^{-\Lambda})},
\end{eqnarray*}
where in the last equality we used that $n_j=\frac{\log\left(\frac{1}{p_j}\right)}{\Lambda}$ for $j=1,\ldots,m$ and the definition of entropy. Hence, it follows that 
$$
G(n) = \frac{\Lambda  e^{\Lambda n}}{H(1-e^{-\Lambda })} + 
\sum_{i=1}^h P_i(n) \omega_i^{-n} - \frac 1{m-1}.
$$

Note that in view of (\ref{eqMrrel}) the constant term $-1/(m-1)$ 
disappears when we translate the asymptotics of $G(n)$ to $M_r$. Next we study the error term (without the constant term $-1/(m-1)$) in more detail. W.l.o.g.\ we can assume that $\omega_1,\ldots,\omega_k$ (with $k\le h$) are those roots of $f(z)$ with smallest modulus \begin{equation}\label{etaRel}
    |\omega_i| = e^{-\Lambda (1-\eta)}
\end{equation} (for some $\eta > 0$) such that $P_i\ne 0$ and the degrees of $P_i$ are maximal and all equal to $d\ge 0$, for $1\le i \le k$. This means that the difference between $G(n)$ and the asymptotic leading term is bounded by
\[
\delta(n) = \left| G(n) - \frac{\Lambda  e^{\Lambda n}}{H(1-e^{-\Lambda }) }  + \frac 1{m-1} \right|= \left|\sum_{i=1}^h P_i(n)\omega_i^{-n}\right|
\le C n^d e^{\Lambda (1-\eta)n}
\]
for some constant $C>0$.
More precisely $\delta(n)$ can be written as
\[
\delta(n) = \left| n^d \sum_{i=1}^k \widetilde c_i \, \omega_i^{-n} \right| 
+\mathcal{O}\left( n^{d-1} e^{\Lambda (1-\eta)n} \right),
\]
with complex numbers $\widetilde c_i \ne 0$, $1\le i\le k$. Since all roots of $f(z)$
are either real or appear in conjugate pairs of complex numbers 
we can rewrite the sum $\sum\limits_{i=1}^k \widetilde c_i \, \omega_i^{-n}$ to
\[
n^d e^{\Lambda (1-\eta)n} \sum_{i=1}^{k'} c_i' \cos(2\pi \theta_i n + \alpha_i)
\]
with real numbers $c_i'\ne 0$ for $1\le i\le k'$. From Lemma~\ref{LeApp1} 
it follows that there exists $\delta > 0$ and infinitely many $n$ such that 
$\left|\sum\limits_{i=1}^{k'} c_i' \cos(2\pi \theta_i n + \alpha_i)\right| \ge \delta$.
This shows that 
\[
\delta(n) \ge C' n^d e^{\Lambda (1-\eta)n}
\]
for infinitely many $n$ and some constant $C'> 0$.
This means that the error term in~(\ref{ratCase}) is optimal.

Finally, we have 
$$
G(n) = \frac{\Lambda  e^{\Lambda n}}{H(1-e^{-\Lambda })}-\frac{1}{m-1} +\mathcal{O}\left(n^d e^{\Lambda (1-\eta)n}\right)
$$
for some $\eta > 0$. Obviously, this implies the representation (\ref{ratCase}) of Theorem~\ref{Thm1}. In fact, since $A(v)=G\left(\left[\frac{\log v}{\Lambda}\right]\right)$ we have
\begin{eqnarray*}
A(v)& =& \frac{\Lambda}{H(1-e^{-\Lambda })}\cdot e^{\Lambda \left[\frac{\log v}{\Lambda}\right]}-\frac{1}{m-1} +\mathcal{O}\left(\left[\frac{\log v}{\Lambda}\right]^d e^{\Lambda (1-\eta)\left[\frac{\log v}{\Lambda}\right]}\right)\\
&=& \frac{\Lambda v }{H(1-e^{-\Lambda })}\cdot e^{-\Lambda \left\{\frac{\log v}{\Lambda}\right\}}-\frac{1}{m-1} +\mathcal{O}\left((\log v)^d v^{(1-\eta)}\right)
\end{eqnarray*}
and so (\ref{eqMrrel}) implies that
$$
M_r=\frac{(m-1)}{rH}\cdot\frac{\Lambda}{(1-e^{-\Lambda })} 
e^{-\Lambda \left\{\frac{\log\left(\frac 1r\right) }{\Lambda}\right\}}
+\mathcal{O}\left((\log r)^d r^{-(1-\eta)}\right).
$$

\paragraph{\emph{Irrational case}}\ \\
The analysis in the irrational case is much more involved.
Instead of using power series we use the Mellin transform of the function $A(v)$ (see \cite{Szp}), i.e.
\[
A^*(s)= \int_0^\infty A(v) v^{s-1}\, dv.
\]
By using the fact that the Mellin transform of $A(av)$ is $a^{-s}A^*(s)$, 
a simple analysis of recurrence (\ref{eqrecA}) 
reveals that the Mellin transform $A^*(s)$ of $A(v)$ is given by
\begin{equation}\label{MellTransf}
A^*(s) = \frac{-1}{s \left( 1- p_1^{-s} - \cdots - p_m^{-s} \right) },
\qquad \Re(s) < -1.
\end{equation}
In fact, for $\Re(s)<-1$ we have that
\begin{eqnarray*}
A^*(s)&=&\int_0^\infty A(v)v^{s-1}dv=\int_1^\infty\left(1+\sum_{j=1}^mA(p_jv)\right)v^{s-1}dv\\
		&=&\int_1^\infty v^{s-1}dv+\sum_{j=1}^m\int_1^\infty A(p_jv)v^{s-1}dv \\
	&=&-\frac{1}{s}+\sum_{j=1}^mp_j^{-s}\int_1^\infty A(v)v^{s-1}dv\\
&=&-\frac{1}{s}+\left(\sum_{j=1}^mp_j^{-s}\right)\int_0^\infty A(v)v^{s-1}dv\\
&=&-\frac{1}{s}+\left(\sum_{j=1}^mp_j^{-s}\right) A^*(s),
\end{eqnarray*}
which implies the relation (\ref{MellTransf}).

In order to find asymptotics of $A(v)$ as $v\to \infty$ one can directly use the Tauberian theorem (for the Mellin transform) by Wiener-Ikehara \cite[Theorem 4.1]{Korevaar}. For this purpose we have to check that $s_0 = -1$ is the only (polar) singularity on the line $\Re(s) = -1$ and that $(s+1)A^*(s)$ can be analytically extended to a region that contains the line $\Re(s) = -1$. However, in the irrational case this follows by a lemma of Schachinger \cite{Schachinger}. In particular, one finds $$ A(v) \sim \frac vH$$ but this procedure does not provide any information about the error term.

For making our presentation as simple as possible we will restrict
ourselves to the case $m=2$ and we will also assume certain conditions
on the Diophantine properties of the irrational number
\[
\gamma=\frac{\log p_1 }{\log p_2}.
\]
We use the simplified notation $p=p_1$ and $q = p_2$.

The principle idea to obtain error terms for $A(v)$, in this case, is using the
formula for the inverse Mellin transfrom
\begin{equation}\label{invMell}
A(v)=\frac{1}{2\pi i}\lim_{T\to\infty}\int_{\sigma-iT}^{\sigma+iT}A^*(s)v^{-s}ds,
\qquad \sigma < -1,
\end{equation}
and shifting the line of integration to the right. Of course, all polar 
singularities of $A^*(s)$, which are given by the solutions of 
the equation $p^{-s} + q^{-s} = 1$ and $s=0$, give rise to 
a polar singularity of $A(v)$. Unfortunately, the order of magnitude of $A^*(s)$ 
is $\mathcal{O}(1/s)$. Hence the integral
in (\ref{invMell}) is not absolutely convergent. Therefore, it is convenient
to smooth the problem and to study the function $$A_1(v) = \int_0^v A(w)\, dw,$$ which is given by
\begin{eqnarray*}
A_1(v) &=&\frac{1}{2\pi i} \int_{\sigma-i\infty}^{\sigma+i\infty} A^*(s)\int_0^v w^{-s}\,dw\,ds\\
		&=& \frac{1}{2\pi i} \int_{\sigma-i\infty}^{\sigma+i\infty}
A^*(s)\cdot\frac{v^{-s+1}}{1-s}\,ds  \\
&=& \frac{1}{2\pi i} \int_{\sigma-i\infty}^{\sigma+i\infty}
\frac{v^{-s+1}}{s(s-1) (1-p^{-s}- q^{-s})}\,ds
,
\qquad \sigma < -1.
\end{eqnarray*}
By Proposition \ref{SchachProp} we know that all zeros of the equation
$p^{-s} + q^{-s} = 1$ that are different from $-1$ satisfy 
$-1 < \Re(s) \le \sigma_0$ for some $\sigma_0$. Furthermore,
there exists $\tau> 0$ such that in each box of the form
\[
B_k = \{s\in \mathbb{C} : -1 < \Re(s) \le \sigma_0,\, 
(2k-1)\tau \le \Im(s) < (2k+1)\tau \}, \quad k \in 
\mathbb{Z}\setminus \{ 0 \},
\]
there is precisely one zero of $p^{-s} + q^{-s} = 1$ that we denote by $s_k$.

Now, in the evaluation of $A_1(v)$ we shift the line of integration to the right, namely to $\Re(s)=\sigma_1$ with $\sigma_1 > \max\{\sigma_0+1,1\}$. Let us denote by $\mathcal{S}$ the set of all the singularities $s'\in\mathbb{C}$ of the integrand function $g(s):=A^*(s)\frac{v^{-s+1}}{1-s}$ such that $\sigma\leq\Re(s')\leq\sigma_1$, i.e.
$$\mathcal{S}=\left\{-1,\ 0,\ 1,\ s_k \ \text{for} \ k\in\mathbb{Z}\setminus \{ 0 \}\right\}.$$ 

Then by applying Cauchy's residue theorem, we get
\begin{equation}\label{ResThm}
A_1(v)=\frac{1}{2\pi i} \int_{\sigma-i\infty}^{\sigma+i\infty}
g(s)\,ds=-\sum_{s'\in\mathcal{S}} Res\left(g(s),\ s=s'\right) + \frac{1}{2\pi i} \int_{\sigma_1-i\infty}^{\sigma_1+i\infty}
g(s)\,ds.
\end{equation}
So, we have to consider the following residues:
\begin{itemize}
\item $Res\left(g(s),\ s=-1\right)=-\frac{v^2}{2H}$ 
\item $Res\left(g(s),\ s=s_k\right)=\frac{v^{1-s_k}}{s_k(s_k-1)\left(p^{-s_k} \log\left(\frac 1p\right) + q^{-s_k} \log\left(\frac 1q\right)\right)},\ \text{for}\ k\in\mathbb{Z}\setminus \{ 0 \}$
\item $Res\left(g(s),\ s=0\right)=v$
\item $Res\left(g(s),\ s=1\right)=\frac{1}{1-p^{-1}-q^{-1}}$.
\end{itemize}
Hence, by collecting all residues and using (\ref{ResThm}) we obtain
\begin{align*}
A_1(v) &= \frac{v^2}{2H} - \sum_{k\in \mathbb{Z}\setminus \{ 0 \} } 
\frac{v^{1-s_k}}{s_k(s_k-1)H(s_k)} 
-  v - \frac 1{1 - p^{-1} - q^{-1}}  \\
&+\frac{1}{2\pi i} \int_{\sigma_1-i\infty}^{\sigma_1+i\infty}
\frac{v^{-s+1}}{s(s-1) (1-p^{-s}- q^{-s})}\,ds,
\end{align*}
where we put $H(s) = p^{-s} \log\left(\frac 1p\right) + q^{-s} \log\left(\frac 1q\right)$.

Now, it is easy to see that the integral can be estimated by 
\begin{equation}\label{Int}
\frac{1}{2\pi i} \int_{\sigma_1-i\infty}^{\sigma_1+i\infty}
\frac{v^{-s+1}}{s(s-1) (1-p^{-s}- q^{-s})}\,ds = 
\mathcal{O}\left( v^{-\sigma_1 +1} \right).
\end{equation}
In fact, we can write
\begin{eqnarray*}
&&\frac{1}{2\pi i} \int_{\sigma_1-i\infty}^{\sigma_1+i\infty}
\frac{v^{-s+1}}{s(s-1) (1-p^{-s}- q^{-s})}\,ds \\
&=&\frac{1}{2\pi} \int_{-\infty}^{\infty}
\frac{v^{-\sigma_1-it+1}}{(\sigma_1+it)(\sigma_1+it-1) (1-p^{-\sigma_1-it}- q^{-\sigma_1-it})}\,dt. 
\end{eqnarray*}
Since
\begin{eqnarray*}
&&\left|\frac{v^{-\sigma_1-it+1}}{(\sigma_1+it)(\sigma_1+it-1) (1-p^{-\sigma_1-it}- q^{-\sigma_1-it})}\right|\\
&=&\frac{|v^{-\sigma_1+1}|\cdot|v^{-it}|}{\sqrt{\sigma_1^2+t^2}
\cdot\sqrt{(\sigma_1-1)^2+t^2}|1-p^{-\sigma_1-it}- q^{-\sigma_1-it}|}\\
&\leq& c\frac{v^{-\sigma_1+1}}{(1+t^2)}
\end{eqnarray*}
for some constant $c>0$, then we have
\begin{eqnarray*}
\left|\frac{1}{2\pi} \int_{-\infty}^{\infty}
\frac{v^{-\sigma_1-it+1}}{(\sigma_1+it)(\sigma_1+it-1) (1-p^{-\sigma_1-it}- q^{-\sigma_1-it})}\,dt\right|&\leq& \frac{c}{2\pi} \int_{-\infty}^{\infty}\left(\frac{v^{-\sigma_1+1}}{1+t^2}\right)\,dt\\
&=&\frac{c}{2}\cdot v^{-\sigma_1+1} 
\end{eqnarray*}
and so the relation (\ref{Int}) holds. Hence, we just have to deal with the sum of residues 
\begin{equation}\label{sumRes}
    \sum\limits_{k\in \mathbb{Z}\setminus \{ 0 \} } 
\frac{v^{1-s_k}}{s_k(s_k-1)H(s_k)}.
\end{equation}

First, let us show that there exists $\delta > 0$ such that $|H(s_k)|\ge \delta$
for all $k\in \mathbb{Z}\setminus \{ 0 \}$. Without loss of generality we can assume that $q<p$ and so $\frac{\log q}{\log p}<1$. Since $s_k\in\mathcal{S}$ we have that
\begin{eqnarray*}
H(s_k)&=&p^{-s_k}\log\left(\frac 1p\right)+q^{-s_k}\log\left(\frac 1q\right)\\
&=&(1-q^{-s_k})\log\left(\frac 1p\right)+q^{-s_k}\log\left(\frac 1q\right)\\
&=&\log\left(\frac 1p\right)\left(1-q^{-s_k}\left(1-\frac{\log q}{\log p}\right)\right).
\end{eqnarray*}
So it follows that
$$|H(s_k)|=\log\left(\frac 1p\right)\left|1-q^{-s_k}\left(1-\frac{\log q}{\log p}\right)\right|\geq\log\left(\frac 1p\right)\left|1-|q^{-s_k}|\cdot\left|1-\frac{\log q}{\log p}\right|\right|.$$
Therefore, since $\left(1-\frac{\log q}{\log p}\right)<1$ and $\Re(s_k)<1$, there exists $c_0>0$ such that
$$|H(s_k)|\geq c_0\log\left(\frac 1p\right)=\delta.$$
Thus, we do not have to care about this factor in (\ref{sumRes}).

Next assume that $\gamma$ is a badly approximable irrational number.
Here Lemma~\ref{freeZeroReg} shows that all zeros $s_k\in\mathcal{S}$ satisfy 
$\Re(s_k) > -1+ c/\Im(s_k)^2$ for some constant $c> 0$. Hence
it follows that $\Re(s_k) > -1 + c_1/k^2$ for some constant $c_1 > 0$
and we can estimate the sum of residues by
\begin{align*}
\left|\sum_{k\in \mathbb{Z}\setminus \{ 0 \} } 
\frac{v^{1-s_k}}{s_k(s_k-1)H(s_k)}  \right| 
&\le\left|\sum_{0< |k| \le K } 
\frac{v^{1-s_k}}{s_k(s_k-1)H(s_k)}  \right| +
\left|\sum_{|k| > K} 
\frac{v^{1-s_k}}{s_k(s_k-1)H(s_k)}  \right| \\
&\le\sum_{0< |k| \le K } 
\frac{v^{1-\Re(s_k)}}{\left|s_k(s_k-1)\right|\delta} +
\sum_{|k| > K} 
\frac{v^{1-\Re(s_k)}}{\left|s_k(s_k-1)\right|\delta}\\
&\le C_1 v^{2- c_1/K^2} \sum_{0< |k| \le K}  \frac 1{k^2}  
+ C_2 v^2 \sum_{|k| > K}  \frac 1{k^2} \\
&\le C_3 v^2 \left( v^{-c_1/K^2} + \frac 1K  \right).
\end{align*}
where $C_1,C_2,C_3$ are appropriate positive constants.

Thus, by choosing $K = \sqrt{c_1(\log v)/(\log\log v)}$, we obtain the
upper bound
\[
\sum_{k\in \mathbb{Z}\setminus \{ 0 \} } 
\frac{v^{1-s_k}}{s_k(s_k-1)H(s_k)}  
= \mathcal{O}\left( v^2 \frac{\sqrt{\log\log v}}{\sqrt{\log v}} \right)
\]
and consequently
\[
A_1(v) = \frac{v^2}{2H} \left( 1 + 
\mathcal{O}\left( \frac{\sqrt{\log\log v}}{\sqrt{\log v}} \right) \right).
\]
Finally, by an application of Lemma~\ref{Taub}, the previous relation implies
\[
A(v) = \frac vH \left( 1 + \mathcal{O}
\left( \frac{ (\log\log v)^{1/4} }{ (\log v)^{1/4}}\right) \right).
\]

Similarly, when $p$ and $q$ are algebraic, we deal with the case when all solutions
of the equation $p^{-s} + q^{-s} = 1$ (that are different from $-1$)
satisfy $\Re(s_k) > -1 + \frac{D}{\Im(s_k)^{2C}}$ for some positive constants
$C,D$ (see Lemma~\ref{freeZeroRegBak}). Then with the same procedure as above we get
\[
\left|\sum_{k\in \mathbb{Z}\setminus \{ 0 \} } 
\frac{v^{1-s_k}}{s_k(s_k-1)H(s_k)}  \right| 
\le C_4 v^2 \left( v^{-c_2K^{-2C}} + \frac 1K  \right)
\]
for some constant $C_4>0$.
Hence, if we choose $K = (c_2 (\log v)/(\log\log v))^{1/(2C)}$,
we obtain (after a second application of Lemma~\ref{Taub})
\[
A(v) = \frac vH \left( 1 + \frac{ (\log\log v)^{\kappa} }{ (\log v)^{\kappa}} \right).
\]
where $\kappa=\frac{1}{4C}$.
This completes the proof of Theorem~\ref{Thm1}.
\\
\qed

\subsection{Discrepancy bounds in the rational case}
First, we want to study the rational case. So in this subsection we are going to consider a partition $\rho$ of $[0,1]$ consisting of $m$ intervals of lengths $p_1,\ldots,p_m$ such that $\log\left(\frac{1}{p_1}\right),\ldots,\log\left(\frac{1}{p_m}\right)$ are rationally related.

By Theorem \ref{Thm1} we know that asymptotically
\begin{equation}\label{ratCase2}
M_{r_n} =  \frac{c'}{r_n} 
+ \mathcal{O}\left((\log r_n)^d r_n^{-(1-\eta)}\right), \quad 
r_n = e^{-\Lambda n},
\end{equation}
for some $\eta > 0$ and some integer $d\ge 0$, 
where $c' = \frac{{(m-1)\Lambda }}{H(1-e^{-\Lambda })}$ 
and the error term is optimal. 
Recall also that $k(n) = M_{r_{n}}$, which gives an 
asymptotic expansion for $k(n)$ of the form
\begin{equation}\label{relKappan}
k(n) \sim \frac{(m-1)\Lambda }{H(1- e^{-\Lambda})}\, e^{\Lambda n}.
\end{equation}

\begin{theorem}\label{disRatCase}\ \\
Suppose that the lengths of the intervals of a partition $\rho$
are $p_1,\ldots,p_m$ and assume that  
$\log\left(\frac{1}{p_1}\right),\ldots,\log\left(\frac{1}{p_m}\right)$ are rationally related. Furthermore, let $\eta > 0$ and $d\ge 0$ be given as 
in Theorem~\ref{Thm1}. Then the discrepancy of the sequence of partitions $(\rho^n\omega)$ is bounded by
\begin{equation}\label{uppBound}
D_n= \left\{ \begin{array} {ll} 
\mathcal{O}\left( (\log k(n))^d k(n)^{-\eta} \right) & \mbox{if $0< \eta < 1$,}\\
\mathcal{O}\left( (\log k(n))^{d+1} k(n)^{-1} \right) & \mbox{if $\eta = 1$,}\\
\mathcal{O}\left( k(n)^{-1} \right) & \mbox{if $\eta > 1$.}
\end{array} \right.
\end{equation}
Moreover, there exist $\delta> 0$ and infinitely many $n$ such that
\begin{equation}\label{lowBound}
D_n\ge \left\{ \begin{array} {ll} 
\delta \,(\log k(n))^d k(n)^{-\eta} & \mbox{if $0< \eta < 1$,}\\
\delta \, (\log k(n))^{d} k(n)^{-1} & \mbox{if $\eta = 1$,}\\
\delta \, k(n)^{-1} & \mbox{if $\eta > 1$.}
\end{array} \right.
\end{equation}
\end{theorem}
\proof\ \\
For notational convenience we set
\[
\Delta_n = \sup_{0< y \le 1} \left| 
\sum_{i=1}^{k(n)}
\chi_{[0,y[}\left(t_i^{(n)}\right)- k(n)y
 \right|,
\]
where $t_i^{(n)}$ are the points defining the partition $\rho^n\omega$. Then we have $D_n \le 2 \Delta_{n} / k(n)$, since $\Delta_{n} / k(n)= D_n^*$ and Theorem \ref{relDisStarDis} holds.

Fix a step in the algorithm corresponding to a certain parameter $r$ of the form $r=e^{-n\Lambda }$ for some integer $n\ge 0$, and consider an interval $A=[0,y[\subset[0,1]$. We want to estimate the number of elementary intervals belonging to $\mathscr{E}_r$ which are contained in $A$. For this purpose, let us fix another parameter $\overline{r}$ of the form $\overline{r}=e^{-\overline{n}\Lambda }$ with an integer $0\le \overline{n}\le n$ corresponding to a previous step in Khodak's construction. At this previous step, we have $M_{\overline{r}}$ intervals $I_j$ generated by the construction. Now, the lenghts of the intervals $I_j$ are given by $\lambda(I_j)$ and we have that 
\begin{equation}\label{lj}
p_{min}\overline{r}\leq \lambda(I_j)<\overline{r}  \quad \textrm{for}\quad j=1,\ldots,M_{\overline{r}},
\end{equation}
since the lengths of the intervals $\mathscr{E}_{\overline r}$ correspond to the values $P(y)$ of the external nodes $y$ in $\mathcal{E}({\overline r})$.

Suppose that precisely the first $h$ of these intervals $I_j$ are contained in $A$, so
$U=I_1\cup\ldots\cup I_h\subset A$. Now, we want to estimate the number of elementary intervals in $\mathscr{E}_r$ contained in each~$I_j$. Khodak's construction shows that this equals precisely the number of external nodes in the subtree of the node $x$ that is related to the interval~$I_j$. An important feature of Khodak's construction is that subtrees of $\mathcal{T}(r)$ rooted at an internal node $x\in\mathcal{I}(r)$ are parts of a self-similar infinite tree and therefore they are constructed in the same way as the whole tree. So, one just has to replace~$r$ by~$\frac{r}{P(x)}$. Hence, by using this remark in (\ref{ratCase2}), the number $N_{I_j}$ of subintervals of $I_j$ (corresponding to the value $r$) equals 
\[
N_{I_j}=M_{\frac{r}{\lambda(I_j)}}=\frac{c'}{r}\lambda(I_j)+
\mathcal{O}\left( |\log r|^d \frac{\lambda(I_j)^{1-\eta}}{r^{1-\eta}}\right).
\]
Therefore, we have that the number $N_U$ of elementary intervals in $\mathscr{E}_r$ contained in $U$ is 
\[
N_U=N_{I_1}+\ldots+N_{I_h}=\frac{c'}{r}(\lambda(I_1)+\ldots+\lambda(I_h))
+\mathcal{O}\left(   \frac{|\log r|^d }{r^{1-\eta}}  
{\sum\limits_{j=1}^h   {\lambda(I_j)}^{1-\eta}}\right).
\]
By using (\ref{lj}) and the fact that $h\leq  M_{\overline{r}} = \mathcal{O}(1/\overline r)$  we obtain
\begin{align*}
N_U &= \frac{c'}{r}(\lambda(I_1)+\ldots+\lambda(I_h))+\mathcal{O}\left( |\log r|^d \frac{{h\overline{r}}^{(1-\eta)}}{r^{(1-\eta)}}\right)\\
&= \frac{c'}{r}(\lambda(I_1)+\ldots+\lambda(I_h))
+\mathcal{O}\left( |\log r|^d \frac{{\overline{r}}^{(-\eta)}}{r^{(1-\eta)}}\right).
\end{align*}
Since the total number of intervals equals $M_r = {c'}/{r} + 
\mathcal{O}(|\log r|^d r^{-1+\eta})$ it follows that
\[
N_U - M_r \lambda(U) = \mathcal{O}\left( |\log r|^d \frac{{\overline{r}}^{(-\eta)}}{r^{(1-\eta)}}\right) + \mathcal{O}\left(\frac{|\log r|^d}{ r^{1-\eta}}\right) =
\mathcal{O}\left( |\log r|^d \frac{{\overline{r}}^{(-\eta)}}{r^{(1-\eta)}}\right).
\]
Since $N_A - M_r \lambda(A) = (N_U - M_r \lambda(U)) + (N_{A\setminus U} - 
M_r \lambda(A \setminus U) )$ it remains to study the difference
\begin{align*}
N_{A\setminus U} - M_r \lambda(A \setminus U) &= \left(
N_{A\setminus U} - M_{r/\lambda(I_{h+1})} \frac{\lambda(A \setminus U)}{\lambda(I_{h+1})}\right) \\
&+\left(M_{r/\lambda(I_{h+1})} \frac{\lambda(A \setminus U)}{\lambda(I_{h+1})} - M_r \lambda(A \setminus U)\right).
\end{align*}
The second term can be directly estimated by
\[
\left|M_{r/\lambda(I_{h+1})} \frac{\lambda(A \setminus U)}{\lambda(I_{h+1})} - M_r \lambda(A \setminus U) \right|
= \mathcal{O}\left( |\log r|^d \,\frac{{\overline{r}}^{(1-\eta)}}{r^{(1-\eta)}}     \right),
\]
whereas the first term is bounded by
\[
\left| N_{A\setminus U} - M_{r/\lambda(I_{h+1})} \frac{\lambda(A \setminus U)}{\lambda(I_{h+1})}
\right| \le \Delta_{n-\overline n}
\]
since $\frac{r}{{\overline{r}}}=e^{-\Lambda(n-{\overline{n}})}$.

Summing up and taking the supremum over all sets $A = [0,y[$, we obtain the recurrence relation
\begin{equation*}\label{eqDeltarec}
\Delta_n \le \Delta_{n-\overline n} +
 \mathcal{O}\left( |\log r|^d \frac{{\overline{r}}^{(-\eta)}}{r^{(1-\eta)}}\right).
\end{equation*}
We now set $\overline n = 1$ and recall that 
$r = e^{-\Lambda n}$ and also $\overline r = e^{-\Lambda \overline n} = e^{-\Lambda }$.
Thus, by the previous relation we get
\begin{equation}\label{eqDeltarec-2}
\Delta_n \le \Delta_{n-1} +
 \mathcal{O}\left( n^d e^{\Lambda n(1-\eta)} \right).
\end{equation}
We distinguish between three cases.
\begin{enumerate}
\item $0< \eta < 1$. In this case we get
\[
\Delta_n = \mathcal{O}\left( \sum_{k\le n} k^d e^{\Lambda k(1-\eta)} \right) 
= \mathcal{O}\left( n^d e^{\Lambda n(1-\eta)} \right).
\]
By taking into account also the relation (\ref{relKappan}), it follows that 
$$D_n \leq 2\frac{\Delta_n}{k(n)}=\mathcal{O}\left( (\log k(n))^d \frac{k(n)^{(1-\eta)}}{k(n)} \right)=\mathcal{O}\left( (\log k(n))^d k(n)^{-\eta} \right).$$
\item $\eta = 1$. In this case we get 
$$\Delta_n =\mathcal{O}\left( \sum_{k\le n} k^d \right)= \mathcal{O}( n^{d+1})$$ and
consequently 
$$D_n \leq 2\frac{\Delta_n}{k(n)} = \mathcal{O}\left( (\log k(n))^{d+1} k(n)^{-1} \right).$$
\item $\eta > 1$. Here we have
\[
\Delta_n = \mathcal{O}\left( \sum_{k\le n} k^d e^{-\Lambda k(\eta-1)} \right) 
= \mathcal{O}\left( 1 \right)
\]
which rewrites to $D_n =\mathcal{O}\left(k(n)^{-1} \right)$.
\end{enumerate}

In order to give a lower bound of the discrepancy it is sufficient to handle the case $0< \eta \le 1$. In fact, if $\eta > 1$ we just use the trivial lower bound $D_n \ge \frac{ 1}{k(n)}$
which meets the upper bound. For the remaining case $0< \eta \le 1$ we consider the interval $A = [0,p_1[$. We also recall (see the proof of Theorem \ref{Thm1}) that we can write $M_r$, for $r = r_n = e^{-\Lambda n}$, as
\[
M_r = c'\,e^{\Lambda n} + \delta_n
\]
where $\delta_n$ has a representation of the form
\[
\delta_n =  n^d e^{\Lambda n(1-\eta)} \sum_{i=1}^k c_i \cos(2\pi \theta_i n + \alpha_i)
+ \mathcal{O}\left( n^{d-1} e^{\Lambda n(1-\eta)} \right).
\]
Similarly to the above we obtain 
\begin{align*}
N_A - M_r \lambda(A) &= M_{r/p_1} - M_r p_1 \\
&=c'\,e^{\Lambda n}p_1 + \delta_{n-n_1}-p_1\left(c'\,e^{\Lambda n} + \delta_n\right)\\
&= \delta_{n-n_1} - p_1 \delta_{n} \\
&= n^d e^{\Lambda n(1-\eta)} \left( 
\sum_{i=1}^k c_i \cos(2\pi \theta_i n + \alpha_i-2\pi \theta_in_1)\right.\\
&-\left.p_1 \sum_{i=1}^k c_i \cos(2\pi \theta_i n + \alpha_i) \right)+  \mathcal{O}\left( n^{d-1} e^{\Lambda n(1-\eta)} \right).
\end{align*}
By applying Lemma~\ref{LeApp1} it follows that there 
exist $\delta > 0$ and infinitely many
$n$ with
\[
|N_A - M_r \lambda(A)| \ge \delta n^d e^{\Lambda n(1-\eta)}.
\]
Consequently
\[
D_n \ge \frac 1{M_r} |N_A - M_r \lambda(A)| \ge \frac{\delta n^d e^{\Lambda n(1-\eta)}}{e^{\Lambda n}}= \delta n^d e^{-\Lambda n\eta}
\]
for some $\delta> 0$.
This completes the proof of the lower bound (\ref{lowBound}).\\
\endproof
 
\subsection{Discrepancy bounds in the irrational case}\label{subsec: Irr}
As mentioned above, the case when $\log\left(\frac{1}{p_1}\right),\ldots,\log\left(\frac{1}{p_m}\right)$ are irrationally related is much more difficult to handle since the error term in the asymptotic expansion for $M_r$ is not explicit in general (see (\ref{irratCase}) in Theorem \ref{Thm1}). Nevertheless, we can provide
upper bounds in some cases of interest.

Suppose that $m=2$, set $p = p_1$ and $q = p_2$ and $\gamma = (\log p)/(\log q)$.
First, let us show that 
\begin{prop}\ \\
The number of intervals $k(n)$ of the partition $\rho^n\omega$ is asymptotically given by
\[
k(n) \sim \left(\frac{m-1}{H}\right) \exp\left( \sqrt{2 n \log\left( \frac 1p\right)\,  \log\left( \frac 1q\right)}\right).
\]
\end{prop}
\proof\ \\
Let $r$ be the parameter in Khodak's construction that corresponds to the step $n$, then $M_r=k(n)$.
By (\ref{irratCase}) in Theorem \ref{Thm1} we have that 
\begin{equation}\label{irrAsymp}
k(n)\sim \frac{m-1}{H}\cdot\frac 1r.
\end{equation}
Note that there is a one-to-one correspondence between the probability of each node $x$, that is $P(x)=p^kq^l$ and the non-negative integral lattice points $(k,\ell)$. So the number $n$ 
of steps corresponding to the value $r$ is approximatively given by the cardinality of the set  
\begin{eqnarray*}
\{x\in\mathcal{T}(r): P(x)\geq r\}&=&\{(k,l): p^kq^l\geq r\}\\
&=&\left\{(k,l):\frac{1}{p^kq^l}\leq\frac{1}{r}\right\}\\
&=&\left\{(k,l):k\log\left(\frac{1}{p}\right)+l\log\left(\frac{1}{q}\right)\leq\log\left(\frac{1}{r}\right)\right\}
\end{eqnarray*}
Now, the equation $k\log p + \ell \log q = \log r$ has at most one solution in integer pairs $(k,\ell)$. Hence, we have
$$n\sim\frac{1}{2}\cdot\frac{\log\left(\frac{1}{r}\right)}{\log\left(\frac{1}{p}\right)}
\cdot\frac{\log\left(\frac{1}{r}\right)}{\log\left(\frac{1}{q}\right)}
=\frac{\left(\log\left(\frac{1}{r}\right)\right)^2}{2\log\left(\frac{1}{p} \right)\log\left(\frac{1}{q}\right)}$$
and so
$$
\frac 1r \sim \exp\left( \sqrt{2 n \log \left(\frac{1}{p}\right)\,  \log\left( \frac{1}{q} \right) }\right).
$$
The conclusion follows by using this relation in (\ref{irrAsymp}).
\endproof

In Theorem \ref{Thm1} we have considered the case when $\gamma$ is badly approximable and the case when $p$ and $q$ are algebraic. By using these results we can show the following theorem for the discrepancy in the irrational case. 
\begin{theorem}\label{EstBadlyApprx}\ \\
Suppose that the lenghts of the intervals of a partition $\rho$ of $[0,1]$ are $p$ and $q=1-p$ and let $\gamma=\frac{\log p}{\log q}$. If $\gamma\notin\mathbb{Q}$ and it is badly approximable, then the discrepancy of $(\rho^n\omega)$ is bounded by
\[
D_n=\mathcal{O}\left(\left(\frac{\log\log{(k(n))}}{\log{(k(n))}}\right)^\frac{1}{4}\right),\quad\textrm{as}\ n\to \infty.
\]
Furthermore, if $p$ and $q$ are algebraic and $\gamma\notin\mathbb{Q}$ then 
\[
D_n=\mathcal{O}\left(\left(\frac{\log\log{(k(n))}}{\log{(k(n))}}\right)^\kappa\right),\quad\textrm{as}\ n\to \infty,
\]
where $\kappa > 0$ is an effectively computable constant (see Theorem~\ref{Thm1}).
\end{theorem}

\proof\ \\
We use a procedure similar to the proof of Theorem~\ref{disRatCase} but now we consider the asymptotic expansion
\[
M_r = \frac{c''}r +\mathcal{O}\left( \frac 1r
\left(  \frac{\log\log \frac 1r}{\log \frac 1r} \right)^\xi
  \right)
\]
where $c'' = (m-1)/H$. Moreover, we have that $\xi=\frac 14$ when $\gamma$ is badly approximable, while $\xi=\kappa$ when $p,q$ are algebraic (see Theorem~\ref{Thm1}).

Fix a step $n$ in the algorithm corresponding to a certain parameter $r$ and consider an interval $A=[0,y[\subset[0,1]$. We want to estimate the number of elementary intervals belonging to $\mathscr{E}_r$ which are contained in $A$. For this purpose, let us fix another parameter $\overline{r}$ corresponding to the step $\overline{n}$ in Khodak's construction such that $0\le \overline{n}\le n$. At this previous step, we have $M_{\overline{r}}$ intervals $I_j$ generated by the construction. Now, the lenghts of the intervals $I_j$ are given by $\lambda(I_j)$ and we have that the relation (\ref{lj}) holds.

Suppose that precisely the first $h$ of these intervals $I_j$ are contained in $A$, so
$U=I_1\cup\ldots\cup I_h\subset A$. We want to estimate the number of elementary intervals in $\mathscr{E}_r$ contained in $I_j$. Similarly to the rational case we have that
the number $N_{I_j}$ of subintervals of $I_j$ (corresponding to the value~$r$) equals
\[
N_{I_j}=M_{\frac{r}{\lambda(I_j)}}=\frac{c''}r\lambda(I_j) +\mathcal{O}\left( \frac {\lambda(I_j)}r
\left(  \frac{\log\log \left(\frac {\lambda(I_j)}r\right)}{\log \left(\frac {\lambda(I_j)}r\right)} \right)^\xi
  \right).
\]
By using (\ref{lj}) and the fact that $h\leq M_{\overline{r}}=\mathcal{O}\left(\frac 1{\overline{r}}\right)$, it follows that

\begin{align*}
N_U - M_r \lambda(U) &= \frac{c''}{r}(\lambda(I_1)+\ldots+\lambda(I_h))
+\mathcal{O}\left( \frac {h\overline r}{r} 
\left(  \frac{\log\log \frac {\overline r}r}
{\log \frac {\overline r}r} \right)^\xi   \right) \\
&- \frac{c''}{r}(\lambda(I_1)+\ldots+\lambda(I_h)) +
\mathcal{O}\left( \frac 1r 
\left(  \frac{\log\log \frac 1r}{\log \frac 1r} \right)^\xi
  \right) \\
&= \mathcal{O}\left( \frac 1r 
\left(  \frac{\log\log \frac 1r}{\log \frac 1r} \right)^\xi
  \right).
\end{align*}
For the remaining interval $A\setminus U$ we use the bounds
$N_{A\setminus U} \le M_{r/\lambda(I_{h+1})} = \mathcal{O}(\overline r/r)$ 
and $\lambda(I_{h+1}) = \mathcal{O}(\overline r)$
to end up with the upper bound
\[
D_n = \mathcal{O}\left( \left(  \frac{\log\log \frac 1r}{\log \frac 1r} \right)^\xi
  \right) + \mathcal{O}\left( \overline r \right).
\]
Hence, by choosing 
\[
\overline r = \left(  \frac{\log\log \frac 1r}{\log \frac 1r} \right)^{\xi}
\]
we finally obtain
\[
D_n = \mathcal{O}\left(
\left(\frac{\log\log \frac 1r}{\log \frac 1r} \right)^{\xi} \right).
\]
This completes the proof of the theorem.\\
\endproof
Note that the upper bounds for the discrepancy we obtained are worse than $k(n)^{-\beta}$ for any $\beta>0$.  Actually, it seems that we cannot do really better in the irrational case.
This is due to the fact that $\liminf\limits_{k\ne 0} \Re(s_k) = -1$
where $s_k$, $k\ne 0$,  runs through all the zeros of the equation
$p^{-s} + q^{-s} = 1$ different from $s_0 = -1$.
Indeed, it seems that the continued fractional expansion of $\gamma = (\log p)/(\log q)$ could be used to obtain more explicit upper bounds. However, since they are all rather poor it 
is probably not worth working them out in detail. 

Moreover, the case $m> 2$ is even more involved, as we can see by comparing with the discussion of \cite{Flajolet-Vallee}.
In fact, in this paper the authors study the asymptotic structure of the main parameters of interest for digital trees. They represent collections of words over some finite alphabet, so the parsing tree constructed by Khodak's algorithm can be included in this class. In particular, in this paper digital trees are assumed to be under the simplest of all probabilistic models; namely, the memoryless source, where letters of words are drawn independently according to a fixed distribution. If $\mathcal{A}=\{a_1,\ldots,a_m\}$ is the alphabet, the model is determined by the basic quantities $p_j=P(a_j)$ with $p_1+\cdots+p_m=1$. As it turned out by an analysis based on Mellin transform, quantifying the main parameters of digital trees (such us expected number of internal nodes, expected path lenght, etc.) is strongly dependent on the location of poles in the complex plane of the fundamental Dirichlet series associated with the $p_j$'s, which is given by $$\Delta(s)=\frac{1}{1-p_1^s-\cdots-p_m^s}.$$ Neverthless, the results obtained in this paper relatively to the aperiodic case, which corresponds to our irrational one, show how the geometry of the set of poles of $\Delta(s)$ depends on the approximation properties of the ratios $\frac{\log p_i}{\log p_j}$ and therefore how much complicated is the study of this case.

\section{Applications}
In this section, we intend to present some examples and applications of the results of this chapter. In particular, we want to stress the application to fractals because 
the technique introduced in this chapter on $[0,1]$ allows to get discrepancy bounds for the elementary discrepancy of u.d.\! sequence of partitions of fractals belonging to a class wider than the one considered in Section~\ref{UDFractals}.
\subsection{Kakutani's sequences}
The procedure introduced in the Section \ref{sec:Discr} can be used to obtain bounds of the discrepancy for a family of classical Kakutani's sequences of partitions. In fact, if we fix $\alpha\in]0,1[$ then the corresponding Kakutani's sequence of partitions of $[0,1]$ is constructed by successive $\alpha$-refinements of the trivial partition $\omega=\{[0,1]\}$. So according to the notation used in the previous sections, in this case we have that $p_1=\alpha$ and $p_2=1-\alpha$.

For a Kakutani's sequence of parameter $\alpha$ we have that $\log\left(\frac{1}{\alpha}\right)$ and $\log\left(\frac{1}{1-\alpha}\right)$ are rationally related if and only if $\frac{\log\alpha}{\log{(1-\alpha)}}\in\mathbb{Q}\ $ (see Definition \ref{defRR}).

Let us denote by $\alpha_{n,m}$ the unique solution in $]0,1[$ of the following equation
$$\frac{\log x }{\log{(1-x)}}=\frac{n}{m}$$
with $\ n,m\in\mathbb{N}\ $. Since the function
$$f(x)=\frac{\log(x)}{\log{(1-x)}}$$ 
is continuous and strictly decreasing on $]0,1[$ and it attains all positive values, the countable set
$$\left\{\alpha_{n,m}=f^{-1}\left(\frac{n}{m}\right),\ n,m\in\mathbb{N}\right\}$$
is dense in $]0,1[$. 

The density of the values of the parameter $\alpha$ for which $\log\left(\frac{1}{\alpha}\right)$ and $\log\left(\frac{1}{1-\alpha}\right)$ are rationally related shows that we have interesting bounds of the discrepancy for a countable set of Kakutani's sequences (see Theorem \ref{disRatCase}). On the other hand, there are much more values of $\alpha$ in $]0,1[$ for which $\log\left(\frac{1}{\alpha}\right)$ and $\log\left(\frac{1}{1-\alpha}\right)$ are irrationally related and for which the discrepancy bounds are weaker. Neverthless, our technique allows to get quantitive results about the discrepancy of a large class of Kakutani's sequences, not known in the existing literature.

\subsection{$LS$-sequences}
$LS$-sequences are a special class of sequences of partitions constructed by successive $\rho-$refinements of the trivial partition $\omega = [0,1]$. We have already introduced these sequences in Section \ref{sec:Discr}, but let us recall their definition.
\begin{definition}\ \\
Fixed two positive integers $L$ and $S$, let $0<\alpha<1$ be the real number given by the equation $L\alpha + S\alpha^2 = 1$. The $LS$-sequence is the sequence of partitions obtained by successive $\rho$-refinements of $\omega$, when $\rho$ consists of $L$ subintervals of $[0,1]$ of length~$\alpha$ and $S$ subintervals
of length~$\alpha^2$.
\end{definition}
For instance, if $L = S = 1$ then $\alpha=\frac{\sqrt{5}-1}{2}$
and we obtain the so-called \emph{Kakutani-Fibonacci sequence}. This term was used first in \cite{Car}, since this is a particular Kakutani's sequence and the sequence $(k(n))_{n\in\mathbb{N}}$ of the number of intervals of the $n-$th partition is the sequence of Fibonacci numbers. The Kakutani-Fibonacci sequence is important because it is the only Kakutani's sequence for which the exact discrepancy is known (apart from the trivial Kakutani's sequence with $\alpha=\frac 12$). Here we have $p_1=\alpha$ and $p_2=1-\alpha=\alpha^2$ and consequently, using our techinique we have
$$
\log\left(\frac{1}{\alpha}\right)=n_1\Lambda \quad\mbox{and}\quad
\log\left(\frac{1}{\alpha^2}\right)=n_2\Lambda
$$  
with $\Lambda =-\log\alpha$, $n_1=1$ and $n_2=2$.
By following the lines of the proof of Theorem~\ref{Thm1} and in particular~(\ref{etaRel}) we can explicitely get the value of $\eta$. In fact, since the roots of the equation $1-z-z^2=0$ are given by 
$z_1=\frac{\sqrt{5}-1}{2}=\alpha=e^{-\Lambda }$ and $z_2=\frac{-\sqrt{5}-1}{2}$,
it follows that $d=0$ and
$$\eta=1+\frac{\log|z_2|}{\Lambda }= 1+\frac{\log\left|\frac{-\sqrt{5}-1}{2}\right|}{-\log\left(\frac{\sqrt{5}-1}{2}\right)}=2.$$ 
According to Theorem \ref{disRatCase}, this shows that the discrepancy is of the order $\mathcal{O}\left({1}/{k(n)}\right)$ and therefore it is optimal. 

In the general case set $m=L+S$. Of course we are in the rational case
since $p_i = \alpha$ or $p_i = \alpha^2$ for $i=1,\ldots,m$. More precisely, according to Definition \ref{defRR} we have
$\Lambda  = \log(1/\alpha)$ and $n_i \in \{1,2\}$ corresponding to $p_i = \alpha^{n_i}$.
The zeros of the equation  $$1-Lz-Sz^2=0$$ are given by $z_1=\frac{-L+\sqrt{L^2+4S}}{2S}=\alpha$ and $z_2=\frac{-L-\sqrt{L^2+4S}}{2S}$. 
Hence, 
$$
\eta=1+\frac{\log\left|\frac{-L-\sqrt{L^2+4S}}{2S}\right|}{\Lambda}=
1+\frac{\log\left(\frac{L+\sqrt{L^2+4S}}{2S}\right)}{\Lambda}.
$$ 
Consequently, we have $\eta<1$ if and only if $\frac{L+\sqrt{L^2+4S}}{2S}<1$
or if $S>L+1$. Similarly we have $\eta = 1$ if and only if $S = L+1$ and
$\eta > 1$ if and only if $S< L+1$. This is in perfect accordance with
the results of Carbone \cite{Car}. The discrepancy bounds are (of course)
also of the same kind.

\subsection{Sequences related to Pisot numbers}
Let us introduce a class of sequences having optimal discrepancy. This kind of sequences is related to Pisot numbers, so let us recall some classical definitions.
\begin{definition}\ 
\begin{itemize}
\item An \emph{algebraic number} is a number which is a root of a non-zero polynomial in one variable with rational (or equivalently, integer) coefficients.

\item Given an algebraic number $\alpha$, there is a unique monic polynomial with rational coefficients of least degree that has $\alpha$ as a root. This polynomial is called \emph{minimal polynomial} of $\alpha$. 

\item If an algebraic number $\alpha$ has its minimal polynomial of degree $n$, then the algebraic number is said to be of degree $n$. 

\item The \emph{conjugates} of an algebraic number $\alpha$ are the other roots of its minimal polynomial.

\item An \emph{algebraic integer} is an algebraic number which is a root of a monic polynomial with integer coefficients.
\end{itemize}
\end{definition}
\begin{definition}\ \\
A Pisot number $\beta$ is an algebraic integer larger than $1$ 
with the property that all its conjugates have modulus smaller than $1$. 
\end{definition}
A prominent example of Pisot numbers are the real roots of a polynomial of the form
\begin{equation}\label{eqPisot}
z^k - a_1 z^{k-1} - a_2 z^{k-2} - \cdots - a_k = 0,
\end{equation}
where $a_j$ are positive integers with $a_1\ge a_2 \ge \cdots \ge a_k$ (see \cite{Pisot}). In this case the polynomial in (\ref{eqPisot}) is
also irreducible over the rationals.

Suppose now that $\rho$ is a partition of $m = a_1+a_2 + \cdots +a_k$ 
intervals, where $a_j$ intervals have length $\alpha^j$, $1\le j \le k$,
$\alpha = 1/\beta$, and $\beta$ is the Pisot number
related to the polynomial (\ref{eqPisot}).  Note that we have
\[
a_1 \alpha + a_2 \alpha^2 + \cdots + a_k \alpha^k = 1.
\]
Since all conjugates of $\alpha$ have now modulus larger than $1$
it follows that $\eta > 1$. This means that the order of
magnitude of the discrepancy is optimal, namely $1/k(n)$.
Moreover, we can note that $LS$-sequences are a special instance for $k=2$, $a_1 = L$ and $a_2 = S$ with $L\ge S$.

Note that in the Pisot case all complex zeros of the polynomial are simple, since
the polynomial is irreducible over the rationals. However, this is not necessarily
true in less restrictive cases than Pisot numbers. For example, let $\alpha = 1/5$
and consider one interval of length $\alpha = 1/5$, $16$ intervals of lengths $\alpha^2 = 1/25$ and $20$ intervals of lengths $\alpha^3 = 1/125$. 
Since $\alpha + 16 \alpha^2 + 20 \alpha^3 = 1$ we have a proper
partition $\rho$. Here the roots of the polynomial
$z + 16 z^2 + 20 z^3 = 1$ are $z_1=\alpha = 1/5$ and 
$z_2 = z_3 = - 1/2$ (which is a double root). Hence, we obtain
$\eta = 1 - (\log 2)/(\log 5) = 0.56932\ldots < 1 $ and $d=1$. Consequently the 
discrepancy is bounded by 
\[
D_n = \mathcal{O}( (\log k(n))\,  k(n)^{-\eta}),
\]
and this upper bound is optimal.

\subsection{The rational case on fractals}
The same procedure of $\rho$-refinements can be used to obtain an extension of the results introduced in Section \ref{UDFractals} to fractals defined by similarities which do not have the same ratio and which satisfy the OSC. In fact, we will describe an analogous of the method of successive $\rho-$refinements which allows to produce sequences of partitions on this new class of fractals. Actually, we will introduce a new correspondence between nodes of the tree associated to Khodak's algorithm and the subsets belonging to the partitions generated on the fractal.\\

Let $\varphi=\{\varphi_1,\ldots,\varphi_m\}$ be a system of $m$ similarities on $\mathbb{R}^h$ which have ratios $c_1,\ldots,c_m\in]0,1[$ respectively and which verify the OSC. Let $F$ be the attractor of $\varphi$ and let $s$ be its Hausdorff dimension. Moreover, we will consider the normalized $s$-dimensional Hausdorff measure $P$  on the fractal $F$, that is given by~(\ref{ProbFrac}).

Start with a tree having a root node of probability 1, which corresponds to the fractal $F$, and $m$ leaves corresponding to the $m$ imagines of $F$ through the $m$ similarities, i.e. $\varphi_1(F),\ldots,\varphi_m(F)$. The probability of each node is given by the probability of the corresponding subset, that is $p_i=P(\varphi_i(F))=c_i^s$. At each iteration we select the leaves having the highest probability and grow $m$ children out of each of them. On the fractal this corresponds to apply successively the $m$ similarities only to those subsets of $F$ having the highest probability at this certain step.
By iterating this procedure we obtain a tree associated to the sequence of partitions on the fractal $F$, which is the same tree generated by Khodak's algorithm. 

Let us denote by $(\pi_n)$ the sequence of partitions of $F$ constructed by this technique, i.e. 
$$
\pi_n=\big\{\psi_{j_{k(n)}}\psi_{j_{({k(n)}-1)}}\cdots \ \psi_{j_1}(F) : {j_1},\ldots,{j_{k(n)}}\in\{1,\ldots,m\}\big\}.
$$
where $k(n)$ is the number of sets constructed at the step $n$.

Let us denote by $\mathscr{E}_n$ the collection of the $k(n)$ sets $E_i^n$ belonging to the partition~$\pi_n$ and by~$\mathscr{E}$ the union of the families $\mathscr{E}_n$ $\forall n\in\mathbb{N}$. The sets of the class $\mathscr{E}$ are called \emph{elementary sets} because they are exactly constructed as the ones defined in Subsection~\ref{subsec: Algo}.

Lemma \ref{det-elementary} and Lemma \ref{ElemContSets} guarantee that the class $\mathscr{E}$ is determining and consisting of $P$-continuity sets. Now, if we choose a point $t_i^{(n)}$ in each $E_i^n\in\pi_n$, we can consider the elementary discrepancy of this set of points on the fractal, i.e.
$$D_n^{\mathscr{ E}}=\sup_{E\in\mathscr{ E}}\Bigg|\frac{1}{k(n)}\sum_{i=1}^{k(n)}\chi_E\left(t_i^{(n)}\right)-P(E)\Bigg|.
$$ 
By using a procedure similar to the one used in the proof of the 
Theorem~\ref{disRatCase} we can prove the following theorem.
\begin{theorem}\label{Thm: RatCaseDiscr}\ \\
Let $(\pi_i)$ be the sequence of partitions of $F$ just constructed. Assume that $\log\left(\frac{1}{p_1}\right),\ldots$ $\ldots,\log\left(\frac{1}{p_m}\right)$ are rationally related. Then we have the following bounds for the elementary discrepancy
\begin{equation}\label{uppBound-fractal}
D_n^{\mathscr{ E}}= \left\{ \begin{array} {ll} 
\mathcal{O}\left( (\log k(n))^d k(n)^{-\eta} \right) & \mbox{if $0< \eta \le 1$,}\\
\mathcal{O}\left( k(n)^{-1} \right) & \mbox{if $\eta > 1$.}
\end{array} \right.
\end{equation}
Furthermore, both upper bounds are best possible.
\end{theorem}
\proof\ \\
Fix a step in the algorithm corresponding to a certain parameter $r$ of the form $r=e^{-n\Lambda }$ for some integer $n\ge 0$. First we observe that the number $N_E^{(n)}$ of elementary sets if $\mathscr{E}_n$ that are contained in a fixed elementary set $E$
is given by $M_{r/P(E)}$. This fact implies that
\begin{equation}\label{eqNE}
N_E^{(n)} = \frac{c'}r P(E) + \mathcal{O}\left( |\log r|^d \,r^{-1+\eta} P(E)^{1-\eta}\right).
\end{equation}
This proves (\ref{uppBound-fractal}) directly for $\eta \le 1$ and also shows that
this bound is optimal.

If $\eta > 1$ then we argue recursively. The elementary set $E$ is
either contained in $\mathscr{E}_1  
= \{ \varphi_1(F), \ldots, \varphi_m(F)\}$, 
which means that we can use (\ref{eqNE}) for
$P(E) \in \{p_1,\ldots, p_m\}$, or it is part of 
$E_j=\varphi_j(F)$ for some $j$. In the latter case we can rewrite 	
$N_E^{(n)} - k(n) P(E)$ to
\[
N_E^{(n)} - k(n) P(E) = \left( N_E^{(n)} - k(n-1) \frac{P(E)}{P(E_j)}\right) +
\left( k(n-1)\frac{P(E)}{P(E_j)} - k(n)P(E) \right),  
\]
which leads to a recurrence of the form
\[
\Delta_n^{\mathscr{ E}} = 
\sup_{E \in \mathscr{ E}} \left| N_E^{(n)} - k(n) P(E) \right|
\le \Delta_{n-1}^{\mathscr{ E}} + \mathcal{O}\left( n^d e^{\Lambda n (1-\eta)} \right).
\]
Here we have
\[
\Delta_n^{\mathscr{ E}} = \mathcal{O}\left( \sum_{k\le n} k^d e^{-\Lambda k(\eta-1)} \right) 
= \mathcal{O}\left( 1 \right)
\]
and consequently $D_n^{\mathscr{ E}} = \mathcal{O}(1/k(n))$ (which is also optimal).\\
\endproof

In particular it follows that the sequence of partitions $(\pi_n)$ is u.d.\! with respect to $P$. Actually, this remains true in the irrationally related case but we can derive effective upper bounds for the discrepancy only in very specific cases.

%% file: Concl.tex
\chapter{Conclusions and open problems}
\markboth{4.CONCLUSIONS AND OPEN PROBLEMS}{4.CONCLUSIONS AND OPEN PROBLEMS}
The main objectives of this work were to construct new classes of u.d.\! sequences of partitions and of points on fractals and on $[0,1]$ and to study their discrepancy. In fact with reference to fractals, we were able to find a general algorithm for producing sequences of van der Corput type on fractals generated by an IFS consisting of similarities which have the same ratio and which satisfy the OSC. Secondly, we got bounds for the discrepancy of a class of generalized Kakutani's sequences of partitions on $[0,1]$ constructed through the recent technique of successive $\rho-$refinements. Moreover, these last results enabled us to introduce a new family of u.d.\! sequences on a wider class of fractals. Although the techniques we used are direct and explicit procedures, the quantitive analysis of the distribution properties of the sequences constructed was proved to be more involved than we expected. Indeed, we now want to discuss major problems as far as the study of discrepancy is concerned. 
\\

A first problem was finding a unifying approach to the discrepancy on the class of fractals considered in this work. In fact, the only kind of discrepancy which makes sense for all the fractals generated by IFS and satisfying the OSC appears to be the so-called elementary discrepancy, as we have already noted in Subsection \ref{subsec: ElemDiscr}. We got estimates for the elementary discrepancy of the sequences generated by our techniques in Theorem \ref{Thm: ElemDiscr} and Theorem \ref{Thm: RatCaseDiscr}. In both cases the particular properties given by our procedures to the sequences allowed to obtain estimates for their discrepancy with a quite direct proof. Neverthless, it is not surprising that these properties are no more sufficient to provide estimates for another kind of discrepancy. This is due to the fact that the elementary sets are intrisically related to the construction of the fractals belonging to our class and not to the specific geometry of each fractal. In literature, there are few papers devoted to u.d.\! sequences on fractals and to estimates of their discrepancy (see \cite{CT,CT2,Grab-Tichy}). In these articles, the various types of discrepancy considered depend heavily on the geometric features of the particular fractal considered. At the moment the choice of elementary discrepancy seems to be the only one which allows to overcome the problem of the peculiar complexity of each fractal and to give explicit results for the fractals of the whole class taken in consideration. So a still open problem is finding a definition of discrepancy, different from the elementary one, which is general at least for the fractals of our class. An attempt in this direction has been proposed by Albrecher, Matou\v sek and Tichy in \cite{AMT}, but it concerns the average discrepancy.  
\\

This work also leaves open several questions concerning the discrepancy of the generalized Kakutani's sequences constructed by the $\rho-$refinements technique on $[0,1]$. As we have already said in Section \ref{sec:Discr}, the natural problem of studying the behaviour of the discrepancy $D_n$ (\ref{discrPart}) of the sequence of partitions constructed as $n$ tends to infinity was posed in \cite{Vol}. In this thesis we partially answered to that question. Indeed, we were able to find a class of partitions $\rho$ such that the speed of convergence to zero of the discrepancy of the sequences generated by successive $\rho-$refinements is quite high (see Theorem \ref{disRatCase}). The strategy used to get bounds for the discrepancy of these sequences exploits the correspondence between the procedure of $\rho-$refinements and Khodak's algorithm. Although this analogy is crucial in our analysis, it is not sufficient to give effective estimates in the irrational case, too. In fact, as we have already discussed in Subsection \ref{subsec: Irr}, the irrational case is too closely related to the Diophantine approximation properties of the quotients $\frac{\log p_i}{\log p_j}$ where the $p_i$'s are the lenghts of the starting partition $\rho$. This problem is mainly caused by the necessity of having an explicit error term in the asymptotic expansion (\ref{irratCase}), which is not provided by the results related to Khodak's algorithm in \cite{Drmota}. In particular, we gave more precise information on the error term when the initial probabilities are only two $p, q=1-p$ and $\gamma = \log p/\log q$ is badly approximable or when $p$ and $q$ are algebraic numbers. In these two cases, the correspondent upper bounds for the discrepancy are weaker than the ones in the rational case. Consequently, even if the continued fractional expansion of $\gamma$ could be useful for getting more explicit upper bounds for discrepancy, the relative estimates will probably be again poor. Moreover, as we mentioned in Subsection \ref{subsec: Irr} the paper \cite{Flajolet-Vallee} shows the difficulties presented by the irrational case when $\rho$ consists of more than two intervals. Therefore it could be interesting to develop different approaches to the problem, which allow to improve the results at least for a class of instances of the irrational case.\\

A further interesting problem is finding explicit algorithms to provide low discrepancy sequences of points associated to a low discrepancy u.d.\! sequence of partitions constructed by successive $\rho-$refinements. In Subsection \ref{subsec: AssSeqPoints} we analyzed this problem in general and we reported the main result about this question developed in \cite{Vol}. In fact, in this paper it has been proved that a random reordering of the points of a u.d.\! sequence of partitions gives with probability one a u.d.\! sequence of points (see Theorem \ref{RandReord}). Neverthless, this important theoretical result does not give any information about the discrepancy of the sequences of points obtained by the random reordering. Some results in this direction have been already obtained concerning to the LS-sequences. In fact, in \cite{Car} the author presents an explicit procedure for associating to LS-sequences with low discrepancy sequences of points having low discrepancy, too. It would be of great interest trying to extend this algorithm to the whole class of sequences of partitions constructed by successive $\rho-$refinements included in the rational case.